\documentclass[12pt,leqno]{article}
\usepackage{amssymb,amsmath,amsfonts,esint}
\usepackage{relsize}
\newtheorem{theorem}{Theorem}[section]
\newtheorem{definition}[theorem]{Definition}
\newtheorem{corollary}[theorem]{Corollary}
\newtheorem{proposition}[theorem]{Proposition}
\newtheorem{remark}[theorem]{Remark}
\newtheorem{lemma}[theorem]{Lemma}
\newtheorem{problem}[theorem]{Problem}
\newtheorem{question}[theorem]{Question}
\newcommand {\Kc}      {{\mathcal K}}
\newcommand {\Lc}      {{\mathcal L}}
\newcommand {\Mc}      {{\mathcal M}}
\newcommand {\Hc}      {{\mathcal H}}
\newcommand {\Ac}      {{\mathcal A}}
\newcommand {\Bc}      {{\mathcal B}}
\newcommand {\Qc}      {{\mathcal Q}}
\newcommand {\Ic}      {{\mathcal I}}
\newcommand {\Sc}      {{\mathcal S}}
\newcommand {\Gc}      {{\mathcal G}}
\newcommand {\R}       {{\bf R}}
\newcommand {\N}       {{\bf N}}
\newcommand {\mZ}      {{\mathbb Z}}
\newcommand {\tQ}      {\widetilde{Q}}
\newcommand {\tH}      {\widetilde{H}}
\newcommand {\tI}      {\widetilde{I}}
\newcommand {\tB}      {\widetilde{B}}
\newcommand {\tK}      {\widetilde{K}}
\newcommand {\tU}      {\widetilde{U}}
\newcommand {\tL}      {\widetilde{{\mathcal L}}}
\newcommand {\tAc}     {\widetilde{{\mathcal A}}}
\newcommand {\tQc}     {\widetilde{{\mathcal Q}}}
\newcommand {\tf}      {\tilde{f}}
\newcommand {\tx}      {\tilde{x}}
\newcommand {\hK}      {\widehat{K}}
\newcommand {\RN}      {\R^n}
\newcommand {\RT}      {\R^2}
\newcommand {\ve}      {\varepsilon}
\newcommand {\LOP}     {L^1_p(\RN)}
\newcommand {\LPM}     {L_p(\RN;\mu)}
\newcommand {\LMP}     {L^m_p(\RN)}
\newcommand {\SUM}     {L^1_p(\RN)+L_p(\RN;\mu)}
\newcommand {\LPRN}    {L_p(\RN)}
\newcommand {\intl}    {\int\limits}
\newcommand {\av}[1]   {f_{#1}}     % average symbol
\newcommand {\emp}     {\emptyset}
\newcommand {\row}     {\rho_w}
\newcommand {\smed}    {\mathlarger{\sum}}
\newcommand {\sbig}    {\mathlarger{\mathlarger{\sum}}}
\newcommand {\shuge}    {\mathlarger{\mathlarger{\mathlarger{\sum}}}}
\newcommand {\ssmall}  {\mathsmaller{\sum}}
\newcommand {\q}       {90}
\newcommand {\LE}      {\Lc_E}
\newcommand {\hL}      {\hat{{\mathcal L}}}
\newcommand {\QL}      {Q^{(L)}}
\newcommand {\PRL}     {\mathcal{PR}}
\newcommand {\lr}      {\leftrightarrow}
\newcommand {\card}    {\#\,}
\newcommand {\diam}    {\operatorname{diam}}
\newcommand {\dist}    {\operatorname{dist}}
\newcommand {\supp}    {\operatorname{supp}}
\newcommand {\bx}      {\hfill$\blacksquare$}
\newcommand {\rbx}     {\hfill$\vartriangleleft$}
\newcommand {\nn}      {\nonumber}
\newcommand {\rf}[1]    {(\ref{#1})}      %no references
\newcommand {\reff}[1] {\ref{#1}}         %no references
\newcommand{\lbl}[1]      {\label{#1}}       %no ref
\newcommand{\be}          {\begin{eqnarray}}
\newcommand{\bel}[1]      {\begin{eqnarray} \label{#1}}
\newcommand{\ee}           {\end{eqnarray}}
\newcommand {\SECT}[2] {\section*{\centerline{\normalsize
{\bf #1}}} \setcounter{section}{#2}
\setcounter{theorem}{0}\setcounter{equation}{0}}
\begin{document}
%----------------------------------------------------------
\medskip
%@@@@@@@@@@@@@@@@@@@@@@@@@@@@@@@@@@@@@@@@@@@@@@@@@@@@@@@@@@
%@@@@@@@@@@@@@@@@@@@@@@@@@@@@@@@@@@@@@@@@@@@@@@@@@@@@@@@@@@
%----------------------------------------------------------
\centerline{{\bf ON THE SUM OF A SOBOLEV SPACE AND A WEIGHTED $L_p$-SPACE}} \vspace*{5mm}
%----------------------------------------------------------
\vspace*{8mm}
%----------------------------------------------------------
\centerline{By~ {\it Pavel Shvartsman}} \vspace*{5 mm}
%----------------------------------------------------------
\centerline {\it Department of Mathematics, Technion - Israel Institute of Technology}\vspace*{2 mm}
%----------------------------------------------------------
\centerline{\it 32000 Haifa, Israel}\vspace*{2 mm}
%----------------------------------------------------------
\centerline{\it e-mail: pshv@tx.technion.ac.il}
%----------------------------------------------------------
\vspace*{10 mm}
%----------------------------------------------------------
\renewcommand{\thefootnote}{ }
%----------------------------------------------------------
\footnotetext[1]{{\it\hspace{-6mm}Math Subject
Classification} 46E35}
%@@@@@@@@@@@@@@@@@@@@@@@@@@@@@@@@@@@@@@@@@@@@@@@@@@@@@@@@@@
%@@@@@@@@@@@@@@@@@@@@@@@@@@@@@@@@@@@@@@@@@@@@@@@@@@@@@@@@@@
%----------------------------------------------------------
\begin{abstract} Let $p>n$ and let $\LOP$ be a homogeneous Sobolev space. For an arbitrary Borel measure $\mu$ on $\RN$ we give a constructive characterization of the space $$\Sigma=\LOP+\LPM.$$ We express the norm in this space in terms of certain oscillations with respect to the measure $\mu$. This enables us to describe the $K$-functional for the couple $(\LPM,\LOP)$ in terms of these oscillations, and to prove that this couple is quasi-linearizable.
\end{abstract}
%----------------------------------------------------------
\renewcommand{\contentsname}{ }
\tableofcontents
%----------------------------------------------------------
\addtocontents{toc}{{\centerline{\sc{Contents}}}
\vspace*{10mm}\par}
%----------------------------------------------------------
%@@@@@@@@@@@@@@@@@@@@@@@@@@@@@@@@@@@@@@@@@@@@@@@@@@@@@@@@@@
%@@@@@@@@@@@@@@@@@@@@@@@@@@@@@@@@@@@@@@@@@@@@@@@@@@@@@@@@@@
%@@@@@@@@@@@@@@@@@@@@@@@@@@@@@@@@@@@@@@@@@@@@@@@@@@@@@@@@@@
%----------------------------------------------------------
%@@@@@@@@@@@@@@@@@@@@@@@@@@@@@@@@@@@@@@@@@@@@@@@@@@@@@@@@@@
%@@@@@@@@@@@@@@@@@@@@@@@@@@@@@@@@@@@@@@@@@@@@@@@@@@@@@@@@@@
%@@@@@@@@@@@@@@@@@@@@@@@@@@@@@@@@@@@@@@@@@@@@@@@@@@@@@@@@@@
%@@@@@@@@@@@@@@@@@@@@@@@@@      @@@@@@@@@@@@@@@@@@@@@@@@@@@
%@@@@@@@@@@@@@@@@@@@@@@@          @@@@@@@@@@@@@@@@@@@@@@@@@
%@@@@@@@@@@@@@@@@@@@@@              @@@@@@@@@@@@@@@@@@@@@@@
%@@@@@@@@@@@@@@@@@@@     SECTION 1    @@@@@@@@@@@@@@@@@@@@@
%@@@@@@@@@@@@@@@@@@@@@              @@@@@@@@@@@@@@@@@@@@@@@
%@@@@@@@@@@@@@@@@@@@@@@@          @@@@@@@@@@@@@@@@@@@@@@@@@
%@@@@@@@@@@@@@@@@@@@@@@@@@      @@@@@@@@@@@@@@@@@@@@@@@@@@@
%@@@@@@@@@@@@@@@@@@@@@@@@@@@@@@@@@@@@@@@@@@@@@@@@@@@@@@@@@@
%@@@@@@@@@@@@@@@@@@@@@@@@@@@@@@@@@@@@@@@@@@@@@@@@@@@@@@@@@@
%----------------------------------------------------------
%@@@@@@@@@@@@@@@@@@@@@@@@@@@@@@@@@@@@@@@@@@@@@@@@@@@@@@@@@@
%----------------------------------------------------------
\SECT{1. Introduction.}{1}
%----------------------------------------------------------
\addtocontents{toc}{~~~~1. Introduction. \hfill \thepage\\\par}
%----------------------------------------------------------
\indent
%@@@@@@@@@@@@@@@@@@@@@@@@@@@@@@@@@@@@@@@@@@@@@@@@@@@@@@@@@@
%----------------------------------------------------------
\par Let $\mu $ be a non-trivial non-negative Borel measure on ${\bf R}^{n}$ and let $\LPM$, $1\le p\le\infty,$ be the space $L_p$ on $\RN$ with respect to the measure $\mu$, with the standard norm
%----------------------------------------------------------
$$ \|f\|_{\LPM}=\left(\,\,\intl\limits_{\RN}|f|^{p}d\mu\right)^
{\frac{1}{p}}.
$$
%----------------------------------------------------------
%@@@@@@@@@@@@@@@@@@@@@@@@@@@@@@@@@@@@@@@@@@@@@@@@@@@@@@@@@@
%@@@@@@@@@@@@@@@@@@@@@@@@@@@@@@@@@@@@@@@@@@@@@@@@@@@@@@@@@@
%----------------------------------------------------------
\par By $\LOP$ we denote the homogeneous Sobolev space
consisting of all (equivalence classes of) real valued
functions $f\in L_{p,loc}(\RN)$ whose distributional partial derivatives of the first order belong to the space $\LPRN$. We equip the space $\LOP$ with the seminorm
%----------------------------------------------------------
$$
\|f\|_{\LOP}:=\|\nabla f\|_{\LPRN}.
$$
%----------------------------------------------------------
%@@@@@@@@@@@@@@@@@@@@@@@@@@@@@@@@@@@@@@@@@@@@@@@@@@@@@@@@@@
%@@@@@@@@@@@@@@@@@@@@@@@@@@@@@@@@@@@@@@@@@@@@@@@@@@@@@@@@@@
%@@@@@@@@@@@@@@@@@@@@@@@@@@@@@@@@@@@@@@@@@@@@@@@@@@@@@@@@@@
%@@@@@@@@@@@@@@@@@@@@@@@@@@@@@@@@@@@@@@@@@@@@@@@@@@@@@@@@@@
%----------------------------------------------------------
\par This paper is devoted to the following main
%@@@@@@@@@@@@@@@@@@@@@@@@@@@@@@@@@@@@@@@@@@@@@@@@@@@@@@@@@@
%@@@@@@@@@@@@@@@@@@@@@@@@@@@@@@@@@@@@@@@@@@@@@@@@@@@@@@@@@@
%----------------------------------------------------------
\begin{problem}\lbl{MAIN-PR1} Given a function $f\in L_{p,loc}(\RN;\mu)$, how can we tell whether $f$ belongs to $\SUM$, i.e., whether there exist functions $f_1\in\LOP$ and $f_2\in\LPRN$ such that $f=f_1+f_2~?$
%----------------------------------------------------------
\end{problem}
%----------------------------------------------------------
%@@@@@@@@@@@@@@@@@@@@@@@@@@@@@@@@@@@@@@@@@@@@@@@@@@@@@@@@@@
\par We also consider a quantitative version of Problem \reff{MAIN-PR1} related to calculation of the norm of $f$ in the space
%----------------------------------------------------------
$$
\mathsmaller{\sum}:=\LOP+\LPM.
$$
%----------------------------------------------------------
As usual, the space $\sum$ is normed by
%----------------------------------------------------------
$$
\|f\|_{\ssmall}:=\inf\{\|f_1\|_{\LOP}+\|f_2\|_{\LPM}:
f_1+f_2=f, f_1 \in \LOP ,f_2\in\LPM\}.
$$
%----------------------------------------------------------
%@@@@@@@@@@@@@@@@@@@@@@@@@@@@@@@@@@@@@@@@@@@@@@@@@@@@@@@@@@
%@@@@@@@@@@@@@@@@@@@@@@@@@@@@@@@@@@@@@@@@@@@@@@@@@@@@@@@@@@
%----------------------------------------------------------
\begin{problem}\lbl{MAIN-PR2} What is the order of magnitude of the norm of a function $f$ in the space $\sum=\LOP+\LPM$~?
%----------------------------------------------------------
\end{problem}
%----------------------------------------------------------
%@@@@@@@@@@@@@@@@@@@@@@@@@@@@@@@@@@@@@@@@@@@@@@@@@@@@@@@@@@
%@@@@@@@@@@@@@@@@@@@@@@@@@@@@@@@@@@@@@@@@@@@@@@@@@@@@@@@@@@
%@@@@@@@@@@@@@@@@@@@@@@@@@@@@@@@@@@@@@@@@@@@@@@@@@@@@@@@@@@
%----------------------------------------------------------
\par In this paper we solve Problems \reff{MAIN-PR1} and \reff{MAIN-PR2} by presenting a constructive formula for calculation of the order of magnitude of the norm in the space $\SUM$. This formula is  expressed in terms of certain local oscillations of functions with respect to the measure $\mu$.
%----------------------------------------------------------
\par Before we formulate the main result of the paper
we need to define several notions and fix some notation:
%----------------------------------------------------------
\par Throughout this paper, the word ``cube'' will mean a closed cube in ${\bf R}^{n}$ whose sides are parallel to the coordinate axes. We let $Q(x,r)$ denote the cube in $\RN$ centered at $x$ with side length $2r$. Given $\alpha >0$ and a cube $Q$ we let $\alpha Q$ denote the dilation of $Q$ with respect to its center by a factor of $\alpha $. (Thus $\alpha\,Q(x,r)=Q(x,\alpha r)$.) The Lebesgue measure of a measurable set $A\subset \RN$ will be denoted by $\left|A\right|$.
%----------------------------------------------------------
%@@@@@@@@@@@@@@@@@@@@@@@@@@@@@@@@@@@@@@@@@@@@@@@@@@@@@@@@@@
%@@@@@@@@@@@@@@@@@@@@@@@@@@@@@@@@@@@@@@@@@@@@@@@@@@@@@@@@@@
%----------------------------------------------------------
\par Here now is the main result of our paper:
%@@@@@@@@@@@@@@@@@@@@@@@@@@@@@@@@@@@@@@@@@@@@@@@@@@@@@@@@@@
%@@@@@@@@@@@@@@@@@@@@@@@@@@@@@@@@@@@@@@@@@@@@@@@@@@@@@@@@@@
%@@@@@@@@@@@@@@@@@@@@@@@@@@@@@@@@@@@@@@@@@@@@@@@@@@@@@@@@@@
%@@@@@@@@@@@@@@@@@@@@@@@@@@@@@@@@@@@@@@@@@@@@@@@@@@@@@@@@@@
%@@@@@@@@@@@@@@@@@@@@@@@@@@@@@@@@@@@@@@@@@@@@@@@@@@@@@@@@@@
%@@@@@@@@@@@@@@@@@@@@@@@@@@@@@@@@@@@@@@@@@@@@@@@@@@@@@@@@@@
%----------------------------------------------------------
\begin{theorem}\lbl{MAIN-CR} Let $n<p<\infty$ and let $\mu$ be a non-trivial non-negative Borel measure on $\RN$.  A function $f\in L_{p,loc}(\RN;\mu)$ belongs to the space $\LOP+\LPM$ if and only if there exists a positive constant $\lambda$ which satisfies the following conditions for a certain absolute positive constant $\gamma$:
%----------------------------------------------------------
\par Let $\Qc$ be an arbitrary finite family of pairwise disjoint cubes in $\RN$. Suppose that to each cube $Q\in\Qc$ we have arbitrarily assigned two cubes $Q',Q''\in\Qc$ such that
%----------------------------------------------------------
\bel{IM-QP}
Q'\cup Q''\subset \gamma Q.
\ee
%----------------------------------------------------------
\par Then the following inequality
%----------------------------------------------------------
\bel{CR}
\shuge_{Q\in\Qc}\,\,
\frac{(\diam Q)^{n-p}\iint \limits_{Q'\times Q''}
|f(x)-f(y)|^p\, d\mu(x)d\mu(y)}
{ \{(\diam Q')^{n-p}+\mu(Q')\} \{(\diam Q'')^{n-p}+\mu(Q'')\}}
\le \lambda
\ee
%----------------------------------------------------------
holds. Furthermore,
%----------------------------------------------------------
$$
\|f\|_{\sum}\sim \inf \lambda^{\frac{1}{p}}
$$
%----------------------------------------------------------
with constants of equivalence depending only on $n$ and $p$.
%----------------------------------------------------------
\end{theorem}
%----------------------------------------------------------
%@@@@@@@@@@@@@@@@@@@@@@@@@@@@@@@@@@@@@@@@@@@@@@@@@@@@@@@@@@
%@@@@@@@@@@@@@@@@@@@@@@@@@@@@@@@@@@@@@@@@@@@@@@@@@@@@@@@@@@
%@@@@@@@@@@@@@@@@@@@@@@@@@@@@@@@@@@@@@@@@@@@@@@@@@@@@@@@@@@
%@@@@@@@@@@@@@@@@@@@@@@@@@@@@@@@@@@@@@@@@@@@@@@@@@@@@@@@@@@
%@@@@@@@@@@@@@@@@@@@@@@@@@@@@@@@@@@@@@@@@@@@@@@@@@@@@@@@@@@
%----------------------------------------------------------
%@@@@@@@@@@@@@@@@@@@@@@@@@@@@@@@@@@@@@@@@@@@@@@@@@@@@@@@@@@
%@@@@@@@@@@@@@@@@@@@@@@@@@@@@@@@@@@@@@@@@@@@@@@@@@@@@@@@@@@
%@@@@@@@@@@@@@@@@@@@@@@@@@@@@@@@@@@@@@@@@@@@@@@@@@@@@@@@@@@
%@@@@@@@@@@@@@@@@@@@@@@@@@@@@@@@@@@@@@@@@@@@@@@@@@@@@@@@@@@
%@@@@@@@@@@@@@@@@@@@@@@@@@@@@@@@@@@@@@@@@@@@@@@@@@@@@@@@@@@
%----------------------------------------------------------
\begin{remark}\lbl{KF} {\em The topic under consideration can be referred to as the Real Interpolation Method for the Banach couple $\vec{A}=(\LPM,\LOP)$, or, more specifically, as the calculation of the $K$-functional
%----------------------------------------------------------
$$
K(t; f :\vec{A}):=\inf\{\|f_1\|_{\LPM}+t\|f_2\|_{\LOP}:
f_1+f_2=f, f_1\in\LPM,f_2 \in \LOP\}.
$$
%----------------------------------------------------------
Here $t$ is a positive number. (See, e.g. \cite{BL}.)
%----------------------------------------------------------
\par Thus $\|f\|_{\sum}=K(1; f :\vec{A})$ and
%----------------------------------------------------------
$$
K(t; f :\vec{A})= t\|f\|_{\Sigma_t}~~~\text{where}~~~
\Sigma_t:=\LOP+L_p(\RN;\tfrac{1}{t^p}\mu).
$$
%----------------------------------------------------------
\par We recall the classical result of Peetre \cite{P} (see also \cite{BSh}, p. 339), which states that whenever $1\le p\le\infty$ and $\mu$ {\it is Lebesgue measure on $\RN$},
%----------------------------------------------------------
$$
K(t;f:(L_p(\RN),\LOP))\sim\omega_1(t,f)_{L_p(\RN)}
$$
%----------------------------------------------------------
with constants depending only on $p$ and $n$. Here
%----------------------------------------------------------
$$
\omega_1(t,f)_{L_p(\RN)}=
\sup_{\|h\|\le t}\left\{\,\,\intl_{\RN}
|f(x+h)-f(x)|^p\,dx\right\}^p
$$
%----------------------------------------------------------
is the modulus of smoothness of $f$ in $\LPRN$.
%----------------------------------------------------------
\par This result leads us to a solution of Problem \reff{MAIN-PR2} for the particular case where $\mu$ is Lebesgue measure multiplied by an arbitrary positive parameter $s$. In this case
%----------------------------------------------------------
$$
\|f\|_{\sum}\sim s^{\frac1p}\,\omega_1
\left(s^{-\frac{1}{p}},f\right)_{L_p(\RN)}
$$
%----------------------------------------------------------
with constants of equivalence depending only on $n$. To the best of our knowledge, this measure $\mu$ is the only example of a measure for which a constructive criterion for the norm of a function in the sum $\SUM$ is known so far.
\rbx}
%----------------------------------------------------------
\end{remark}
%----------------------------------------------------------
%@@@@@@@@@@@@@@@@@@@@@@@@@@@@@@@@@@@@@@@@@@@@@@@@@@@@@@@@@@
%@@@@@@@@@@@@@@@@@@@@@@@@@@@@@@@@@@@@@@@@@@@@@@@@@@@@@@@@@@
%@@@@@@@@@@@@@@@@@@@@@@@@@@@@@@@@@@@@@@@@@@@@@@@@@@@@@@@@@@
%@@@@@@@@@@@@@@@@@@@@@@@@@@@@@@@@@@@@@@@@@@@@@@@@@@@@@@@@@@
%@@@@@@@@@@@@@@@@@@@@@@@@@@@@@@@@@@@@@@@@@@@@@@@@@@@@@@@@@@
%@@@@@@@@@@@@@@@@@@@@@@@@@@@@@@@@@@@@@@@@@@@@@@@@@@@@@@@@@@
%----------------------------------------------------------
%@@@@@@@@@@@@@@@@@@@@@@@@@@@@@@@@@@@@@@@@@@@@@@@@@@@@@@@@@@
%@@@@@@@@@@@@@@@@@@@@@@@@@@@@@@@@@@@@@@@@@@@@@@@@@@@@@@@@@@
%@@@@@@@@@@@@@@@@@@@@@@@@@@@@@@@@@@@@@@@@@@@@@@@@@@@@@@@@@@
%@@@@@@@@@@@@@@@@@@@@@@@@@@@@@@@@@@@@@@@@@@@@@@@@@@@@@@@@@@
%@@@@@@@@@@@@@@@@@@@@@@@@@@@@@@@@@@@@@@@@@@@@@@@@@@@@@@@@@@
%@@@@@@@@@@@@@@@@@@@@@@@@@@@@@@@@@@@@@@@@@@@@@@@@@@@@@@@@@@
%----------------------------------------------------------
\par Our second main result, Theorem \reff{MainLinear}, states that the Banach couple $(\LOP,\LPM)$ is {\it quasi-linearizable}, see \cite{BL}. In other words, for every function $f\in \sum=\LOP+\LPM$ the functions $f_1\in\LOP$ and $f_2\in\LPM$ of an almost optimal decomposition $f=f_1+f_2$ can be chosen to {\it depend linearly} on $f$.
%----------------------------------------------------------
%@@@@@@@@@@@@@@@@@@@@@@@@@@@@@@@@@@@@@@@@@@@@@@@@@@@@@@@@@@
%@@@@@@@@@@@@@@@@@@@@@@@@@@@@@@@@@@@@@@@@@@@@@@@@@@@@@@@@@@
%@@@@@@@@@@@@@@@@@@@@@@@@@@@@@@@@@@@@@@@@@@@@@@@@@@@@@@@@@@
%@@@@@@@@@@@@@@@@@@@@@@@@@@@@@@@@@@@@@@@@@@@@@@@@@@@@@@@@@@
%@@@@@@@@@@@@@@@@@@@@@@@@@@@@@@@@@@@@@@@@@@@@@@@@@@@@@@@@@@
%@@@@@@@@@@@@@@@@@@@@@@@@@@@@@@@@@@@@@@@@@@@@@@@@@@@@@@@@@@
%----------------------------------------------------------
\begin{theorem}\lbl{MainLinear} Let $n<p<\infty$ and let $\mu$ be a non-trivial non-negative Borel measure on $\RN$. There exist continuous linear operators
%----------------------------------------------------------
$$
T_1:\LOP+\LPM\to \LOP~~~~\text{and}~~~~T_2:\LOP+\LPM\to \LPM
$$
%----------------------------------------------------------
such that
%----------------------------------------------------------
$$
T_1+T_2=Id_{\Sigma}
$$
%----------------------------------------------------------
and
%----------------------------------------------------------
$$
\|T_1\|_{\ssmall\to\LOP}+\|T_2\|_{\ssmall\to\LPM}\le C.
$$
%----------------------------------------------------------
Here $C=C(n,p)$ is a constant depending only on $n$ and $p$.
%----------------------------------------------------------
\end{theorem}
%----------------------------------------------------------
%@@@@@@@@@@@@@@@@@@@@@@@@@@@@@@@@@@@@@@@@@@@@@@@@@@@@@@@@@@
%@@@@@@@@@@@@@@@@@@@@@@@@@@@@@@@@@@@@@@@@@@@@@@@@@@@@@@@@@@
%@@@@@@@@@@@@@@@@@@@@@@@@@@@@@@@@@@@@@@@@@@@@@@@@@@@@@@@@@@
%@@@@@@@@@@@@@@@@@@@@@@@@@@@@@@@@@@@@@@@@@@@@@@@@@@@@@@@@@@
%@@@@@@@@@@@@@@@@@@@@@@@@@@@@@@@@@@@@@@@@@@@@@@@@@@@@@@@@@@
%----------------------------------------------------------
\par Let us briefly describe the main ideas of the proof of Theorem \reff{MAIN-CR}.  The necessity part of the proof, which we present in {\bf Section 2}, is based on the classical Sobolev-Poincar\'e inequality for $\LOP$-functions whenever $p>n$ and the Hardy-Littlewood maximal theorem.\medskip
%----------------------------------------------------------
\par We prove the sufficiency part of the theorem in two steps. The first step is presented in {\bf Section 3} where we construct a closed subset $E\subset\RN$  and a certain family $\Kc_E$ of pairwise disjoint ``well separated'' cubes of $\RN$ with centers in $E$.
%----------------------------------------------------------
\par This family of cubes possesses certain measure concentration properties (with respect to the measure $\mu$). In particular, $\mu(K)\sim (\diam K)^{n-p}$ for every cube $K\in\Kc_E$. We also prove that, on the other hand,  if $Q$ is a cube in $\RN$, $\theta>0$ and $
\diam Q\le\theta\,\dist(Q,E)$, then
$\mu(Q)\le C(\diam Q)^{n-p}$ where $C$ is a constant depending only on $p$ and $\theta$.\medskip
%----------------------------------------------------------
\par In {\bf Section 4}, given a function $f:\RN\to\R$ satisfying the sufficiency conditions, we construct the functions $f_1$ and $f_2=f-f_1$ of an almost optimal decomposition of $f$. We start by defining a function $\tf$ on $E$ by the formula
%----------------------------------------------------------
$$
\tf(x):=\frac{1}{\mu(K^{(x)})}\intl_{K^{(x)}}f\,d\mu,~~~x\in E,
$$
%----------------------------------------------------------
where $K^{(x)}$ denotes the (unique) cube from $\Kc_E$ centered at $x$. Then we extend $\tf$ from $E$ to all of $\RN$ using {\it the classical Whitney's extension method}. This gives us $f_1$ (and therefore of course also $f_2=f-f_1$). Section 4 also includes a proof  that the function $f_1$ satisfies the inequality $\|f_1\|_{\LOP}\le C(n,p)\,\lambda^{\frac1p}$. \medskip
%----------------------------------------------------------
\par In {\bf Section 5} we show that the function $f_2$ satisfies the inequality $\|f_2\|_{\LPM}\le C(n,p)\,\lambda^{\frac1p}$. This and the previous inequality prove the sufficiency part of Theorem \reff{MAIN-CR}.
%----------------------------------------------------------
\par Note that {\it the Whitney extension operator is linear}, so that {\it the functions $f_1$ and $f_2$ depend linearly on $f$}. This proves Theorem \reff{MainLinear}.\medskip
%----------------------------------------------------------
\par In {\bf Section 6} we prove several refinements of Theorem \reff{MAIN-CR}. Note that the criterion for the norm in the space $\Sigma=\LOP+\LPM$ given in this theorem describes the structure of $\Sigma$ and shows which properties of a function $f$ on $\RN$ control its almost optimal decomposition into a sum of a function from  $\LOP$ and a function from $\LPM$. At the same time it is not quite clear how one could check the conditions \rf{CR} of Theorem \reff{MAIN-CR} for a given function $f$ on $\RN$. In fact, these conditions depend on an infinite number of families $\Qc$ of cubes and all possible choices of cubes $Q',Q''\in \Qc$ satisfying condition \rf{IM-QP}.
%----------------------------------------------------------
\par Nevertheless a careful examination of our proof of Theorem \reff{MAIN-CR} shows that it constructs {\it a particular family} $\Qc$ of cubes and {\it particular mappings} $Q\mapsto Q'$ and $Q\mapsto Q''$ satisfying  \rf{IM-QP} depending only on $p$ and the measure $\mu$, and that it is enough to examine the behavior of $f$ only on this particular family and these particular mappings.
%----------------------------------------------------------
\par We express this fact by Theorem \reff{REF-MAIN-CR} which refines one part of the criterion of Theorem \reff{MAIN-CR}.
%----------------------------------------------------------
\par The next refinement of this result, Theorem \reff{REF-2}, enables us to express the norm of an arbitrary function $f\in L_{p,loc}(\RN;\mu)$ as a linear combination of $p$-oscillations of $f$ over a certain family of subsets in $\RN$ with fixed covering multiplicity. Note that the coefficients of this linear combination and the family of subsets depend only on $n,p,$ and the measure $\mu$.
%----------------------------------------------------------
\par We prove this result in Subsection 6.3. Remark that this rather specifical refinement of the main result has important applications to problems of characterizations of restrictions of Sobolev functions to closed subsets of $\RN$. (See a discussion at the end of this section.)
%----------------------------------------------------------
\par The proof of Theorem \reff{REF-2} is based on a new approach to extensions of functions which we call {\it a lacunary modification} of the Whitney extension method. We present this approach in Subsection 6.2. The main idea of this modification is to use certain {\it families} of Whitney's cubes rather than to treat each Whitney cube separately. We call these families of Whitney cubes {\it lacunae}. Each lacuna characterizes a certain ``hole'' in the complement $\RN\setminus E$.
%----------------------------------------------------------
\par In Subsection 6.2 we present main definitions and main properties of lacunae. For the proof of these properties we refer the reader to the paper \cite{S6}, Sections 4-5.
%----------------------------------------------------------
\par In Subsection 7.1 of {\bf Section 7} we prove several variants of the main result. Let us formulate one of them.
%----------------------------------------------------------
%@@@@@@@@@@@@@@@@@@@@@@@@@@@@@@@@@@@@@@@@@@@@@@@@@@@@@@@@@@
%@@@@@@@@@@@@@@@@@@@@@@@@@@@@@@@@@@@@@@@@@@@@@@@@@@@@@@@@@@
%@@@@@@@@@@@@@@@@@@@@@@@@@@@@@@@@@@@@@@@@@@@@@@@@@@@@@@@@@@
%@@@@@@@@@@@@@@@@@@@@@@@@@@@@@@@@@@@@@@@@@@@@@@@@@@@@@@@@@@
%@@@@@@@@@@@@@@@@@@@@@@@@@@@@@@@@@@@@@@@@@@@@@@@@@@@@@@@@@@
%@@@@@@@@@@@@@@@@@@@@@@@@@@@@@@@@@@@@@@@@@@@@@@@@@@@@@@@@@@
%----------------------------------------------------------
\begin{theorem}\lbl{CR-V1} Let $n<p<\infty$ and let $\mu$ be a non-trivial non-negative Borel measure on $\RN$.  A function $f\in L_{p,loc}(\RN;\mu)$ belongs to the space $\LOP+\LPM$ if and only if there exists a  positive constant $\lambda$ which satisfies the following conditions for a certain absolute positive constant $\gamma$: Let $\Qc$ be an arbitrary finite family of pairwise disjoint cubes in $\RN$. Suppose that to each cube $Q\in\mathcal{Q}$ we have arbitrarily assigned two cubes $Q',Q''\in\Qc$ such that $Q'\cup Q''\subset \gamma Q$ and
%----------------------------------------------------------
\bel{M-DM-V1}
(\diam Q')^{p-n}\mu(Q')+(\diam Q'')^{p-n}\mu(Q'')\le 1.
\ee
%----------------------------------------------------------
Then the following inequality
%----------------------------------------------------------
\bel{IN-V1}
\sbig_{Q\in\Qc}\,\,
\left(\frac{\diam Q' \diam Q''}{\diam Q}\right)^{p-n} \iint \limits_{Q'\times Q''}
|f(x)-f(y)|^p\, d\mu(x)d\mu(y)
\le \lambda
\ee
%----------------------------------------------------------
holds. Furthermore,
$\|f\|_{\sum}\sim \inf \lambda^{\frac{1}{p}}$ with constants of equivalence depending only on $n$ and $p$.
%----------------------------------------------------------
\end{theorem}
%----------------------------------------------------------
%@@@@@@@@@@@@@@@@@@@@@@@@@@@@@@@@@@@@@@@@@@@@@@@@@@@@@@@@@@
%@@@@@@@@@@@@@@@@@@@@@@@@@@@@@@@@@@@@@@@@@@@@@@@@@@@@@@@@@@
%@@@@@@@@@@@@@@@@@@@@@@@@@@@@@@@@@@@@@@@@@@@@@@@@@@@@@@@@@@
%@@@@@@@@@@@@@@@@@@@@@@@@@@@@@@@@@@@@@@@@@@@@@@@@@@@@@@@@@@
%@@@@@@@@@@@@@@@@@@@@@@@@@@@@@@@@@@@@@@@@@@@@@@@@@@@@@@@@@@
%@@@@@@@@@@@@@@@@@@@@@@@@@@@@@@@@@@@@@@@@@@@@@@@@@@@@@@@@@@
%----------------------------------------------------------
\par Note that the hypotheses of this theorem are equivalent to the hypotheses of Theorem \reff{MAIN-CR} provided the cubes $Q',Q''$ from its formulation satisfy inequality \rf{M-DM-V1}. Thus the sufficiency part of Theorem \reff{CR-V1} is slightly stronger than the sufficiency part of Theorem \reff{MAIN-CR}: it asserts that it suffices to verify \rf{CR} only for cubes satisfying inequality \rf{M-DM-V1} rather than for {\it all} cubes, as required in Theorem \reff{MAIN-CR}.
%----------------------------------------------------------
%@@@@@@@@@@@@@@@@@@@@@@@@@@@@@@@@@@@@@@@@@@@@@@@@@@@@@@@@@@
%@@@@@@@@@@@@@@@@@@@@@@@@@@@@@@@@@@@@@@@@@@@@@@@@@@@@@@@@@@
%----------------------------------------------------------
\par In Subsection 7.2, we obtain another variant of Theorem\reff{MAIN-CR}, which we use in Subsection 7.3 to prove the following explicit formula for calculation of the $K$-functional for the couple $(\LPM,\LOP)$.
%----------------------------------------------------------
%@@@@@@@@@@@@@@@@@@@@@@@@@@@@@@@@@@@@@@@@@@@@@@@@@@@@@@@@@@
%@@@@@@@@@@@@@@@@@@@@@@@@@@@@@@@@@@@@@@@@@@@@@@@@@@@@@@@@@@
%@@@@@@@@@@@@@@@@@@@@@@@@@@@@@@@@@@@@@@@@@@@@@@@@@@@@@@@@@@
%@@@@@@@@@@@@@@@@@@@@@@@@@@@@@@@@@@@@@@@@@@@@@@@@@@@@@@@@@@
%@@@@@@@@@@@@@@@@@@@@@@@@@@@@@@@@@@@@@@@@@@@@@@@@@@@@@@@@@@
%@@@@@@@@@@@@@@@@@@@@@@@@@@@@@@@@@@@@@@@@@@@@@@@@@@@@@@@@@@
%----------------------------------------------------------
\begin{theorem} Let $n<p<\infty$ and let $f\in L_{p,loc}(\RN;\mu)$. Then, for every $t>0$,
%----------------------------------------------------------
\be
&&K(t; f :(\LPM,\LOP))\nn\\&\sim&
\sup
\left\{
\sbig_{Q\in\Qc}\,\,
\left(\frac{\diam Q' \diam Q''}{\diam Q}\right)^{p-n} \frac{\iint \limits_{Q'\times Q''}
|f(x)-f(y)|^p\, d\mu(x)d\mu(y)}
{(\diam Q')^{p-n}\mu(Q')+(\diam Q'')^{p-n}\mu(Q'')}
\right\}^{\frac1p}\nn
\ee
%----------------------------------------------------------
where the supremum is taken over all finite families $\Qc$ of pairwise disjoint cubes in $\RN$ and all mappings $\Qc\ni Q\mapsto Q'\in\Qc$ and $\Qc\ni Q\mapsto Q''\in\Qc$
such that $Q'\cup Q''\subset \gamma Q$ and
%----------------------------------------------------------
$$
(\diam Q')\left(\frac{\mu(Q')}{|Q'|}\right)^{\frac1p}
+(\diam Q'')\left(\frac{\mu(Q'')}{|Q''|}\right)^{\frac1p}
\le t.
$$
%----------------------------------------------------------
Here $\gamma$ is an absolute constant. Furthermore, the above equivalence holds with constants depending only on $n$ and $p$.
%----------------------------------------------------------
\end{theorem}
%----------------------------------------------------------
%@@@@@@@@@@@@@@@@@@@@@@@@@@@@@@@@@@@@@@@@@@@@@@@@@@@@@@@@@@
%@@@@@@@@@@@@@@@@@@@@@@@@@@@@@@@@@@@@@@@@@@@@@@@@@@@@@@@@@@
%@@@@@@@@@@@@@@@@@@@@@@@@@@@@@@@@@@@@@@@@@@@@@@@@@@@@@@@@@@
%@@@@@@@@@@@@@@@@@@@@@@@@@@@@@@@@@@@@@@@@@@@@@@@@@@@@@@@@@@
%----------------------------------------------------------
\par Using Theorem \reff{MainLinear} we also prove that this formula for the $K$-functional of the couple $(\LPM,\LOP)$ can be quasi-linearized. See Subsection 7.2 for the details.\medskip
%----------------------------------------------------------
%@@@@@@@@@@@@@@@@@@@@@@@@@@@@@@@@@@@@@@@@@@@@@@@@@@@@@@@@@@
%@@@@@@@@@@@@@@@@@@@@@@@@@@@@@@@@@@@@@@@@@@@@@@@@@@@@@@@@@@
%@@@@@@@@@@@@@@@@@@@@@@@@@@@@@@@@@@@@@@@@@@@@@@@@@@@@@@@@@@
%----------------------------------------------------------
\par Finally, in Subsection 7.3 we give a geometrical interpretation of Theorem \reff{VRN}
and simple geometrical proofs of some particular cases of it, some of which have been kindly provided by V. Dolnikov.
%----------------------------------------------------------
%@@@@@@@@@@@@@@@@@@@@@@@@@@@@@@@@@@@@@@@@@@@@@@@@@@@@@@@@@@
%@@@@@@@@@@@@@@@@@@@@@@@@@@@@@@@@@@@@@@@@@@@@@@@@@@@@@@@@@@
%@@@@@@@@@@@@@@@@@@@@@@@@@@@@@@@@@@@@@@@@@@@@@@@@@@@@@@@@@@
%----------------------------------------------------------
\par Our interest in Problems \reff{MAIN-PR1} and \reff{MAIN-PR2} has been motivated by their intimate connection with
the characterization of the restrictions of Sobolev $L^2_p(\RN)$-functions to arbitrary closed subsets of $\RN$. In particular, Theorem \reff{MAIN-CR} is one of the main ingredients of our approach to this problem in \cite{S6} where it enables us to give a constructive description of the trace space  $L^2_p(\R^2)|_E$ whenever $p>2$ and $E$ is an arbitrary finite set $E\subset\R^2$.
%----------------------------------------------------------
\par Our second main result here, Theorem \reff{MainLinear}, is also used in \cite{S6} in order to prove the existence of a {\it continuous linear extension operator} from $L^2_p(\R^2)|_E$ into $L^2_p(\R^2)$, $p>2$, whose operator norm is bounded by a constant depending only on $p$. A different proof of this latter result has been given earlier by A. Israel \cite{Is}. Quite recently  C. Fefferman, A. Israel and G. K. Luli \cite{FIL} proved the existence of such an operator for the space $\LMP|_E$ whenever $n<p<\infty$ and $E\subset\RN$ is an arbitrary closed set. We refer to \cite{S6} for more details.
%----------------------------------------------------------
\medskip
%----------------------------------------------------------
%@@@@@@@@@@@@@@@@@@@@@@@@@@@@@@@@@@@@@@@@@@@@@@@@@@@@@@@@@@
%@@@@@@@@@@@@@@@@@@@@@@@@@@@@@@@@@@@@@@@@@@@@@@@@@@@@@@@@@@
%----------------------------------------------------------
%@@@@@@@@@@@@@@@@@@@@@@@@@@@@@@@@@@@@@@@@@@@@@@@@@@@@@@@@@@
%@@@@@@@@@@@@@@@@@@@@@@@@@@@@@@@@@@@@@@@@@@@@@@@@@@@@@@@@@@
\par {\bf Acknowledgements.} I am very thankful to M. Cwikel for useful suggestions and remarks. I am pleased to thank V. Dolnikov for very useful discussions of some geometrical aspects of Theorem \reff{VRN}.
I am also very grateful to C. Fefferman, N. Zobin and all the participants of ``Whitney Problems Workshop'', Williamsburg, August 2011, for stimulating discussions and valuable advice.
%----------------------------------------------------------
%@@@@@@@@@@@@@@@@@@@@@@@@@@@@@@@@@@@@@@@@@@@@@@@@@@@@@@@@@@
%@@@@@@@@@@@@@@@@@@@@@@@@@@@@@@@@@@@@@@@@@@@@@@@@@@@@@@@@@@
%@@@@@@@@@@@@@@@@@@@@@@@@@@@@@@@@@@@@@@@@@@@@@@@@@@@@@@@@@@
%----------------------------------------------------------
%@@@@@@@@@@@@@@@@@@@@@@@@@@@@@@@@@@@@@@@@@@@@@@@@@@@@@@@@@@
%@@@@@@@@@@@@@@@@@@@@@@@@@@@@@@@@@@@@@@@@@@@@@@@@@@@@@@@@@@
%@@@@@@@@@@@@@@@@@@@@@@@@@@@@@@@@@@@@@@@@@@@@@@@@@@@@@@@@@@
%@@@@@@@@@@@@@@@@@@@@@@@@@      @@@@@@@@@@@@@@@@@@@@@@@@@@@
%@@@@@@@@@@@@@@@@@@@@@@@          @@@@@@@@@@@@@@@@@@@@@@@@@
%@@@@@@@@@@@@@@@@@@@@@              @@@@@@@@@@@@@@@@@@@@@@@
%@@@@@@@@@@@@@@@@@@@     SECTION 2    @@@@@@@@@@@@@@@@@@@@@
%@@@@@@@@@@@@@@@@@@@@@              @@@@@@@@@@@@@@@@@@@@@@@
%@@@@@@@@@@@@@@@@@@@@@@@          @@@@@@@@@@@@@@@@@@@@@@@@@
%@@@@@@@@@@@@@@@@@@@@@@@@@      @@@@@@@@@@@@@@@@@@@@@@@@@@@
%@@@@@@@@@@@@@@@@@@@@@@@@@@@@@@@@@@@@@@@@@@@@@@@@@@@@@@@@@@
%@@@@@@@@@@@@@@@@@@@@@@@@@@@@@@@@@@@@@@@@@@@@@@@@@@@@@@@@@@
%----------------------------------------------------------
%@@@@@@@@@@@@@@@@@@@@@@@@@@@@@@@@@@@@@@@@@@@@@@@@@@@@@@@@@@
%----------------------------------------------------------
\SECT{2. Proof of Theorem \reff{MAIN-CR}: Necessity.}{2}
%----------------------------------------------------------
\addtocontents{toc}{2. Proof of Theorem \reff{MAIN-CR}: Necessity. \hfill \thepage\\\par}
%----------------------------------------------------------
%@@@@@@@@@@@@@@@@@@@@@@@@@@@@@@@@@@@@@@@@@@@@@@@@@@@@@@@@@@
%----------------------------------------------------------
\indent %----------------------------------------------------------
\par Throughout the paper $C,C_1,C_2,...$ and $\gamma,\gamma_1,\gamma_2,...$ will be generic
positive constants which depend only on  $n$ and $p$. Sometimes these constants can depend on certain parameters (say $\eta,\theta,$ etc.) which we fix in formulations of some auxiliary results. These constants can change even in a single string of estimates. The dependence of a constant on certain parameters is expressed, for example, by the notation $C=C(n,p)$ or $\gamma=\gamma(n)$. We write $A\sim B$ if there is a constant $C\ge 1$ such that $A/C\le B\le CA$.
%----------------------------------------------------------
\par Throughout the paper the words ``a subset of $\RN$'' will mean ``a Borel subset of $\RN$''. For a locally integrable (with respect to the measure $\mu$) function $f$ and a subset $S\subset\RN$ of a positive $\mu$-measure by $f_{S}$ we denote the $\mu$-average of $f$ over $S$:
%----------------------------------------------------------
$$
\av{S}:=\frac{1}{\mu(S)}\intl_S f\,d\mu.
$$
%----------------------------------------------------------
\par By $\|\cdot\|$ we denote the uniform measure in $\RN$. Given a set $A$ by $\card A$ we denote the cardinality of $A$.
%----------------------------------------------------------
\par Let $\Ac$ be a family of sets in $\RN$. By $M(\Ac)$ we denote its covering multiplicity, i.e., the minimal positive integer M such that every point $x\in\RN$ is covered by at most M sets from $\Ac$. Finally, given a function $g\in L_{1,loc}(\RN)$ we let $\Mc[g]$ denote its Hardy-Littlewood maximal function:
%----------------------------------------------------------
\bel{HLM}
\Mc[g](x):=\sup_{K\ni x}\frac{1}{|K|}\intl_{K} g(y)\,dy, ~~~~x\in\RN.
\ee
%----------------------------------------------------------
As usual, in this formula the supremum is taken over all cubes $K$ in $\RN$ containing $x$.
%----------------------------------------------------------
\par When $p>n$, it follows from the Sobolev embedding theorem that every function $F\in\LOP$ coincides almost everywhere with a continuous function.  This fact enables us {\it to identify each element $F\in \LOP$, $p>n$, with its unique continuous representative}.
%----------------------------------------------------------
%@@@@@@@@@@@@@@@@@@@@@@@@@@@@@@@@@@@@@@@@@@@@@@@@@@@@@@@@@@
%---------------------------------------------------------- %@@@@@@@@@@@@@@@@@@@@@@@@@@@@@@@@@@@@@@@@@@@@@@@@@@@@@@@@@@
\par One of the main tools of the proof of Theorem \reff{MAIN-CR} is the following proposition which presents a classical Sobolev imbedding inequality for the case $p>n$, see, e.g. \cite{M}, p. 61, or \cite{MP}, p. 55. This inequality is also known in the literature as Sobolev-Poincar\'e inequality (for $p>n$).
%----------------------------------------------------------
%@@@@@@@@@@@@@@@@@@@@@@@@@@@@@@@@@@@@@@@@@@@@@@@@@@@@@@@@@@
\begin{proposition}\lbl{EDIFF} Let $F\in \LOP$ be
a continuous function and let $n<q\le p<\infty$. Then for
every cube $Q\subset\RN$ and every $x,y\in Q$ the
following inequality
%----------------------------------------------------------
$$ |F(x)-F(y)|\le C(n,q)\, \diam Q \left(\frac{1}{|Q|}\intl_Q
\|\nabla F(z)\|^q\,dz\right) ^{\frac{1}{q}} $$
%----------------------------------------------------------
holds.
\end{proposition}
%@@@@@@@@@@@@@@@@@@@@@@@@@@@@@@@@@@@@@@@@@@@@@@@@@@@@@@@@@@
%@@@@@@@@@@@@@@@@@@@@@@@@@@@@@@@@@@@@@@@@@@@@@@@@@@@@@@@@@@
%@@@@@@@@@@@@@@@@@@@@@@@@@@@@@@@@@@@@@@@@@@@@@@@@@@@@@@@@@@
%----------------------------------------------------------
\par We begin the proof of the necessity part with the following auxiliary lemma.
%----------------------------------------------------------
%@@@@@@@@@@@@@@@@@@@@@@@@@@@@@@@@@@@@@@@@@@@@@@@@@@@@@@@@@@
\begin{lemma}\lbl{N-D} Let $n<q\le p<\infty$, $\gamma>1$, and $\theta',\theta''>0$. Let $f=f_1+f_2$ where $f_1\in\LOP$ and $f_2\in\LPM$. Given a cube $Q\subset\RN$ let $S',S''$ be closed subsets of $Q$ such that $S'\cup S''\subset\gamma Q$.
%----------------------------------------------------------
\par Then  the following inequality
%----------------------------------------------------------
\be
\frac{\frac{1}{\mu(S')\mu(S'')}
\iint \limits_{S'\times S''}
|f(x)-f(y)|^p\, d\mu(x)d\mu(y)}
{(\diam Q)^{p-n}+\theta'/\mu(S')+\theta''/\mu(S'')}&\le&
C\left\{\intl_Q\Mc[(\|\nabla f_1\|)^q]^{\frac{p}{q}}(x)\,dx\right.\nn\\
&+&\left.\frac{1}{\theta'}
\intl_{S'}|f_2(x)|^p\,d\mu(x)+\frac{1}{\theta''}
\intl_{S''}|f_2(x)|^p\,d\mu(x)\right\}\nn
\ee
%----------------------------------------------------------
holds. Here $C=C(n,p,q,\gamma)$ is a constant depending only $n,p,q$ and $\gamma$.
\end{lemma}
%@@@@@@@@@@@@@@@@@@@@@@@@@@@@@@@@@@@@@@@@@@@@@@@@@@@@@@@@@@
%@@@@@@@@@@@@@@@@@@@@@@@@@@@@@@@@@@@@@@@@@@@@@@@@@@@@@@@@@@
%@@@@@@@@@@@@@@@@@@@@@@@@@@@@@@@@@@@@@@@@@@@@@@@@@@@@@@@@@@
\par {\it Proof.} We have
%----------------------------------------------------------
\be
I&:=&\frac{1}{\mu(S')\mu(S'')}
\iint \limits_{S'\times S''}
|f(x)-f(y)|^p\, d\mu(x)d\mu(y)\nn\\&\le&
\frac{1}{\mu(S')\mu(S'')}
\iint \limits_{S'\times S''}
(|f_1(x)-f_1(y)|+|f_2(x)-f_2(y)|)^p\, d\mu(x)d\mu(y)\nn\\&\le&
\frac{2^p}{\mu(S')\mu(S'')}
\iint\limits_{S'\times S''}
|f_1(x)-f_1(y)|^p\, d\mu(x)d\mu(y)\nn\\
&+& \frac{2^p}{\mu(S')\mu(S'')}
\iint\limits_{S'\times S''}
|f_2(x)-f_2(y)|^p\,d\mu(x)d\mu(y)=2^p\{I_1+I_2\}.
\nn
\ee
%----------------------------------------------------------
\par By the Sobolev-Poincar\'e inequality, see Proposition \reff{EDIFF}, for every $x,y\in\gamma Q$ we have
%----------------------------------------------------------
$$ |f_1(x)-f_1(y)|\le C(n,q)\, \diam Q \left(\frac{1}{|\gamma Q|}\intl_{\gamma Q}
\|\nabla f_1(z)\|^q\,dz\right) ^{\frac{1}{q}}
$$
%----------------------------------------------------------
with $C=C(n,q,\gamma)$. Hence
%----------------------------------------------------------
\be
I_1&:=&
\frac{1}{\mu(S')\mu(S'')}
\iint\limits_{S'\times S''}
|f_1(x)-f_1(y)|^p\, d\mu(x)d\mu(y)\nn\\&\le& C(\diam Q)^p
\left(\frac{1}{|\gamma Q|}\intl_{\gamma Q}
\|\nabla f_1(z)\|^q\,dz\right) ^{\frac{p}{q}}.\nn
\ee
%----------------------------------------------------------
Then, by \rf{HLM}, for every $z\in Q$ we have
%----------------------------------------------------------
$$
\left(\frac{1}{|\gamma Q|}\intl_{\gamma Q}
\|\nabla f_1(z)\|^q\,dz\right)^{\frac{p}{q}}\le
\Mc[\,\|\nabla f_1(z)\|^q]^{\frac{p}{q}}(z).
$$
%----------------------------------------------------------
Integrating this inequality on  $Q$ (with respect to $z$) we obtain
%----------------------------------------------------------
\bel{MA-S1}
|Q|\left(\frac{1}{|\gamma Q|}\intl_{\gamma Q}
\|\nabla f_1(z)\|^q\,dz\right)^{\frac{p}{q}}\le
\intl_{Q}\Mc[\,\|\nabla f_1(z)\|^q]^{\frac{p}{q}}(z)\,dz.
\ee
%----------------------------------------------------------
Hence
%----------------------------------------------------------
\bel{NI1}
I_1\le C(\diam Q)^{p-n}
\intl_{Q}\Mc[\,\|\nabla f_1(z)\|^q]^{\frac{p}{q}}(z)\,dz.
\ee
%----------------------------------------------------------
\par Let us estimate the quantity $I_2$. We have
%----------------------------------------------------------
\be
I_2&:=&
\frac{1}{\mu(S')\mu(S'')}
\iint\limits_{S'\times S''}
|f_2(x)-f_2(y)|^p\, d\mu(x)d\mu(y)\nn\\&\le&
\frac{2^p}{\mu(S')\mu(S'')}\left\{
\iint\limits_{S'\times S''}
|f_2(x)|^p\, d\mu(x)d\mu(y)+
\iint\limits_{S'\times S''}
|f_2(y)|^p\, d\mu(x)d\mu(y)\right\}\nn\\
&=&
2^p\left\{\frac{1}{\mu(S')}
\intl_{S'}
|f_2(x)|^p\, d\mu(x)+
\frac{1}{\mu(S'')}
\intl_{S''}
|f_2(x)|^p\, d\mu(x)\right\}.\nn
\ee
%----------------------------------------------------------
Combining this inequality with inequality \rf{NI1} we obtain
%----------------------------------------------------------
\be
I&=&2^p\{I_1+I_2\}\le
C\left\{(\diam Q)^{p-n}\intl_Q\Mc[(\|\nabla f_1\|)^q] ^{\frac{p}{q}}(x)\,dx\right.\nn\\
&+&\left.\frac{1}{\mu(S')}
\intl_{S'}|f_2(x)|^p\,d\mu(x)+\frac{1}{\mu(S'')}
\intl_{S''}|f_2(x)|^p\,d\mu(x)\right\}.
\nn
\ee
%----------------------------------------------------------
Hence
%----------------------------------------------------------
\be
&&I/
\{(\diam Q)^{p-n}+\theta'/\mu(S')+\theta''/\mu(S'')\}\le
C\left\{\intl_Q\Mc[(\|\nabla f_1\|)^q]^{\frac{p}{q}}(x)\,dx\right.\nn\\
&+&\left.\frac{1}{\theta'}
\intl_{S'}|f_2(x)|^p\,d\mu(x)+\frac{1}{\theta''}
\intl_{S''}|f_2(x)|^p\,d\mu(x)\right\}\nn
\ee
%----------------------------------------------------------
proving the lemma.\bx\bigskip
%----------------------------------------------------------
%@@@@@@@@@@@@@@@@@@@@@@@@@@@@@@@@@@@@@@@@@@@@@@@@@@@@@@@@@@
%@@@@@@@@@@@@@@@@@@@@@@@@@@@@@@@@@@@@@@@@@@@@@@@@@@@@@@@@@@
%@@@@@@@@@@@@@@@@@@@@@@@@@@@@@@@@@@@@@@@@@@@@@@@@@@@@@@@@@@
\par We are in a position to prove a slightly more general version of the necessity part of Theorem \reff{MAIN-CR}.
%----------------------------------------------------------
%@@@@@@@@@@@@@@@@@@@@@@@@@@@@@@@@@@@@@@@@@@@@@@@@@@@@@@@@@@
%@@@@@@@@@@@@@@@@@@@@@@@@@@@@@@@@@@@@@@@@@@@@@@@@@@@@@@@@@@
%@@@@@@@@@@@@@@@@@@@@@@@@@@@@@@@@@@@@@@@@@@@@@@@@@@@@@@@@@@
%@@@@@@@@@@@@@@@@@@@@@@@@@@@@@@@@@@@@@@@@@@@@@@@@@@@@@@@@@@
%@@@@@@@@@@@@@@@@@@@@@@@@@@@@@@@@@@@@@@@@@@@@@@@@@@@@@@@@@@
%@@@@@@@@@@@@@@@@@@@@@@@@@@@@@@@@@@@@@@@@@@@@@@@@@@@@@@@@@@
%@@@@@@@@@@@@@@@@@@@@@@@@@@@@@@@@@@@@@@@@@@@@@@@@@@@@@@@@@@
%@@@@@@@@@@@@@@@@@@@@@@@@@@@@@@@@@@@@@@@@@@@@@@@@@@@@@@@@@@
%----------------------------------------------------------
\begin{proposition}\lbl{OP-W} Let $n<p<\infty, \gamma>1,N\ge 1,$ and let $\mu$ be a non-trivial non-negative Borel measure on $\RN$.  Let $\Qc$ be a family of cubes in $\RN$ with covering multiplicity $M(\Qc)\le N$, and let $\Sc$ be a finite family of closed subsets of $\RN$ of positive $\mu$-measure with $M(\Sc)\le N$. Suppose that to each cube $Q\in\Qc$ we have assigned two subsets $S'_Q,S''_Q\in\Sc$ such that
%----------------------------------------------------------
$$
S'_Q\cup S''_Q\subset \gamma Q.
$$
%----------------------------------------------------------
\par If a function $f\in \LOP+\LPM$, then
%----------------------------------------------------------
\bel{N-11}
\sbig_{Q\in\Qc}\,\,
\frac{\frac{1}{\mu(S'_Q)\mu(S''_Q)}\iint \limits_{S'_Q\times S''_Q}
|f(x)-f(y)|^p\, d\mu(x)d\mu(y)}
{(\diam Q)^{p-n}\{1+(\diam S'_Q)^{n-p}/\mu(S'_Q)
+(\diam S''_Q)^{n-p}/\mu(S''_Q)\}}
\le \lambda
\ee
%----------------------------------------------------------
where $\lambda=C\|f\|^p_{\sum}$. Here $C$ is a constant depending only on $n,p,\gamma$ and $N$.
%----------------------------------------------------------
\end{proposition}
%----------------------------------------------------------
%@@@@@@@@@@@@@@@@@@@@@@@@@@@@@@@@@@@@@@@@@@@@@@@@@@@@@@@@@@
%@@@@@@@@@@@@@@@@@@@@@@@@@@@@@@@@@@@@@@@@@@@@@@@@@@@@@@@@@@
%@@@@@@@@@@@@@@@@@@@@@@@@@@@@@@@@@@@@@@@@@@@@@@@@@@@@@@@@@@
\par {\it Proof.} Let $q:=(p+n)/2$. Given a cube $Q\in\Qc$ we put $S':=S'_Q,
S'':=S''_Q,$ and
%---------------------------------------------------------- %@@@@@@@@@@@@@@@@@@@@@@@@@@@@@@@@@@@@@@@@@@@@@@@@@@@@@@@@@@
$$
\theta':=(\diam Q/\diam S'_Q)^{p-n},~~ \theta'':=(\diam Q/\diam S''_Q)^{p-n}.
$$
%---------------------------------------------------------- %@@@@@@@@@@@@@@@@@@@@@@@@@@@@@@@@@@@@@@@@@@@@@@@@@@@@@@@@@@
Since  $f\in \LOP+\LPM$, there exist functions $f_1\in\LOP$ and $f_2\in\LPM$ such that $f=f_1+f_2$ and
%---------------------------------------------------------- %@@@@@@@@@@@@@@@@@@@@@@@@@@@@@@@@@@@@@@@@@@@@@@@@@@@@@@@@@@
$$
\|f_1\|_{\LOP}\le 2\|f\|_{\sum}\,,~~~\|f_2\|_{\LPM}\le 2\|f\|_{\sum}\,.
$$
%---------------------------------------------------------- %@@@@@@@@@@@@@@@@@@@@@@@@@@@@@@@@@@@@@@@@@@@@@@@@@@@@@@@@@@
Then, by Lemma \reff{N-D}, the quantity
%----------------------------------------------------------
$$
J_Q:=\frac{\frac{1}{\mu(S'_Q)\mu(S''_Q)}\iint \limits_{S'_Q\times S''}
|f(x)-f(y)|^p\, d\mu(x)d\mu(y)}
{(\diam Q)^{p-n}\{1+(\diam S'_Q)^{n-p}/\mu(S'_Q)
+(\diam S''_Q)^{n-p}/\mu(S''_Q)\}}
$$
%----------------------------------------------------------
satisfies the following inequality
%----------------------------------------------------------
\be
J_Q&\le&
C\left\{\intl_Q\Mc[(\|\nabla f_1\|)^q]^{\frac{p}{q}}(x)\,dx+\frac{1}{\theta'}
\intl_{S'}|f_2(x)|^p\,d\mu(x)+\frac{1}{\theta''}
\intl_{S''}|f_2(x)|^p\,d\mu(x)\right\}\nn\\
&=& C\left\{\intl_Q\Mc[(\|\nabla f_1\|)^q]^{\frac{p}{q}}(x)\,dx\right.\nn\\
&+&
\left.\left(\frac{\diam S'_Q}{\diam Q}\right)^{p-n}
\intl_{S'_Q}|f_2(x)|^p\,d\mu(x)+
\left(\frac{\diam S''_Q}{\diam Q}\right)^{p-n}
\intl_{S''_Q}|f_2(x)|^p\,d\mu(x)\right\}\nn\\&=& C(J_1+J_2+J_3).\nn
\ee
%----------------------------------------------------------
\par Prove that
%----------------------------------------------------------
$$
J:=\sum_{Q\in\Qc}J_Q\le C\|f\|_{\sum}^p\,.
$$
%----------------------------------------------------------
\par First let us show that the following inequality
%----------------------------------------------------------
\bel{E-J1}
J_1:=\sum_{Q\in\Qc}\,\,\intl_Q\Mc[(\|\nabla f_1\|)^q]^{\frac{p}{q}}(x)\,dx
\le C\|f_1\|_{\LOP}^p
\ee
%----------------------------------------------------------
holds. In fact, since covering multiplicity of the family $\Qc$ is at most $N$ and $p/q>1$, by the Hardy-Littlewood maximal theorem,
%----------------------------------------------------------
$$
J_1\le N\intl_{\RN}\Mc[\|\nabla f_1\|^q]^{\frac{p}{q}}(x)\,dx\le C\intl_{\RN}(\|\nabla f_1\|^q])^{\frac{p}{q}}(x)\,dx\nn=
C\|\nabla f_1\|^p_{\LPRN}.
$$
%----------------------------------------------------------
%@@@@@@@@@@@@@@@@@@@@@@@@@@@@@@@@@@@@@@@@@@@@@@@@@@@@@@@@@@
%@@@@@@@@@@@@@@@@@@@@@@@@@@@@@@@@@@@@@@@@@@@@@@@@@@@@@@@@@@
%----------------------------------------------------------
\par Prove that
%----------------------------------------------------------
\bel{D-J2}
J_2:=\sbig_{Q\in\Qc}
\left(\frac{\diam S'_Q}{\diam Q}\right)^{p-n}
\intl_{S'_Q}|f_2(x)|^p\,d\mu(x)
\le C\|f_2\|^p_{\LPM}\,.
\ee
%----------------------------------------------------------
\par To this end let us fix a set $S\in\Sc$ and prove that the quantity
%----------------------------------------------------------
$$
I(S):=\sum\{(\diam Q)^{n-p}:~Q\in\Qc, S'_Q=S\}
$$
%----------------------------------------------------------
satisfies the following inequality
%----------------------------------------------------------
\bel{IS}
I(S)\le C(\diam S)^{n-p}.
\ee
%----------------------------------------------------------
\par Recall that
%----------------------------------------------------------
$$
S=S'_Q\subset\gamma Q~~~\text{for every}~~~Q\in \Qc,
$$
%----------------------------------------------------------
and that $M(\Qc)\le N$.
%----------------------------------------------------------
\par Fix a point $a\in S$ and put $K_S:=Q(a,\diam S)$. Define three subfamilies of the family $\Qc$:
%----------------------------------------------------------
$$
\Qc_S:=\{Q\in\Qc: S'_Q=S\},
$$
%----------------------------------------------------------
%@@@@@@@@@@@@@@@@@@@@@@@@@@@@@@@@@@@@@@@@@@@@@@@@@@@@@@@@@@
%----------------------------------------------------------
$$
\Qc^{(1)}_S:=\{Q\in\Qc: S'_Q=S,~ Q\cap K_S\ne\emp\},
$$
%----------------------------------------------------------
and
%----------------------------------------------------------
$$
\Qc^{(2)}_S:=\{Q\in\Qc: S'_Q=S,~ Q\cap K_S=\emp\}.
$$
%----------------------------------------------------------
\par Prove that $\Qc^{(1)}_S$ contains at most $N=N(n,\gamma)$ elements. We will make use of the following simple statement: Let $Q_1,Q_2$ be cubes in $\RN$ such that $Q_1\cap Q_2\ne\emp$. Then the set $Q_1\cap(2Q_2)$ contains a cube $\tQ$ such that
%----------------------------------------------------------
$$
\diam \tQ\ge \tfrac{1}{2}\min\{\diam Q_1,\diam Q_2\}.
$$
%----------------------------------------------------------
In fact, suppose that $\diam Q_1\le \frac{1}{2}\diam Q_2$. Since $Q_1\cap Q_2\ne\emp$, we have $Q_1\subset 2Q_2$ so that we can put $\tQ:=Q_1.$
%----------------------------------------------------------
\par Assume that $\diam Q_1>\frac{1}{2}\diam Q_2$. Let $y\in Q_1\cap Q_2$. Then there exists a cube $Q^{(y)}_1\subset Q_1$ such that $Q^{(y)}_1\ni y$ and
%----------------------------------------------------------
$$
\diam Q^{(y)}_1=\tfrac{1}{2}\diam Q_2.
$$
%----------------------------------------------------------
Using the same argument as in the first case we conclude that  $Q^{(y)}_1\subset 2Q_2$. Thus $Q^{(y)}_1\subset Q_1\cap(2Q_2)$ so that we can put $\tQ:=Q^{(y)}_1.$ Then
%----------------------------------------------------------
$$
\diam \tQ=\frac{1}{2}\diam Q_2\ge \min\{\diam Q_1,\diam Q_2\}
$$
%----------------------------------------------------------
proving the statement.
%----------------------------------------------------------
\par Let us prove the required inequality
%----------------------------------------------------------
\bel{CA1}
\card \Qc_S^{(1)}\le C(n,\gamma,N).
\ee
%----------------------------------------------------------
If a cube $Q\in \Qc_S^{(1)}$, then $Q\cap K_S\ne\emp$ so that, by the above statement, there exists a cube
$\tQ$ such that $\tQ\subset Q\cap(2K_S)$ and
%----------------------------------------------------------
$$
\diam \tQ\ge \tfrac{1}{2}\min\{\diam Q,K_S\}.
$$
%----------------------------------------------------------
But $S\subset\gamma Q$ so that $\diam S\le \gamma\diam Q$.
Since $K_S=Q(a,\diam S)$, we obtain
%----------------------------------------------------------
$$
\diam K_S=2\diam S\le 2\gamma\diam Q.
$$
%----------------------------------------------------------
Hence
%----------------------------------------------------------
$$
\diam \tQ\ge (1/4\gamma) K_S.
$$
%----------------------------------------------------------
Thus $2K_S\supset \tQ$ and
$\diam\tQ\ge (1/4\gamma) \diam K_S.$ Note that the family $\Qc^{(1)}_S\subset\Qc$ has covering multiplicity $M(\Qc^{(1)}_S)\le M(\Qc)\le N$ so that $M(\{\tQ:Q\in\Qc^{(1)}_S\})\le N$ as well. Clearly, the cube $2K_S$ can contain at most $C(n,\gamma,N)$ of cubes $\tQ$ of diameter at least $(1/4\gamma) \diam K_S.$ This proves \rf{CA1}.
%----------------------------------------------------------
\par Now we have
%----------------------------------------------------------
$$
I_1(S):=\sum\{(\diam Q)^{n-p}:Q\in\Qc^{(1)}_S\}\le (\card \Qc^{(1)}_S)\max\{(\diam Q)^{n-p}:Q\in \Qc^{(1)}_S\}.
$$
%----------------------------------------------------------
Since $\card \Qc_S^{(1)}\le N(n,\gamma)$ and
%----------------------------------------------------------
$$
\diam S\le \gamma\diam Q, ~~~Q\in \Qc_S,
$$
%----------------------------------------------------------
we conclude that
%----------------------------------------------------------
$$
I_1(S)\le C(n,\gamma)(\diam S)^{n-p}.
$$
%----------------------------------------------------------
\par Let us estimate the quantity
%----------------------------------------------------------
$$
I_2(S):=\sum\{(\diam Q)^{n-p}:Q\in\Qc^{(2)}_S\}.
$$
%----------------------------------------------------------
Recall that $K_S=Q(a,\diam S)$ where $a\in S$ and
%----------------------------------------------------------
$$
K_S\cap Q\ne\emp~~~\text{for every}~~~Q\in\Qc^{(2)}_S.
$$
%----------------------------------------------------------
Since $a\in S\subset\gamma Q$ for each $Q\in \Qc^{(2)}_S$
we have
%----------------------------------------------------------
$$
\|x-a\|\le \diam(\gamma Q)=\gamma\diam Q~~~\text{for every}~~~x\in Q.
$$
%----------------------------------------------------------
Hence
%----------------------------------------------------------
$$
(\diam Q)^{n-p}\le C(\diam Q)^{-p}|Q|\le C\|x-a\|^{-p}|Q|,~~x\in Q.
$$
%----------------------------------------------------------
Integrating this inequality over the cube $Q$ (with respect to $x$), we obtain
%----------------------------------------------------------
$$
(\diam Q)^{n-p}\le C\intl_Q\|x-a\|^{-p}\,dx, ~~~Q\in\Qc^{(2)}_S.
$$
%----------------------------------------------------------
Since $K_S\cap Q\ne\emp$ for every $Q\in\Qc^{(2)}_S$, we have
%----------------------------------------------------------
$$
U_S:=\cup\{Q:Q\in\Qc^{(2)}_S\}\subset\RN\setminus K_S.
$$
%----------------------------------------------------------
Since $M(\Qc^{(2)}_S)\le N$, we obtain
%----------------------------------------------------------
\be
I_2(S)&:=&\sum\{(\diam Q)^{n-p}:Q\in\Qc^{(2)}_S\}\le C\sum_{Q\in\Qc^{(2)}_S}\intl_Q\|x-a\|^{-p}\,dx\nn\\
&=&C\intl_{U_S}\|x-a\|^{-p}\,dx\le
C\intl_{\RN\setminus K_S}\|x-a\|^{-p}\,dx
\le C(\diam K_S)^{n-p}.
\nn
\ee
%----------------------------------------------------------
Since $\diam K_S\sim \diam S$, we have
%----------------------------------------------------------
$$
I_2(S)\le C(\diam S)^{n-p}.
$$
%----------------------------------------------------------
\par Finally we obtain
%----------------------------------------------------------
\bel{F1}
I(S)=I_1(S)+I_2(S)\le C(\diam S)^{n-p}
\ee
%----------------------------------------------------------
proving inequality \rf{IS}.
%----------------------------------------------------------
\par Using this inequality we have the following estimate of the quantity $J_2$ defined in \rf{D-J2}:
%----------------------------------------------------------
\be
J_2&=&\sbig_{Q\in\Qc}
\left(\frac{\diam S'_Q}{\diam Q}\right)^{p-n}
\intl_{S'_Q}|f_2(x)|^p\,d\mu(x)\nn\\
&=&\sbig_{S\in\Sc}\,\,\sbig_{Q\in\Qc,\,S_{Q'}=S}
\left(\frac{\diam S}{\diam Q}\right)^{p-n}
\intl_{S}|f_2(x)|^p\,d\mu(x)\nn\\
&=&\smed_{S\in\Sc}(\diam S)^{p-n} \left(\smed_{Q\in\Qc,\,S_{Q'}=S}
(\diam Q)^{n-p}\right)\intl_{S}|f_2(x)|^p\,d\mu(x).\nn
\ee
%----------------------------------------------------------
By \rf{F1},
%----------------------------------------------------------
\be
J_2&\le& C\,\smed_{S\in\Sc}\,
(\diam S)^{p-n}\, (\diam S)^{n-p} \intl_{S}|f_2(x)|^p\,d\mu(x)\nn\\&\le& C\intl_{\RN}|f_2(x)|^p\,d\mu(x)=C\|f_2\|^p_{\LPM}.\nn
\ee
%----------------------------------------------------------
\par In the same fashion we prove that
%----------------------------------------------------------
$$
J_3:=\sbig_{Q\in\Qc}
\left(\frac{\diam S''_Q}{\diam Q}\right)^{p-n}
\intl_{S''_Q}|f_2(x)|^p\,d\mu(x)
\le C\|f_2\|^p_{\LPM}.
$$
%----------------------------------------------------------
Finally, summarizing estimates for the quantities $J_1$, see \rf{E-J1}, $J_2$ and $J_3$, we obtain
%----------------------------------------------------------
$$
J\le C(J_1+J_2+J_3)\le C(\|f_1\|^p_{\LOP}+\|f_2\|^p_{\LPM})\le C\|f\|_{\sum}^p\,.
$$
%----------------------------------------------------------
The proposition is completely proved.\bx\smallskip
%----------------------------------------------------------
%@@@@@@@@@@@@@@@@@@@@@@@@@@@@@@@@@@@@@@@@@@@@@@@@@@@@@@@@@@
%@@@@@@@@@@@@@@@@@@@@@@@@@@@@@@@@@@@@@@@@@@@@@@@@@@@@@@@@@@
%@@@@@@@@@@@@@@@@@@@@@@@@@@@@@@@@@@@@@@@@@@@@@@@@@@@@@@@@@@
%@@@@@@@@@@@@@@@@@@@@@@@@@@@@@@@@@@@@@@@@@@@@@@@@@@@@@@@@@@
%@@@@@@@@@@@@@@@@@@@@@@@@@@@@@@@@@@@@@@@@@@@@@@@@@@@@@@@@@@
%----------------------------------------------------------
%@@@@@@@@@@@@@@@@@@@@@@@@@@@@@@@@@@@@@@@@@@@@@@@@@@@@@@@@@@
%@@@@@@@@@@@@@@@@@@@@@@@@@@@@@@@@@@@@@@@@@@@@@@@@@@@@@@@@@@
%@@@@@@@@@@@@@@@@@@@@@@@@@@@@@@@@@@@@@@@@@@@@@@@@@@@@@@@@@@
%@@@@@@@@@@@@@@@@@@@@@@@@@@@@@@@@@@@@@@@@@@@@@@@@@@@@@@@@@@
%@@@@@@@@@@@@@@@@@@@@@@@@@@@@@@@@@@@@@@@@@@@@@@@@@@@@@@@@@@
%----------------------------------------------------------
\begin{theorem} Let $n<p<\infty, \gamma>1,$ and let $\mu$ be a non-trivial non-negative Borel measure on $\RN$. Let $\Qc$ be a finite family of pairwise disjoint cubes in $\RN$, and let $\Sc$ be a finite family of pairwise disjoint closed subsets of $\RN$. Suppose that to each cube $Q\in\Qc$ we have assigned two subsets $S'_Q,S''_Q\in\Sc$ such that
$S'_Q\cup S''_Q\subset \gamma Q.$
%----------------------------------------------------------
\par If a function $f\in \LOP+\LPM$, then
%----------------------------------------------------------
$$
\shuge_{Q\in\Qc}\,\,
\frac{(\diam Q)^{n-p}\iint \limits_{S'_Q\times S''_Q}
|f(x)-f(y)|^p\, d\mu(x)d\mu(y)}
{ \{(\diam S'_Q)^{n-p}+\mu(S'_Q)\} \{(\diam S''_Q)^{n-p}+\mu(S''_Q)\}}
\le C\|f\|^p_{\sum}\,.
$$
%----------------------------------------------------------
Here $C$ is a constant depending only on $n,p$ and $\gamma$.
%----------------------------------------------------------
\end{theorem}
%----------------------------------------------------------
%@@@@@@@@@@@@@@@@@@@@@@@@@@@@@@@@@@@@@@@@@@@@@@@@@@@@@@@@@@
%@@@@@@@@@@@@@@@@@@@@@@@@@@@@@@@@@@@@@@@@@@@@@@@@@@@@@@@@@@
%@@@@@@@@@@@@@@@@@@@@@@@@@@@@@@@@@@@@@@@@@@@@@@@@@@@@@@@@@@
\par {\it Proof.} By Proposition \reff{OP-W},
%----------------------------------------------------------
$$
\sbig_{Q\in\Qc}\,\,
(\diam Q)^{n-p}\iint \limits_{S'_Q\times S''_Q}
|f(x)-f(y)|^p\, d\mu(x)d\mu(y)/A(S'_Q,S''_Q)
\le C\|f\|^p_{\sum}
$$
%----------------------------------------------------------
where
%----------------------------------------------------------
\bel{AQQP}
~~~~~A(S'_Q,S''_Q):=\mu(S'_Q)\mu(S''_Q)\{1+(\diam S'_Q)^{n-p}/\mu(S'_Q)
+(\diam S''_Q)^{n-p}/\mu(S''_Q)\}\,.
\ee
%----------------------------------------------------------
But
%----------------------------------------------------------
\bel{INAQ}
A(S'_Q,S''_Q)\le\{(\diam S'_Q)^{n-p}+\mu(S'_Q)\} \{(\diam S''_Q)^{n-p}+\mu(S''_Q)\}\nn,
\ee
%----------------------------------------------------------
and the theorem follows.\bx
%----------------------------------------------------------
\bigskip
%@@@@@@@@@@@@@@@@@@@@@@@@@@@@@@@@@@@@@@@@@@@@@@@@@@@@@@@@@@
%----------------------------------------------------------
\par Finally, we apply this theorem to a function $f\in \LOP+\LPM$ with $\Sc=\Qc,S'_Q=Q',$ and $S''_Q=Q''$ proving the necessity part of Theorem \reff{MAIN-CR}.
%----------------------------------------------------------
%@@@@@@@@@@@@@@@@@@@@@@@@@@@@@@@@@@@@@@@@@@@@@@@@@@@@@@@@@@
%@@@@@@@@@@@@@@@@@@@@@@@@@@@@@@@@@@@@@@@@@@@@@@@@@@@@@@@@@@
%@@@@@@@@@@@@@@@@@@@@@@@@@@@@@@@@@@@@@@@@@@@@@@@@@@@@@@@@@@
%@@@@@@@@@@@@@@@@@@@@@@@@@@@@@@@@@@@@@@@@@@@@@@@@@@@@@@@@@@
%@@@@@@@@@@@@@@@@@@@@@@@@@@@@@@@@@@@@@@@@@@@@@@@@@@@@@@@@@@
%----------------------------------------------------------
%@@@@@@@@@@@@@@@@@@@@@@@@@@@@@@@@@@@@@@@@@@@@@@@@@@@@@@@@@@
%@@@@@@@@@@@@@@@@@@@@@@@@@@@@@@@@@@@@@@@@@@@@@@@@@@@@@@@@@@
%@@@@@@@@@@@@@@@@@@@@@@@@@@@@@@@@@@@@@@@@@@@@@@@@@@@@@@@@@@
%@@@@@@@@@@@@@@@@@@@@@@@@@      @@@@@@@@@@@@@@@@@@@@@@@@@@@
%@@@@@@@@@@@@@@@@@@@@@@@          @@@@@@@@@@@@@@@@@@@@@@@@@
%@@@@@@@@@@@@@@@@@@@@@              @@@@@@@@@@@@@@@@@@@@@@@
%@@@@@@@@@@@@@@@@@@@     SECTION 3    @@@@@@@@@@@@@@@@@@@@@
%@@@@@@@@@@@@@@@@@@@@@              @@@@@@@@@@@@@@@@@@@@@@@
%@@@@@@@@@@@@@@@@@@@@@@@          @@@@@@@@@@@@@@@@@@@@@@@@@
%@@@@@@@@@@@@@@@@@@@@@@@@@      @@@@@@@@@@@@@@@@@@@@@@@@@@@
%@@@@@@@@@@@@@@@@@@@@@@@@@@@@@@@@@@@@@@@@@@@@@@@@@@@@@@@@@@
%@@@@@@@@@@@@@@@@@@@@@@@@@@@@@@@@@@@@@@@@@@@@@@@@@@@@@@@@@@
%----------------------------------------------------------
%@@@@@@@@@@@@@@@@@@@@@@@@@@@@@@@@@@@@@@@@@@@@@@@@@@@@@@@@@@
%----------------------------------------------------------
\SECT{3. A $\mu$-measure concentration set and $\mu$-measure concentration cubes.}{3}
%----------------------------------------------------------
\addtocontents{toc}{3. A $\mu$-measure concentration set and $\mu$-measure concentration cubes. \hfill \thepage\\\par}
%----------------------------------------------------------
\indent
%----------------------------------------------------------
\par We turn to the proof of the sufficiency part of Theorem \reff{MAIN-CR}. Actually, in the next three sections  we prove a more general result, Theorem \reff{S-V2}, which immediately implies the sufficiency in Theorem \reff{MAIN-CR}.
%----------------------------------------------------------
%@@@@@@@@@@@@@@@@@@@@@@@@@@@@@@@@@@@@@@@@@@@@@@@@@@@@@@@@@@
%@@@@@@@@@@@@@@@@@@@@@@@@@@@@@@@@@@@@@@@@@@@@@@@@@@@@@@@@@@
%@@@@@@@@@@@@@@@@@@@@@@@@@@@@@@@@@@@@@@@@@@@@@@@@@@@@@@@@@@
%@@@@@@@@@@@@@@@@@@@@@@@@@@@@@@@@@@@@@@@@@@@@@@@@@@@@@@@@@@
%@@@@@@@@@@@@@@@@@@@@@@@@@@@@@@@@@@@@@@@@@@@@@@@@@@@@@@@@@@
%@@@@@@@@@@@@@@@@@@@@@@@@@@@@@@@@@@@@@@@@@@@@@@@@@@@@@@@@@@
%----------------------------------------------------------
\begin{theorem}\lbl{S-V2} Let $n<p<\infty$. A function  $f\in\sum=\LOP+\LPM$ provided $f\in L_{p,loc}(\RN;\mu)$ and there exists a  positive constant $\lambda$ which satisfies the following conditions for a certain absolute positive constant $\gamma$: Let $\Qc$ be an arbitrary finite family of pairwise disjoint cubes in $\RN$ and let $\Qc\ni Q\mapsto Q'\in\Qc$ and $\Qc\ni Q\mapsto Q''\in\Qc$ be arbitrary mappings such that $Q'\cup Q''\subset \gamma Q$ and
%----------------------------------------------------------
\bel{M-DM}
(\diam Q')^{p-n}\mu(Q')+(\diam Q'')^{p-n}\mu(Q'')\le 1.
\ee
%----------------------------------------------------------
Then the following inequality
%----------------------------------------------------------
\bel{CR-V2}
\sbig_{Q\in\Qc}\,\,
\left(\frac{\diam Q' \diam Q''}{\diam Q}\right)^{p-n} \iint \limits_{Q'\times Q''}
|f(x)-f(y)|^p\, d\mu(x)d\mu(y)
\le \lambda
\ee
%----------------------------------------------------------
holds. Furthermore, $\|f\|_{\sum}\le C(n,p)\, \lambda^{\frac{1}{p}}.$
%----------------------------------------------------------
\end{theorem}
%----------------------------------------------------------
%@@@@@@@@@@@@@@@@@@@@@@@@@@@@@@@@@@@@@@@@@@@@@@@@@@@@@@@@@@
%@@@@@@@@@@@@@@@@@@@@@@@@@@@@@@@@@@@@@@@@@@@@@@@@@@@@@@@@@@
%@@@@@@@@@@@@@@@@@@@@@@@@@@@@@@@@@@@@@@@@@@@@@@@@@@@@@@@@@@
%@@@@@@@@@@@@@@@@@@@@@@@@@@@@@@@@@@@@@@@@@@@@@@@@@@@@@@@@@@
%@@@@@@@@@@@@@@@@@@@@@@@@@@@@@@@@@@@@@@@@@@@@@@@@@@@@@@@@@@
%@@@@@@@@@@@@@@@@@@@@@@@@@@@@@@@@@@@@@@@@@@@@@@@@@@@@@@@@@@
%----------------------------------------------------------
\begin{remark} {\em We first prove a version of Theorem \reff{S-V2} where inequality \rf{M-DM} is  replaced by weaker conditions
%----------------------------------------------------------
\bel{G-MD}
\mu(Q')\le 2^{32p}(\diam Q')^{n-p}~~~\text{and}~~~
\mu(Q'')\le 2^{32p}(\diam Q'')^{n-p}.
\ee
%----------------------------------------------------------
It can be readily seen that Theorem \reff{S-V2} in its original formulation immediately follows from this weaker variant by transition to the measure $\tilde{\mu}=2^{32p}\mu$.
%----------------------------------------------------------
\par Thus throughout the proof of the theorem we will assume that the cubes $Q',Q''$ satisfy inequalities \rf{G-MD} rather than \rf{M-DM}\rbx}.
\end{remark}
%----------------------------------------------------------
%@@@@@@@@@@@@@@@@@@@@@@@@@@@@@@@@@@@@@@@@@@@@@@@@@@@@@@@@@@
%@@@@@@@@@@@@@@@@@@@@@@@@@@@@@@@@@@@@@@@@@@@@@@@@@@@@@@@@@@
%@@@@@@@@@@@@@@@@@@@@@@@@@@@@@@@@@@@@@@@@@@@@@@@@@@@@@@@@@@
%@@@@@@@@@@@@@@@@@@@@@@@@@@@@@@@@@@@@@@@@@@@@@@@@@@@@@@@@@@
%@@@@@@@@@@@@@@@@@@@@@@@@@@@@@@@@@@@@@@@@@@@@@@@@@@@@@@@@@@
%@@@@@@@@@@@@@@@@@@@@@@@@@@@@@@@@@@@@@@@@@@@@@@@@@@@@@@@@@@
%----------------------------------------------------------
\par Let $f$ be a function on $\RN$ satisfying the hypothesis of Theorem \reff{S-V2}. Let us construct its almost optimal decomposition, i.e., functions $f_1\in\LOP$ and $f_2\in \LPM$ with almost minimal norms in the spaces $\LOP$ and $\LPM$ respectively, and such that $f=f_1+f_2$.
%----------------------------------------------------------
\par We do this in two stages. At the first stage which we present in this section we construct a closed set $E\subset\RN$ and a family $\Kc_E$ of pairwise disjoint ``well separated'' cubes of $\RN$ with centers in $E$, see Proposition \reff{RPR} and definition \rf{KE-DF}. These cubes are determined only by the measure $\mu$, and by $n$ and $p$. The set $E$ and the family $\Kc_E$ possess certain $\mu$-measure concentration properties which we prove in Corollary \reff{PR-5K} and Lemma \reff{Q-NE}. In Section 4, using the averages of the function $f$ on cubes from $\Kc_E$ and the Whitney extension method, we define the function $f_1$. See \rf{DEF-TF1} and \rf{DEF-F1}. Finally, we put $f_2:=f-f_1$. %----------------------------------------------------------
%@@@@@@@@@@@@@@@@@@@@@@@@@@@@@@@@@@@@@@@@@@@@@@@@@@@@@@@@@@
%@@@@@@@@@@@@@@@@@@@@@@@@@@@@@@@@@@@@@@@@@@@@@@@@@@@@@@@@@@
%@@@@@@@@@@@@@@@@@@@@@@@@@@@@@@@@@@@@@@@@@@@@@@@@@@@@@@@@@@
%@@@@@@@@@@@@@@@@@@@@@@@@@@@@@@@@@@@@@@@@@@@@@@@@@@@@@@@@@@
%----------------------------------------------------------
\par We begin the first stage with the following theorem which is an important element of our geometrical construction. %----------------------------------------------------------
%@@@@@@@@@@@@@@@@@@@@@@@@@@@@@@@@@@@@@@@@@@@@@@@@@@@@@@@@@@
%@@@@@@@@@@@@@@@@@@@@@@@@@@@@@@@@@@@@@@@@@@@@@@@@@@@@@@@@@@
%@@@@@@@@@@@@@@@@@@@@@@@@@@@@@@@@@@@@@@@@@@@@@@@@@@@@@@@@@@
%@@@@@@@@@@@@@@@@@@@@@@@@@@@@@@@@@@@@@@@@@@@@@@@@@@@@@@@@@@
%@@@@@@@@@@@@@@@@@@@@@@@@@@@@@@@@@@@@@@@@@@@@@@@@@@@@@@@@@@
%@@@@@@@@@@@@@@@@@@@@@@@@@@@@@@@@@@@@@@@@@@@@@@@@@@@@@@@@@@
%@@@@@@@@@@@@@@@@@@@@@@@@@@@@@@@@@@@@@@@@@@@@@@@@@@@@@@@@@@
%----------------------------------------------------------
\begin{theorem}\lbl{VRN} Let $w:\RN\to(0,\infty)$ be a positive function on $\RN$ such that for every $x\in\RN$ the following inequality
%----------------------------------------------------------
\bel{AZ}
\liminf_{y\to x}w(y)>0
\ee
%----------------------------------------------------------
holds. Then there exists a set $S\subset \RN$ which satisfies all of the following conditions:
%----------------------------------------------------------
\par (i). For every $x\in\RN$ there exists a point $\tx\in S$ such that
%----------------------------------------------------------
\bel{BPS}
\|x-\tx\|+w(\tx)\le 83\,w(x);
\ee
%----------------------------------------------------------
%@@@@@@@@@@@@@@@@@@@@@@@@@@@@@@@@@@@@@@@@@@@@@@@@@@@@@@@@@@
%----------------------------------------------------------
\par (ii). For every $z_1,z_2\in S, z_1\neq z_2,$ we have
%----------------------------------------------------------
$$
w(z_1)+w(z_2)\le \|z_1-z_2\|/6.
$$
%----------------------------------------------------------
\end{theorem}
%----------------------------------------------------------
%@@@@@@@@@@@@@@@@@@@@@@@@@@@@@@@@@@@@@@@@@@@@@@@@@@@@@@@@@@
%@@@@@@@@@@@@@@@@@@@@@@@@@@@@@@@@@@@@@@@@@@@@@@@@@@@@@@@@@@
%@@@@@@@@@@@@@@@@@@@@@@@@@@@@@@@@@@@@@@@@@@@@@@@@@@@@@@@@@@
\par {\it Proof.} Given an integer $j$ we define a set
%---------------------------------------------------
\bel{DAJ}
A_{j}:=\{y \in\RN:2^{-j-1} < w(y) \le 2^{-j}\}.
\ee
%---------------------------------------------------
We introduce a metric on $\RN$ by letting
%----------------------------------------------------------
$$
\row(x,y):=\left \{
%----------------------------------------------------------
\begin{array}{ll}
\|x-y\|+w(x)+w(y),& x\ne y,\\
\\
0,& x=y.
\end{array}
%----------------------------------------------------------
\right.
$$
%----------------------------------------------------------
%@@@@@@@@@@@@@@@@@@@@@@@@@@@@@@@@@@@@@@@@@@@@@@@@@@@@@@@@@@
Let $\ve_j:=14\cdot 2^{-j}$ and let $B_j$ be a maximal $\ve_j$-net in $A_j$ with respect to the metric $\row$. Thus, if $A_j\ne\emp$ and $\card B_j>1$, the following conditions are satisfied:
%@@@@@@@@@@@@@@@@@@@@@@@@@@@@@@@@@@@@@@@@@@@@@@@@@@@@@@@@@@
\par (1). For every $z_1,z_2\in B_j, z_1\ne z_2,$ we have
%----------------------------------------------------------
$$
\row(z_1,z_2)\ge \ve_j~;
$$
%----------------------------------------------------------
%@@@@@@@@@@@@@@@@@@@@@@@@@@@@@@@@@@@@@@@@@@@@@@@@@@@@@@@@@@
\par (2). For every $x\in A_j$ there exists a point $x'\in B_j$ such that
%----------------------------------------------------------
$$
\row(x,x')<\ve_j~;
$$
%----------------------------------------------------------
Since
%----------------------------------------------------------
$$
14w(x)\le 14\cdot 2^{-j}=\ve_j<28w(x)~~~\text{for every}~~~z\in A_j\,,
$$
%----------------------------------------------------------
we have
%----------------------------------------------------------
\bel{IW}
\row(z_1,z_2)\ge 7\{w(z_1)+w(z_2)\}, ~~~z_1,z_2\in B_j\,,
\ee
%----------------------------------------------------------
and
%----------------------------------------------------------
\bel{CP}
\row(x,x')< 28w(x)
\ee
%----------------------------------------------------------
for some $x'\in B_j$.
%----------------------------------------------------------
\par Given $\ve>0$ and a set $B\in\RN$ we let $[B]_{\ve}$ denote the closed $\ve$-neighborhood of $B$ with respect to the metric $\row$:
%----------------------------------------------------------
$$
[B]_{\ve}:=\{x\in \RN:\exists\,\, y\in B~~\text{such that}~~\row(x,y)\le \ve\}.
$$
%----------------------------------------------------------
Let us define a set $\tB_j$ by letting
%----------------------------------------------------------
\bel{DBT}
\tB_j:=B_j\setminus \left[\bigcup_{i>j}B_i\right]_{\ve_j}.
\ee
%----------------------------------------------------------
Finally we put
%----------------------------------------------------------
\bel{DSF}
S:=\bigcup_{j=-\infty}^\infty\tB_j.
\ee
%----------------------------------------------------------
\par Prove that $S$ satisfies all the conditions of the proposition. We do this in three steps.
%----------------------------------------------------------
%@@@@@@@@@@@@@@@@@@@@@@@@@@@@@@@@@@@@@@@@@@@@@@@@@@@@@@@@@@
%@@@@@@@@@@@@@@@@@@@@@@@@@@@@@@@@@@@@@@@@@@@@@@@@@@@@@@@@@@
%@@@@@@@@@@@@@@@@@@@@@@@@@@@@@@@@@@@@@@@@@@@@@@@@@@@@@@@@@@
%----------------------------------------------------------
\medskip
%----------------------------------------------------------
\par {\it The first step.} Prove that $S\ne\emp.$
%----------------------------------------------------------
\par Suppose that $S=\emp$ and prove that this contradicts to the condition \rf{AZ}. Since $\bigcup\limits_{j=-\infty}^\infty A_j=\RN$, there exists $j_0\in\mZ$ such that $A_{j_0}\ne\emp.$ Hence $B_{j_0}\ne\emp$ as well so that there exists a point $x_0\in B_{j_0}$.
%----------------------------------------------------------
\par By the assumption
%----------------------------------------------------------
$$
S:=\bigcup_{j=-\infty}^\infty\tB_j=\emp,
$$
%----------------------------------------------------------
so that the set
%----------------------------------------------------------
$$
\tB_{j_0}:=B_{j_0}\setminus \left[\bigcup_{i>j_0}B_i\right]_{\ve_{j_0}}=\emp.
$$
%----------------------------------------------------------
Therefore there exist an integer $j_1>j_0$ and a point $x_1\in B_{j_1}$ such that
%----------------------------------------------------------
$$
\row(x_0,x_1)\le \ve_{j_0}.
$$
%----------------------------------------------------------
But $\tB_{j_1}=\emp$ as well so that there exist an integer $j_2>j_1$ and a point $x_2\in B_{j_2}$ such that
%----------------------------------------------------------
$$
\row(x_1,x_2)\le \ve_{j_1}.
$$
%----------------------------------------------------------
Continuing this process we get a sequence of points $\{x_k\}_{k=0}^\infty$ such that
%----------------------------------------------------------
$$
\row(x_k,x_{k+1}):=\|x_k-x_{k+1}\|+w(x_k)+w(x_{k+1})\le \ve_{j_k}:=14\cdot 2^{-j_k},~~~k=0,1,...~.
$$
%----------------------------------------------------------
Hence $w(x_k)\le \ve_{j_k}$ so that
$w(x_k)\to\ 0$ as $k\to \infty.$ Furthermore, since
%----------------------------------------------------------
$$
\|x_k-x_{k+1}\|\le 14\cdot 2^{-j_k},~~~k=0,1,...~,
$$
%----------------------------------------------------------
$\{x_k\}_{k=0}^\infty$ is a Cauchy sequence so that there exist $\bar{x}\in\RN$ such that
%----------------------------------------------------------
$$
\lim_{k\to\infty}x_k=\bar{x}.
$$
%----------------------------------------------------------
Since $w\ge 0$, we obtain
%----------------------------------------------------------
$$
\liminf_{x\to\bar{x}}w(x_k)=0,
$$
%----------------------------------------------------------
a contradiction.
%----------------------------------------------------------
\medskip
%@@@@@@@@@@@@@@@@@@@@@@@@@@@@@@@@@@@@@@@@@@@@@@@@@@@@@@@@@@
%@@@@@@@@@@@@@@@@@@@@@@@@@@@@@@@@@@@@@@@@@@@@@@@@@@@@@@@@@@
%@@@@@@@@@@@@@@@@@@@@@@@@@@@@@@@@@@@@@@@@@@@@@@@@@@@@@@@@@@
%----------------------------------------------------------
\par {\it The second step.} Prove the property (ii) of the proposition which is equivalent to the inequality
%----------------------------------------------------------
\bel{STT}
\row(z_1,z_2)\ge 7\{w(z_1)+w(z_2)\}.
\ee
%----------------------------------------------------------
\par Suppose that $z_1\ne z_2$ and $w(z_2)\le w(z_1)$. If
$z_1,z_2\in \tB_j$ for some integer $j$, then \rf{STT} follows from \rf{IW}. Suppose that $z_1\in \tB_i,z_2\in \tB_j$ for some $i>j$. Since $z_2\in \tB_j$, by \rf{DBT}, $z_2\notin [\tB_i]_{\ve_j}$ so that
%----------------------------------------------------------
$$
\row(z_1,z_2)\ge \ve_j=14\cdot 2^{-j}.
$$
%----------------------------------------------------------
On the other hand, since  $z_1\in \tB_i\subset A_i,z_2\in \tB_j\subset A_j$, by \rf{DAJ},
%----------------------------------------------------------
$$
w(z_1)\le 2^{-i},~~w(z_2)\le 2^{-j}.
$$
%----------------------------------------------------------
Hence
%----------------------------------------------------------
$$
\row(z_1,z_2)\ge 14\cdot 2^{-j}\ge 7(2^{-j}+2^{-i})\ge 7\{w(z_1)+w(z_2)\}
$$
%----------------------------------------------------------
proving \rf{STT}.
%----------------------------------------------------------
\medskip
%@@@@@@@@@@@@@@@@@@@@@@@@@@@@@@@@@@@@@@@@@@@@@@@@@@@@@@@@@@
%@@@@@@@@@@@@@@@@@@@@@@@@@@@@@@@@@@@@@@@@@@@@@@@@@@@@@@@@@@
%@@@@@@@@@@@@@@@@@@@@@@@@@@@@@@@@@@@@@@@@@@@@@@@@@@@@@@@@@@
%----------------------------------------------------------
\par {\it The third step.} Prove the property (i) of the proposition. Clearly, if $x\in S$, then we can put $\tx:=x$.
%----------------------------------------------------------
\par Let $x\in\RN\setminus S$. Put  $x_0:=x$. Then there exist an integer $j_0$ such that $x\in A_{j_0}$ so that
$x\in A_{j_0}\setminus S$. By \rf{CP}, there exist a point $x_0\in B_{j_0}$ such that
%----------------------------------------------------------
\bel{RDX}
\row(x_0,x)< 28w(x).
\ee
%----------------------------------------------------------
If $x_0\in \tB_{j_0}$, see \rf{DBT}, then we put $\tx:=x_0$ and stop. If $x_0\notin\tB_{j_0}$, then, by \rf{DSF} and \rf{DBT}, there exist an integer $j_1> j_0$ and a point $x_1\in B_{j_1}$ such that
%----------------------------------------------------------
$$
\row(x_0,x_1)\le \ve_{j_0}.
$$
%----------------------------------------------------------
If $x_1\in \tB_{j_1}$, then we put $\tx:=x_1$ and stop. If $x_1\notin\tB_{j_1}$, then there exist an integer $j_2> j_1$ and a point $x_2\in B_{j_2}$ such that
$\row(x_1,x_2)\le \ve_{j_1}$.
%----------------------------------------------------------
\par We continue this process and after $k+1$ stages of the procedure we obtain $k+1$ integers $j_0<j_1<...<j_k$ and points $x_m\in B_{j_m}$, $m=0,...,k$, such that
%----------------------------------------------------------
\bel{RK}
\row(x_m,x_{m+1})\le \ve_{j_m},~~~m=0,...,k-1.
\ee
%----------------------------------------------------------
If $x_k\in \tB_{j_k}$, then we put $\tx:=x_k$ and stop. If $x_k\notin\tB_{j_k}$, then, by \rf{DSF} and \rf{DBT}, there exist an integer $j_{k+1}> j_k$ and a point $x_{k+1}\in B_{j_{k+1}}$ such that
$\row(x_k,x_{k+1})\le \ve_{j_{k}}$.
%----------------------------------------------------------
\par Let us prove that this procedure is finite, i.e., $x_k\in \tB_{j_k}$ for some $k\ge 1$. In fact, otherwise there exists an infinite sequence of points $\{x_m\}_{m=0}^\infty$ such that $x_m\in B_{j_m}$, $m=0,1...,$ and
%----------------------------------------------------------
$$
\row(x_m,x_{m+1}):=\|x_m-x_{m+1}\|+w(x_m)+w(x_{m+1})
\le\ve_{j_m}=14\cdot 2^{-j_m}.
$$
%----------------------------------------------------------
Hence
%----------------------------------------------------------
$$
0\le w(x_m)\le 14\cdot 2^{-j_m}.
$$
%----------------------------------------------------------
so that $w(x_m)\to 0$ as $k\to\infty$. Furthermore,
%----------------------------------------------------------
$$
\|x_m-x_{m+1}\|\le 14\cdot 2^{-j_m}
$$
%----------------------------------------------------------
so that $\{x_m\}_{m=0}^\infty$ is a Cauchy sequence. Consequently  $\{x_m\}_{m=0}^\infty$ converges to a point $\bar{x}\in\RN$. Hence
%----------------------------------------------------------
$$
\liminf_{x\to\bar{x}}w(x)=0
$$
%----------------------------------------------------------
which contradicts \rf{AZ}.
%----------------------------------------------------------
\par Thus we have proved that there exists a positive integer $k$ such that for all $x_k\in S$ and for all $m=0,...,k$ inequality \rf{RK} is satisfied. We have
%----------------------------------------------------------
$$
\row(x,x_{k})\le \row(x,x_{0})+\sum_{m=0}^{k-1}\row(x_m,x_{m+1})
$$
%----------------------------------------------------------
so that, by \rf{RDX} and \rf{RK},
%----------------------------------------------------------
$$
\row(x,x_{k})\le 28w(x)+\sum_{m=0}^{k-1}\ve_{j_m}\le
28w(x)+14\,\sum_{j\ge j_0} 2^{-j_m}\le 28w(x)+28\cdot 2^{-j_0}.
$$
%----------------------------------------------------------
\par Recall that $x\in A_{j_0}$ so that $w(x)\ge 2^{-j_0-1}$. Hence
%----------------------------------------------------------
$$
\row(x,x_{k})\le 28w(x)+28\cdot(2w(x))= 84w(x).
$$
%----------------------------------------------------------
Thus the point $\tx:=x_k\in S$ and the following inequality %----------------------------------------------------------
$$
\row(x,\tx):=\|x-\tx\|+w(x)+w(\tx)\le 84w(x)
$$
%----------------------------------------------------------
holds. This proves inequality \rf{BPS} and the theorem.\bx
%----------------------------------------------------------
\bigskip
%----------------------------------------------------------
%@@@@@@@@@@@@@@@@@@@@@@@@@@@@@@@@@@@@@@@@@@@@@@@@@@@@@@@@@@
%@@@@@@@@@@@@@@@@@@@@@@@@@@@@@@@@@@@@@@@@@@@@@@@@@@@@@@@@@@
%@@@@@@@@@@@@@@@@@@@@@@@@@@@@@@@@@@@@@@@@@@@@@@@@@@@@@@@@@@
%@@@@@@@@@@@@@@@@@@@@@@@@@@@@@@@@@@@@@@@@@@@@@@@@@@@@@@@@@@
%@@@@@@@@@@@@@@@@@@@@@@@@@@@@@@@@@@@@@@@@@@@@@@@@@@@@@@@@@@
%@@@@@@@@@@@@@@@@@@@@@@@@@@@@@@@@@@@@@@@@@@@@@@@@@@@@@@@@@@
%----------------------------------------------------------
\par Note that Theorem \reff{VRN} can be reformulated in a purely geometrical way. We discuss a geometrical background of this theorem and other related geometrical problems in Subsection 7.3.
%----------------------------------------------------------
\par Fix a point $x\in\RN$ and consider two functions of a positive parameter $r$: a function
%----------------------------------------------------------
$$
s_x(r):=\mu(Q(x,r)),~~~r\in(0,+\infty),
$$
%---------------------------------------------------------- %@@@@@@@@@@@@@@@@@@@@@@@@@@@@@@@@@@@@@@@@@@@@@@@@@@@@@@@@@@
and a function
%----------------------------------------------------------
$$
v(r):=r^{n-p},~~~r\in(0,+\infty).
$$
%----------------------------------------------------------
%@@@@@@@@@@@@@@@@@@@@@@@@@@@@@@@@@@@@@@@@@@@@@@@@@@@@@@@@@@
%@@@@@@@@@@@@@@@@@@@@@@@@@@@@@@@@@@@@@@@@@@@@@@@@@@@@@@@@@@
Clearly, $s=s_x(r)$ is a {\it non-decreasing} function on $(0,+\infty)$. Since $\mu$ is a non-trivial measure on $\RN$,
%----------------------------------------------------------
$$
\lim\limits_{r\to +\infty}s_x(r)=\mu(\RN)>0.
$$
%---------------------------------------------------------- %@@@@@@@@@@@@@@@@@@@@@@@@@@@@@@@@@@@@@@@@@@@@@@@@@@@@@@@@@@
\par On the other hand, since $p>n$, the function $v=v(r)$
is {\it strictly decreasing} on $(0,+\infty)$. Clearly,
%----------------------------------------------------------
$$
\lim\limits_{r\to 0} v(r)=+\infty~~~\text{and}~~~\lim\limits_{r\to +\infty}v(r)=0.
$$
%---------------------------------------------------------- %@@@@@@@@@@@@@@@@@@@@@@@@@@@@@@@@@@@@@@@@@@@@@@@@@@@@@@@@@@
\par These properties of the functions $s_x(r)$ and $v(r)$ imply  {\it the existence of a unique number} $R(x)\in (0,\infty)$ such that
%----------------------------------------------------------
$$
s_x(r)>v(R(x))~~\text{if}~~r>R(x),~~~\text{and}~~~
s_x(r)<v(R(x)),~~\text{if}~~r<R(x).
$$
%---------------------------------------------------------- %@@@@@@@@@@@@@@@@@@@@@@@@@@@@@@@@@@@@@@@@@@@@@@@@@@@@@@@@@@
\par Thus, for every $x\in\RN$ we have
%----------------------------------------------------------
\bel{M-RX1}
\mu(Q(x,r))>R(x)^{n-p}~~~\text{for every}~~r>R(x),
\ee
%---------------------------------------------------------- %@@@@@@@@@@@@@@@@@@@@@@@@@@@@@@@@@@@@@@@@@@@@@@@@@@@@@@@@@@
and
%----------------------------------------------------------
\bel{M-RX2}
\mu(Q(x,r))<R(x)^{n-p}~~~\text{for every}~~r<R(x).
\ee
%---------------------------------------------------------- %@@@@@@@@@@@@@@@@@@@@@@@@@@@@@@@@@@@@@@@@@@@@@@@@@@@@@@@@@@
Since $\mu$ is a Borel measure, the function $s_x(r):=\mu(Q(x,r))$ is right continuous on $(0,\infty)$ so that, by \rf{M-RX1},
%----------------------------------------------------------
\bel{M-RX3}
\mu(Q(x,R(x)))\ge R(x)^{n-p}.
\ee
%----------------------------------------------------------  %@@@@@@@@@@@@@@@@@@@@@@@@@@@@@@@@@@@@@@@@@@@@@@@@@@@@@@@@@@
\par We also recall that the number $R(x)$ satisfies the inequality
%----------------------------------------------------------
\bel{RXP}
0<R(x)<+\infty.
\ee
%---------------------------------------------------------- %@@@@@@@@@@@@@@@@@@@@@@@@@@@@@@@@@@@@@@@@@@@@@@@@@@@@@@@@@@
%@@@@@@@@@@@@@@@@@@@@@@@@@@@@@@@@@@@@@@@@@@@@@@@@@@@@@@@@@@
%@@@@@@@@@@@@@@@@@@@@@@@@@@@@@@@@@@@@@@@@@@@@@@@@@@@@@@@@@@
\begin{lemma}\lbl{PRC} The function $R=R(x)$ satisfies the Lipschitz condition on $\RN$:
%----------------------------------------------------------
\bel{LIPR}
|R(x)-R(y)|\le\|x-y\|~~~~\text{for every}~~~~x,y\in\RN.
\ee
%----------------------------------------------------------
\end{lemma}
%----------------------------------------------------------
%@@@@@@@@@@@@@@@@@@@@@@@@@@@@@@@@@@@@@@@@@@@@@@@@@@@@@@@@@@
%@@@@@@@@@@@@@@@@@@@@@@@@@@@@@@@@@@@@@@@@@@@@@@@@@@@@@@@@@@
%@@@@@@@@@@@@@@@@@@@@@@@@@@@@@@@@@@@@@@@@@@@@@@@@@@@@@@@@@@
\par {\it Proof.} Suppose that $R(x)>R(y)$.
%@@@@@@@@@@@@@@@@@@@@@@@@@@@@@@@@@@@@@@@@@@@@@@@@@@@@@@@@@@
\par Let $r\in(R(y),R(x))$. Then, by \rf{M-RX2},
%----------------------------------------------------------
\bel{QRP}
\mu(Q(x,r))<r^{n-p}.
\ee
%---------------------------------------------------------- %@@@@@@@@@@@@@@@@@@@@@@@@@@@@@@@@@@@@@@@@@@@@@@@@@@@@@@@@@@
\par Prove that $Q(y,R(y)) \varsubsetneq Q(x,r)$. In fact, if $Q(y,R(y))\subset Q(x,r)$, then, by \rf{M-RX3},
%----------------------------------------------------------
$$
\mu(Q(x,r))\ge \mu(Q(y,R(y)))\ge R(y)^{n-p}\ge r^{n-p}
$$
%----------------------------------------------------------
which contradicts \rf{QRP}.
%@@@@@@@@@@@@@@@@@@@@@@@@@@@@@@@@@@@@@@@@@@@@@@@@@@@@@@@@@@
\par Thus for every $r\in(R(y),R(x))$ there exists a point $a_r\in Q(y,R(y))\setminus Q(x,r)$ so that $\|a_r-x\|>r$ and $\|a_r-y\|\le R(y)$. Hence
%----------------------------------------------------------
$$
r<\|a_r-x\|\le \|a_r-y\|+\|y-x\|\le R(y)+\|x-y\|
$$
%----------------------------------------------------------
proving that
%----------------------------------------------------------
$$
|r-R(y)|=r-R(y)<\|x-y\|.
$$
%----------------------------------------------------------
Since $r\in(R(y),R(x))$ is arbitrary, we obtain the required inequality \rf{LIPR}.\bx\smallskip
%----------------------------------------------------------
%@@@@@@@@@@@@@@@@@@@@@@@@@@@@@@@@@@@@@@@@@@@@@@@@@@@@@@@@@@
%@@@@@@@@@@@@@@@@@@@@@@@@@@@@@@@@@@@@@@@@@@@@@@@@@@@@@@@@@@
%@@@@@@@@@@@@@@@@@@@@@@@@@@@@@@@@@@@@@@@@@@@@@@@@@@@@@@@@@@
%----------------------------------------------------------
\begin{proposition}\lbl{RPR} There exists a subset $E\subset \RN$ such that:
%@@@@@@@@@@@@@@@@@@@@@@@@@@@@@@@@@@@@@@@@@@@@@@@@@@@@@@@@@@
\par (i). For every $x,y\in E, x\ne y,$
%----------------------------------------------------------
\bel{DRS}
6(R(x)+R(y))\le\|x-y\|;
\ee
%----------------------------------------------------------
%@@@@@@@@@@@@@@@@@@@@@@@@@@@@@@@@@@@@@@@@@@@@@@@@@@@@@@@@@@
\par (ii). For every $x\in\RN$ there exists a point
$\tx\in E$ such that
%----------------------------------------------------------
\bel{IRW}
R(\tx)\le 83\,R(x)
\ee
%----------------------------------------------------------
and
%----------------------------------------------------------
\bel{IDW}
\|\tx-x\|\le 83\,R(x).
\ee
%----------------------------------------------------------
\end{proposition}
%----------------------------------------------------------
%@@@@@@@@@@@@@@@@@@@@@@@@@@@@@@@@@@@@@@@@@@@@@@@@@@@@@@@@@@
%@@@@@@@@@@@@@@@@@@@@@@@@@@@@@@@@@@@@@@@@@@@@@@@@@@@@@@@@@@
%@@@@@@@@@@@@@@@@@@@@@@@@@@@@@@@@@@@@@@@@@@@@@@@@@@@@@@@@@@
\par {\it Proof.} By Lemma \reff{PRC}, the function $R=R(x)$ is Lipschitz continuous on $\RN$ so that it is continuous. Hence for every $x\in\RN$ we have
%----------------------------------------------------------
$$
\liminf_{y\to x} R(y)=\lim_{y\to x} R(y)=R(x).
$$
%----------------------------------------------------------
Since $R(x)>0$ (see \rf{RXP}), condition \rf{AZ} of Theorem \reff{VRN} for the function $w(x):=R(x), x\in\RN,$ is satisfied. By this theorem, there exists a set $E\subset\RN$ satisfying the required inequalities \rf{DRS},\rf{IRW}, and \rf{IDW}.\bx
%----------------------------------------------------------
\bigskip
%----------------------------------------------------------
%@@@@@@@@@@@@@@@@@@@@@@@@@@@@@@@@@@@@@@@@@@@@@@@@@@@@@@@@@@
%@@@@@@@@@@@@@@@@@@@@@@@@@@@@@@@@@@@@@@@@@@@@@@@@@@@@@@@@@@
%@@@@@@@@@@@@@@@@@@@@@@@@@@@@@@@@@@@@@@@@@@@@@@@@@@@@@@@@@@
%----------------------------------------------------------
\par Note that, by inequality \rf{DRS}, the set $E$ consists of {\it isolated points} of\, $\RN$.  %@@@@@@@@@@@@@@@@@@@@@@@@@@@@@@@@@@@@@@@@@@@@@@@@@@@@@@@@@@
%----------------------------------------------------------
\par Given $x\in E$ we let $K^{(x)}$ denote the cube
%----------------------------------------------------------
\bel{D-KX}
K^{(x)}:=Q(x,R(x)).
\ee
%----------------------------------------------------------
We put
%----------------------------------------------------------
\bel{KE-DF}
\Kc_E:=\{K^{(x)}:x\in E\}.
\ee
%----------------------------------------------------------
\par Recall that, by the inequality \rf{M-RX3}, for every $K\in \Kc_E$ we have
%----------------------------------------------------------
\bel{KXM}
\mu(K)\ge 2^{p-n}(\diam K)^{n-p}.
\ee
%----------------------------------------------------------
\par Let us present several properties of the set $E$ and the cubes of the family  $\Kc_E$.
%---------------------------------------------------------- %@@@@@@@@@@@@@@@@@@@@@@@@@@@@@@@@@@@@@@@@@@@@@@@@@@@@@@@@@@
%@@@@@@@@@@@@@@@@@@@@@@@@@@@@@@@@@@@@@@@@@@@@@@@@@@@@@@@@@@
%@@@@@@@@@@@@@@@@@@@@@@@@@@@@@@@@@@@@@@@@@@@@@@@@@@@@@@@@@@
\begin{lemma}\lbl{PR-KE} (i). For every two cubes $K,K'\in \Kc_E$, $K\ne K'$, we have
%----------------------------------------------------------
$$
\diam K +\diam K'\le \dist(K,K')/2;
$$
%----------------------------------------------------------
\par (ii). Let $\tau\ge 1$ be a constant and let $x,x'\in E$, $x\ne x'$. Let $Q,Q'$ be cubes in $\RN$ such that $\tau Q\ni x, \tau Q'\ni x'$, and $Q\cap Q'\ne\emp$. Then
%----------------------------------------------------------
$$
\diam K^{(x)} +\diam K^{(x')}\le \tau(\diam Q +\diam Q').
$$
%----------------------------------------------------------
\end{lemma}
%----------------------------------------------------------
%@@@@@@@@@@@@@@@@@@@@@@@@@@@@@@@@@@@@@@@@@@@@@@@@@@@@@@@@@@
%@@@@@@@@@@@@@@@@@@@@@@@@@@@@@@@@@@@@@@@@@@@@@@@@@@@@@@@@@@
%@@@@@@@@@@@@@@@@@@@@@@@@@@@@@@@@@@@@@@@@@@@@@@@@@@@@@@@@@@
\par {\it Proof.} (i). Let $K=K^{(a)},K'=K^{(a')}$ for some $a,a'\in E$. By part (i) of Proposition \reff{RPR},
%----------------------------------------------------------
$$
6(R(a)+R(a'))\le \|a-a'\|.
$$
%----------------------------------------------------------
On the other hand
%----------------------------------------------------------
$$
\|a-a'\|\le\dist(K,K')+R(a)+R(a').
$$
%----------------------------------------------------------
Hence
%----------------------------------------------------------
\be
\dist(K,K')&\ge&\|a-a'\|-R(a)-R(a')\nn\\
&\ge& 6(R(a)+R(a'))-R(a)-R(a')=5(R(a)+R(a'))\nn\\
&=&\tfrac52(\diam K+\diam K')\nn
\ee
%----------------------------------------------------------
proving the statement (i).
%----------------------------------------------------------
\par (ii). By part (i) of Proposition \reff{RPR},
%----------------------------------------------------------
$$
3(\diam K^{(x)} +\diam K^{(x')})\le \|x-x'\|.
$$
%----------------------------------------------------------
Since $\tau Q\ni x, \tau Q'\ni x'$, and $Q\cap Q'\ne\emp$,
%----------------------------------------------------------
$$
\|x-x'\|\le\tau r_Q+\tau r_{Q'}+r_Q+r_{Q'}=
(\tau+1)( r_Q+ r_{Q'})\le 2\tau( r_Q+r_{Q'})
$$
%----------------------------------------------------------
so that
%----------------------------------------------------------
$$
\diam K^{(x)} +\diam K^{(x')}\le \tfrac13\|x-x'\|\le
\tfrac23\tau(\diam Q +\diam Q').
$$
%----------------------------------------------------------
\par The lemma is proved.\bx\medskip
%----------------------------------------------------------
%@@@@@@@@@@@@@@@@@@@@@@@@@@@@@@@@@@@@@@@@@@@@@@@@@@@@@@@@@@
%@@@@@@@@@@@@@@@@@@@@@@@@@@@@@@@@@@@@@@@@@@@@@@@@@@@@@@@@@@
%@@@@@@@@@@@@@@@@@@@@@@@@@@@@@@@@@@@@@@@@@@@@@@@@@@@@@@@@@@
%@@@@@@@@@@@@@@@@@@@@@@@@@@@@@@@@@@@@@@@@@@@@@@@@@@@@@@@@@@
%----------------------------------------------------------
\begin{lemma}\lbl{Q-NE} For every cube $Q\subset\RN$ and every $\theta>0$ such that  %----------------------------------------------------------
$$
\diam Q\le\theta\,\dist(Q,E)
$$
%----------------------------------------------------------
the following inequality
%----------------------------------------------------------
$$
\mu(Q)\le 42^p(1+\theta)^p\, r_Q^{n-p}
$$
%----------------------------------------------------------
holds.
%----------------------------------------------------------
\end{lemma}
%----------------------------------------------------------
%@@@@@@@@@@@@@@@@@@@@@@@@@@@@@@@@@@@@@@@@@@@@@@@@@@@@@@@@@@
%@@@@@@@@@@@@@@@@@@@@@@@@@@@@@@@@@@@@@@@@@@@@@@@@@@@@@@@@@@
%@@@@@@@@@@@@@@@@@@@@@@@@@@@@@@@@@@@@@@@@@@@@@@@@@@@@@@@@@@
\par {\it Proof.} By Proposition \reff{RPR}, for every $x\in \RN$ there exists a point $\tilde{x}\in E$ such that
%----------------------------------------------------------
$$
\|\tilde{x}-x\|\le 83\, R(x).
$$
%----------------------------------------------------------
Hence,
%----------------------------------------------------------
\bel{ND-RF}
\dist(x,E)\le 83\,R(x),~~x\in \RN \setminus E.
\ee
%----------------------------------------------------------
\par We let $]\theta[$ denote the (unique) positive integer such that $\theta\le\,]\theta[\,<\theta+1$. Let $m:=42]\theta[$. Consider a partition $\Kc_Q$ of the cube $Q$ into $m^n$ equal cubes $\{K_1,K_2,...,K_{m^n}\}$ of diameter $\diam Q/m$. Clearly, for every $K\in \Kc_Q$ we have
%----------------------------------------------------------
$$
\dist (Q,E)\le \dist (K,E).
$$
%----------------------------------------------------------
\par Let $K=Q(c_K,r_K)$. We have
%----------------------------------------------------------
$$
m\diam K=\diam Q \le \theta \,\dist (Q,E)\le \theta \,\dist(K,E)\le \theta \,\dist(c_K,E).
$$
%----------------------------------------------------------
By \rf{ND-RF},
%----------------------------------------------------------
$$
\dist(c_K,E)\le 83\,R(c_K)
$$
%----------------------------------------------------------
so that
%----------------------------------------------------------
$$
m\diam K = 2m\, r_K\le \theta \,\dist(c_K,E)
\le 83\,\theta \,R(c_K).
$$
%----------------------------------------------------------
Hence
%----------------------------------------------------------
$$
r_K \le (83\,\theta/2m)\,R(c_K).
$$
%----------------------------------------------------------
Since $m=42\,]\theta[\,\ge 42\,\theta$, we have $2m\ge 84\,\theta$ so that
%----------------------------------------------------------
$$
r_K < R(c_K).
$$
%----------------------------------------------------------
\par Now, by \rf{M-RX2},
%----------------------------------------------------------
$$
\mu(K)=\mu(Q(c_K,r_K)))<R(c_K)^{n-p}
$$
%----------------------------------------------------------
so that
%----------------------------------------------------------
$$
\mu(K)< r_K^{n-p}=(r_Q/m)^{n-p}.
$$
%----------------------------------------------------------
Hence,
%----------------------------------------------------------
$$
\mu(Q)= \sum_{i=1}^{m^n}\mu(K_i)\le m^n\,m^{p-n}\,r_Q^{n-p}=m^p\,r_Q^{n-p}\le
42^p(1+\theta)^p\, r_Q^{n-p}.
$$
%----------------------------------------------------------
The lemma is proved. \bx
%----------------------------------------------------------
\bigskip
%----------------------------------------------------------
%@@@@@@@@@@@@@@@@@@@@@@@@@@@@@@@@@@@@@@@@@@@@@@@@@@@@@@@@@@
%@@@@@@@@@@@@@@@@@@@@@@@@@@@@@@@@@@@@@@@@@@@@@@@@@@@@@@@@@@
%@@@@@@@@@@@@@@@@@@@@@@@@@@@@@@@@@@@@@@@@@@@@@@@@@@@@@@@@@@
%@@@@@@@@@@@@@@@@@@@@@@@@@@@@@@@@@@@@@@@@@@@@@@@@@@@@@@@@@@
%@@@@@@@@@@@@@@@@@@@@@@@@@@@@@@@@@@@@@@@@@@@@@@@@@@@@@@@@@@
%@@@@@@@@@@@@@@@@@@@@@@@@@@@@@@@@@@@@@@@@@@@@@@@@@@@@@@@@@@
%----------------------------------------------------------
\par The next lemma states that an inequality which is converse to the inequality \rf{M-RX3} is also true.
%----------------------------------------------------------
%@@@@@@@@@@@@@@@@@@@@@@@@@@@@@@@@@@@@@@@@@@@@@@@@@@@@@@@@@@
%@@@@@@@@@@@@@@@@@@@@@@@@@@@@@@@@@@@@@@@@@@@@@@@@@@@@@@@@@@
%@@@@@@@@@@@@@@@@@@@@@@@@@@@@@@@@@@@@@@@@@@@@@@@@@@@@@@@@@@
\begin{lemma} For every $x\in E$ the following inequality
%----------------------------------------------------------
$$
\mu(Q(x,5R(x)))\le 2^{14p}\,R(x)^{n-p}
$$
%----------------------------------------------------------
holds.
%----------------------------------------------------------
\end{lemma}
%----------------------------------------------------------
%@@@@@@@@@@@@@@@@@@@@@@@@@@@@@@@@@@@@@@@@@@@@@@@@@@@@@@@@@@
%@@@@@@@@@@@@@@@@@@@@@@@@@@@@@@@@@@@@@@@@@@@@@@@@@@@@@@@@@@
%@@@@@@@@@@@@@@@@@@@@@@@@@@@@@@@@@@@@@@@@@@@@@@@@@@@@@@@@@@
\par {\it Proof.} By subdividing each edge of the cube $K^{(x)}=Q(x,5R(x))$ into $20$ equal parts we can partition this cube into a family $A$ consisting of $20^n$ congruent cubes of diameter $\tfrac14\diam K^{(x)}$. Clearly, those cubes of the family $A$ which contain the point $x$, the center of the cube $K^{(x)}$, are a partition of the cube
$\tfrac12 K^{(x)}$ into a family of $2^n$ congruent cubes. %----------------------------------------------------------
\par Thus the set $5K^{(x)}\setminus\left(\tfrac12 K^{(x)}\right)$ is partitioned into a family $B\subset A$ consisting of $20^n-2^n$ congruent cubes of diameter $\tfrac14\diam K^{(x)}$. Clearly, for each cube $K\in B$
%----------------------------------------------------------
\bel{DX-1}
\dist(K,\{x\})\ge \diam K\left(=\tfrac14\diam K^{(x)}\right).
\ee
%----------------------------------------------------------
On the other hand, by part (i) of Lemma \reff{RPR}, the family
%----------------------------------------------------------
$$
6\Kc_E:=\{6 K^{(y)}=Q(x,6R(y)):y\in E\}
$$
%----------------------------------------------------------
consists of {\it non-overlapping} cubes. Hence
%----------------------------------------------------------
$$
\dist(5K^{(x)},E\setminus (5K^{(x)}))\ge R(x)=\tfrac12\diam K^{(x)}=2\diam K.
$$
%----------------------------------------------------------
Since $K\subset 5K^{(x)}$,
%----------------------------------------------------------
$$
\dist(K,E\setminus (5K^{(x)}))\ge
\dist(5K^{(x)},E\setminus (5K^{(x)}))\ge 2\diam K.
$$
%----------------------------------------------------------
Combining this inequality with \rf{DX-1} we obtain
%----------------------------------------------------------
$$
\dist(K,E)\ge \diam K.
$$
%----------------------------------------------------------
\par This property of the cube $K$ enables us to apply to $K$ the result of Lemma \reff{Q-NE} with $\theta=1$. By this lemma,
%----------------------------------------------------------
$$
\mu(K)\le 84^{p}\,r_K^{n-p}.
$$
%----------------------------------------------------------
Since $r_K=\tfrac14\,R(x)$, we have
%----------------------------------------------------------
$$
\mu(K)\le 84^{p}\,\left(\tfrac14\,R(x)\right)^{n-p}
=4^{p-n}84^{p}\,R(x)^{n-p}.
$$
%----------------------------------------------------------
\par By inequality \rf{M-RX2},
%----------------------------------------------------------
$$
\mu\left(\tfrac12\, K^{(x)}\right)=\mu\left(Q(x,\tfrac12\,R(x)\right)<R(x)^{n-p}.
$$
%----------------------------------------------------------
Finally we have
%----------------------------------------------------------
$$
\mu(5K^{(x)})\le\mu\left(\tfrac12\, K^{(x)}\right)+\sum_{K\in B}\mu(K)
$$
%----------------------------------------------------------
so that
%----------------------------------------------------------
$$
\mu(5K^{(x)})\le R(x)^{n-p}+(20^n-2^n)4^{p-n}84^{p}\,R(x)^{n-p}\le 2^{14p}\,R(x)^{n-p}
$$
%----------------------------------------------------------
proving the lemma.\bx
%----------------------------------------------------------
\bigskip
%----------------------------------------------------------
%@@@@@@@@@@@@@@@@@@@@@@@@@@@@@@@@@@@@@@@@@@@@@@@@@@@@@@@@@@
%@@@@@@@@@@@@@@@@@@@@@@@@@@@@@@@@@@@@@@@@@@@@@@@@@@@@@@@@@@
%@@@@@@@@@@@@@@@@@@@@@@@@@@@@@@@@@@@@@@@@@@@@@@@@@@@@@@@@@@
%----------------------------------------------------------
\par This lemma and inequality \rf{KXM} imply the following
%----------------------------------------------------------
%@@@@@@@@@@@@@@@@@@@@@@@@@@@@@@@@@@@@@@@@@@@@@@@@@@@@@@@@@@
%@@@@@@@@@@@@@@@@@@@@@@@@@@@@@@@@@@@@@@@@@@@@@@@@@@@@@@@@@@
%@@@@@@@@@@@@@@@@@@@@@@@@@@@@@@@@@@@@@@@@@@@@@@@@@@@@@@@@@@
%----------------------------------------------------------
\begin{corollary}\lbl{PR-5K} For every cube $K\in\Kc_E$ we have
%----------------------------------------------------------
%@@@@@@@@@@@@@@@@@@@@@@@@@@@@@@@@@@@@@@@@@@@@@@@@@@@@@@@@@@
%----------------------------------------------------------
\bel{M-DK}
%----------------------------------------------------------
2^{p-n}(\diam K)^{n-p}\le\mu(K)
\le 2^{15p}(\diam K)^{n-p}
%----------------------------------------------------------
\ee
%----------------------------------------------------------
and
%----------------------------------------------------------
$$
%----------------------------------------------------------
\mu(5K)
\le 2^{14p}\mu(K).
%----------------------------------------------------------
$$
%----------------------------------------------------------
\end{corollary}
%----------------------------------------------------------
%@@@@@@@@@@@@@@@@@@@@@@@@@@@@@@@@@@@@@@@@@@@@@@@@@@@@@@@@@@
%@@@@@@@@@@@@@@@@@@@@@@@@@@@@@@@@@@@@@@@@@@@@@@@@@@@@@@@@@@
%@@@@@@@@@@@@@@@@@@@@@@@@@@@@@@@@@@@@@@@@@@@@@@@@@@@@@@@@@@
%@@@@@@@@@@@@@@@@@@@@@@@@@      @@@@@@@@@@@@@@@@@@@@@@@@@@@
%@@@@@@@@@@@@@@@@@@@@@@@          @@@@@@@@@@@@@@@@@@@@@@@@@
%@@@@@@@@@@@@@@@@@@@@@              @@@@@@@@@@@@@@@@@@@@@@@
%@@@@@@@@@@@@@@@@@@@     SECTION 4    @@@@@@@@@@@@@@@@@@@@@
%@@@@@@@@@@@@@@@@@@@@@              @@@@@@@@@@@@@@@@@@@@@@@
%@@@@@@@@@@@@@@@@@@@@@@@          @@@@@@@@@@@@@@@@@@@@@@@@@
%@@@@@@@@@@@@@@@@@@@@@@@@@      @@@@@@@@@@@@@@@@@@@@@@@@@@@
%@@@@@@@@@@@@@@@@@@@@@@@@@@@@@@@@@@@@@@@@@@@@@@@@@@@@@@@@@@
%@@@@@@@@@@@@@@@@@@@@@@@@@@@@@@@@@@@@@@@@@@@@@@@@@@@@@@@@@@
%----------------------------------------------------------
%@@@@@@@@@@@@@@@@@@@@@@@@@@@@@@@@@@@@@@@@@@@@@@@@@@@@@@@@@@
%----------------------------------------------------------
\SECT{4. Sufficiency: the Sobolev norm of the function $f_1$.} {4}
%----------------------------------------------------------
\addtocontents{toc}{4. Sufficiency: the Sobolev norm of the function $f_1$. \hfill \thepage\\\par}
%----------------------------------------------------------
\indent
%----------------------------------------------------------
\par In this section, given a function $f$ satisfying the sufficiency condition of Theorem \reff{S-V2}, we define functions $f_1\in\LOP$ and $f_2\in\LPM$ such that $f_1+f_2=f$. We prove that
%----------------------------------------------------------
$$
\|f_1\|_{\LOP}\le C\lambda^{\frac1p}
$$
%----------------------------------------------------------
where $\lambda$ is the constant from inequality \rf{CR-V2} and $C=C(n,p)$. In the next section we show that $\|f_2\|_{\LPM}\le C\lambda^{\frac1p}$.
%----------------------------------------------------------
\par Let $E$ be the set constructed in the previous section. Since $E$ is a closed set, the set $\RN\setminus E$ is open so that it admits a Whitney decomposition $W_E$ into a family of non-overlapping cubes. In the next theorem we recall the main properties of this decomposition.
See, e.g. \cite{St}, or \cite{G}.
%----------------------------------------------------------
%@@@@@@@@@@@@@@@@@@@@@@@@@@@@@@@@@@@@@@@@@@@@@@@@@@@@@@@@@@
%@@@@@@@@@@@@@@@@@@@@@@@@@@@@@@@@@@@@@@@@@@@@@@@@@@@@@@@@@@
%@@@@@@@@@@@@@@@@@@@@@@@@@@@@@@@@@@@@@@@@@@@@@@@@@@@@@@@@@@
%@@@@@@@@@@@@@@@@@@@@@@@@@@@@@@@@@@@@@@@@@@@@@@@@@@@@@@@@@@
%@@@@@@@@@@@@@@@@@@@@@@@@@@@@@@@@@@@@@@@@@@@@@@@@@@@@@@@@@@
%@@@@@@@@@@@@@@@@@@@@@@@@@@@@@@@@@@@@@@@@@@@@@@@@@@@@@@@@@@
%@@@@@@@@@@@@@@@@@@@@@@@@@@@@@@@@@@@@@@@@@@@@@@@@@@@@@@@@@@
%----------------------------------------------------------
\begin{theorem}\lbl{Wcov} $W_E=\{Q_k\}$ is a countable family of non-overlapping cubes  such that
%----------------------------------------------------------
\par (i). $\RN\setminus E=\cup\{Q:Q\in W_E\}$;
%----------------------------------------------------------
\par (ii). For every cube $Q\in W_E$ we have
%----------------------------------------------------------
\bel{DQ-E}
\diam Q\le \dist(Q,E)\le 4\diam Q.
\ee
%----------------------------------------------------------
%@@@@@@@@@@@@@@@@@@@@@@@@@@@@@@@@@@@@@@@@@@@@@@@@@@@@@@@@@@
%----------------------------------------------------------
\end{theorem}
%----------------------------------------------------------
%@@@@@@@@@@@@@@@@@@@@@@@@@@@@@@@@@@@@@@@@@@@@@@@@@@@@@@@@@@
%@@@@@@@@@@@@@@@@@@@@@@@@@@@@@@@@@@@@@@@@@@@@@@@@@@@@@@@@@@
%@@@@@@@@@@@@@@@@@@@@@@@@@@@@@@@@@@@@@@@@@@@@@@@@@@@@@@@@@@
%@@@@@@@@@@@@@@@@@@@@@@@@@@@@@@@@@@@@@@@@@@@@@@@@@@@@@@@@@@
%@@@@@@@@@@@@@@@@@@@@@@@@@@@@@@@@@@@@@@@@@@@@@@@@@@@@@@@@@@
%@@@@@@@@@@@@@@@@@@@@@@@@@@@@@@@@@@@@@@@@@@@@@@@@@@@@@@@@@@
%----------------------------------------------------------
\par Let us note an important property of the Whitney cubes.
%----------------------------------------------------------
%@@@@@@@@@@@@@@@@@@@@@@@@@@@@@@@@@@@@@@@@@@@@@@@@@@@@@@@@@@
%@@@@@@@@@@@@@@@@@@@@@@@@@@@@@@@@@@@@@@@@@@@@@@@@@@@@@@@@@@
%----------------------------------------------------------
\begin{lemma}\lbl{WQ-M} For every cube $Q\in W_E$ the following inequality
%----------------------------------------------------------
$$
\mu(Q)\le 84^p\, r^{n-p}_Q
$$
%----------------------------------------------------------
holds.
\end{lemma}
%----------------------------------------------------------
%@@@@@@@@@@@@@@@@@@@@@@@@@@@@@@@@@@@@@@@@@@@@@@@@@@@@@@@@@@
%----------------------------------------------------------
\par {\it Proof.} Since $Q\in W_E$, by Theorem \reff{Wcov}, $\diam Q \le \dist (Q,E).$ It remains to apply Lemma \reff{Q-NE}  with $\theta=1$, and the lemma follows.\bx\bigskip
%----------------------------------------------------------
%@@@@@@@@@@@@@@@@@@@@@@@@@@@@@@@@@@@@@@@@@@@@@@@@@@@@@@@@@@
%@@@@@@@@@@@@@@@@@@@@@@@@@@@@@@@@@@@@@@@@@@@@@@@@@@@@@@@@@@
%@@@@@@@@@@@@@@@@@@@@@@@@@@@@@@@@@@@@@@@@@@@@@@@@@@@@@@@@@@
%@@@@@@@@@@@@@@@@@@@@@@@@@@@@@@@@@@@@@@@@@@@@@@@@@@@@@@@@@@
%----------------------------------------------------------
\par Combining this result with the second inequality in \rf{M-DK} we obtain the following
%----------------------------------------------------------
%@@@@@@@@@@@@@@@@@@@@@@@@@@@@@@@@@@@@@@@@@@@@@@@@@@@@@@@@@@
%@@@@@@@@@@@@@@@@@@@@@@@@@@@@@@@@@@@@@@@@@@@@@@@@@@@@@@@@@@
%@@@@@@@@@@@@@@@@@@@@@@@@@@@@@@@@@@@@@@@@@@@@@@@@@@@@@@@@@@
%@@@@@@@@@@@@@@@@@@@@@@@@@@@@@@@@@@@@@@@@@@@@@@@@@@@@@@@@@@
%----------------------------------------------------------
\begin{corollary}\lbl{P-WKE} Every cube $Q\in W_E\cup\Kc_E$ satisfies the inequality
%----------------------------------------------------------
\bel{WKE-S}
\mu(Q)\le 2^{15p}(\diam Q)^{n-p}.
\ee
%----------------------------------------------------------
\par Thus for every  $Q',Q''\in W_E\cup\Kc_E$ the conditions \rf{G-MD} of Theorem \reff{S-V2} hold.
%----------------------------------------------------------
\end{corollary}
%----------------------------------------------------------
%@@@@@@@@@@@@@@@@@@@@@@@@@@@@@@@@@@@@@@@@@@@@@@@@@@@@@@@@@@
%@@@@@@@@@@@@@@@@@@@@@@@@@@@@@@@@@@@@@@@@@@@@@@@@@@@@@@@@@@
%@@@@@@@@@@@@@@@@@@@@@@@@@@@@@@@@@@@@@@@@@@@@@@@@@@@@@@@@@@
%----------------------------------------------------------
\bigskip
%----------------------------------------------------------
\par We are also needed certain additional properties of
Whitney cubes which we present in the next lemma. These
properties easily follow from constructions of Whitney decomposition presented in \cite{St} and \cite{G}.
%----------------------------------------------------------
%@@@@@@@@@@@@@@@@@@@@@@@@@@@@@@@@@@@@@@@@@@@@@@@@@@@@@@@@@@
\par Given a cube $Q\subset\RN$ let $Q^*:=\frac{9}{8}Q$.
%----------------------------------------------------------
%@@@@@@@@@@@@@@@@@@@@@@@@@@@@@@@@@@@@@@@@@@@@@@@@@@@@@@@@@@
%@@@@@@@@@@@@@@@@@@@@@@@@@@@@@@@@@@@@@@@@@@@@@@@@@@@@@@@@@@
%@@@@@@@@@@@@@@@@@@@@@@@@@@@@@@@@@@@@@@@@@@@@@@@@@@@@@@@@@@
%@@@@@@@@@@@@@@@@@@@@@@@@@@@@@@@@@@@@@@@@@@@@@@@@@@@@@@@@@@
%----------------------------------------------------------
\begin{lemma}\lbl{Wadd}
%----------------------------------------------------------
%@@@@@@@@@@@@@@@@@@@@@@@@@@@@@@@@@@@@@@@@@@@@@@@@@@@@@@@@@@
%----------------------------------------------------------
(1). If $Q,K\in W_E$ and $Q^*\cap K^*\ne\emptyset$, then
%----------------------------------------------------------
$$
\frac{1}{4}\diam Q\le \diam K\le 4\diam Q.
$$
%----------------------------------------------------------
\smallskip
%----------------------------------------------------------
%@@@@@@@@@@@@@@@@@@@@@@@@@@@@@@@@@@@@@@@@@@@@@@@@@@@@@@@@@@
\par (2). For every cube $K\in W_E$ there are at most
$N=N(n)$ cubes from the family
%----------------------------------------------------------
$W_E^*:=\{Q^*:Q\in W_E\}$
%----------------------------------------------------------
which intersect $K^*$.
%----------------------------------------------------------
\medskip
%----------------------------------------------------------
\par (3). If $Q,K\in W_E$, then $Q^*\cap K^*\ne\emptyset$
if and only if  $Q\cap K\ne\emptyset$.
%----------------------------------------------------------
%@@@@@@@@@@@@@@@@@@@@@@@@@@@@@@@@@@@@@@@@@@@@@@@@@@@@@@@@@@
\end{lemma}
%----------------------------------------------------------
\medskip
%----------------------------------------------------------
%@@@@@@@@@@@@@@@@@@@@@@@@@@@@@@@@@@@@@@@@@@@@@@@@@@@@@@@@@@
%@@@@@@@@@@@@@@@@@@@@@@@@@@@@@@@@@@@@@@@@@@@@@@@@@@@@@@@@@@
%@@@@@@@@@@@@@@@@@@@@@@@@@@@@@@@@@@@@@@@@@@@@@@@@@@@@@@@@@@
%@@@@@@@@@@@@@@@@@@@@@@@@@@@@@@@@@@@@@@@@@@@@@@@@@@@@@@@@@@
%@@@@@@@@@@@@@@@@@@@@@@@@@@@@@@@@@@@@@@@@@@@@@@@@@@@@@@@@@@
%@@@@@@@@@@@@@@@@@@@@@@@@@@@@@@@@@@@@@@@@@@@@@@@@@@@@@@@@@@
%----------------------------------------------------------
\par Note that inequality \rf{DQ-E} implies the following property of Whitney cubes:
%----------------------------------------------------------
\bel{9QE}
(9Q)\cap E\ne\emp~~~ \text{for every}~~~ Q\in W_E.
\ee
%----------------------------------------------------------
\par Let us fix a constant $\tau\ge 9$. Then by the above property
%----------------------------------------------------------
$$
(\tau Q)\cap E\ne\emp.
$$
%----------------------------------------------------------
\par To every cube $Q\in W_E$ we assign a point $a_Q\in E$ such that
%----------------------------------------------------------
\bel{AQ-T}
a_Q\in \tau Q.
\ee
%----------------------------------------------------------
For instance, one can choose $a_Q$ to be a point nearest to $Q$ on the set $E$. Then, by the property \rf{9QE},
$a_Q\in \tau Q$ with $\tau=9$.
%@@@@@@@@@@@@@@@@@@@@@@@@@@@@@@@@@@@@@@@@@@@@@@@@@@@@@@@@@@
%@@@@@@@@@@@@@@@@@@@@@@@@@@@@@@@@@@@@@@@@@@@@@@@@@@@@@@@@@@
%@@@@@@@@@@@@@@@@@@@@@@@@@@@@@@@@@@@@@@@@@@@@@@@@@@@@@@@@@@
%----------------------------------------------------------
\par Let $\Phi_E:=\{\varphi_Q:Q\in W_E\}$ be a smooth partition of unity subordinated to the Whitney decomposition $W_E$. Recall the main properties of this partition.
%@@@@@@@@@@@@@@@@@@@@@@@@@@@@@@@@@@@@@@@@@@@@@@@@@@@@@@@@@@
\begin{lemma}\lbl{P-U} The family of functions $\Phi_E$ has the following properties:
%----------------------------------------------------------
%@@@@@@@@@@@@@@@@@@@@@@@@@@@@@@@@@@@@@@@@@@@@@@@@@@@@@@@@@@
%----------------------------------------------------------
\medskip
%----------------------------------------------------------
\par (a). $\varphi_Q\in C^\infty(\RN)$ and
$0\le\varphi_Q\le 1$ for every $Q\in W_E$;\smallskip
%----------------------------------------------------------
\medskip
%----------------------------------------------------------
\par (b). $\supp \varphi_Q\subset Q^*(:=\frac{9}{8}Q),$
$Q\in W_E$;\smallskip
%----------------------------------------------------------
\medskip
%----------------------------------------------------------
\par (c). $\sum\{\varphi_Q(x):Q\in W_E\}=1$ for every
$x\in\RN\setminus S$;\smallskip
%----------------------------------------------------------
\medskip
%----------------------------------------------------------
\par (d). $ \|\nabla\varphi_Q(x)\| \le C(n)/\diam Q\,\,$
for every $Q\in W_E$ and every $x\in\RN$.
%----------------------------------------------------------
%@@@@@@@@@@@@@@@@@@@@@@@@@@@@@@@@@@@@@@@@@@@@@@@@@@@@@@@@@@
\end{lemma}
%----------------------------------------------------------
%@@@@@@@@@@@@@@@@@@@@@@@@@@@@@@@@@@@@@@@@@@@@@@@@@@@@@@@@@@
%@@@@@@@@@@@@@@@@@@@@@@@@@@@@@@@@@@@@@@@@@@@@@@@@@@@@@@@@@@
%----------------------------------------------------------
\medskip
%@@@@@@@@@@@@@@@@@@@@@@@@@@@@@@@@@@@@@@@@@@@@@@@@@@@@@@@@@@
\par We turn to definition of the functions $f_1\in\LOP$ and $f_2\in\LPM$ which provides an almost optimal decomposition of a function $f$ satisfying the sufficiency condition of Theorem \reff{MAIN-CR}.
%----------------------------------------------------------
%@@@@@@@@@@@@@@@@@@@@@@@@@@@@@@@@@@@@@@@@@@@@@@@@@@@@@@@@@@
%@@@@@@@@@@@@@@@@@@@@@@@@@@@@@@@@@@@@@@@@@@@@@@@@@@@@@@@@@@
%----------------------------------------------------------
\par Let $\tf_1:E\to\R$ be a function defined by the following formula:
%----------------------------------------------------------
\bel{DEF-TF1}
\tf_1(x):=\av{Q}=\frac{1}{\mu(Q)}\intl_Q f\,d\mu~~~\text{for every}~~~x\in E.
\ee
%----------------------------------------------------------
Here $Q=K^{(x)}=Q(x,R(x))$, see \rf{D-KX}.
%----------------------------------------------------------
%@@@@@@@@@@@@@@@@@@@@@@@@@@@@@@@@@@@@@@@@@@@@@@@@@@@@@@@@@@
\par Using the Whitney extension formula we extend $\tf_1$ from $E$ to all of $\RN$. We denote this extension by $f_1$. Thus:
%----------------------------------------------------------
\bel{DEF-F1}
f_1(x):=\left \{
%----------------------------------------------------------
\begin{array}{ll}
\tf_1(x),& x\in E,\\\\
\smed\limits_{Q\in W_E}
\varphi_Q(x)\tf_1(a_Q),& x\in\RN\setminus E.
\end{array}
%----------------------------------------------------------
\right.
\ee
%----------------------------------------------------------
Finally we put
%----------------------------------------------------------
$$
f_2:=f-f_1.
$$
%----------------------------------------------------------
\par Our goal is to prove that under Theorem's  \reff{S-V2} conditions the following inequality
%----------------------------------------------------------
\bel{F12}
\|f_1\|_{\LOP}+\|f_2\|_{\LPM}\le C\lambda^{\frac1p}
\ee
%----------------------------------------------------------
%@@@@@@@@@@@@@@@@@@@@@@@@@@@@@@@@@@@@@@@@@@@@@@@@@@@@@@@@@@
holds. Here $C$ is a constant depending only on $n,p,$ and $\tau$.
%@@@@@@@@@@@@@@@@@@@@@@@@@@@@@@@@@@@@@@@@@@@@@@@@@@@@@@@@@@
\par Let us estimate the norm $\|f_1\|_{\LOP}$. Let
$K$ be a cube in $\RN$  and let %----------------------------------------------------------
$$
V_K:=\{Q\in W_E:Q\cap K\ne \emptyset\}.
$$
%----------------------------------------------------------
%@@@@@@@@@@@@@@@@@@@@@@@@@@@@@@@@@@@@@@@@@@@@@@@@@@@@@@@@@@
\begin{lemma}\lbl{GR-C}. For every cube $K\in W_E$ the following inequality
%----------------------------------------------------------
$$
\intl_K\|\nabla f_1(x)\|^pdx\le C(n)\,\smed_{Q\in V_K}\,\frac{|\tilde{f}_1(a_K)-\tilde{f}_1(a_Q)|^p}{(\diam K)^{p-n}}
$$
%----------------------------------------------------------
holds.
%@@@@@@@@@@@@@@@@@@@@@@@@@@@@@@@@@@@@@@@@@@@@@@@@@@@@@@@@@@
\end{lemma}
%@@@@@@@@@@@@@@@@@@@@@@@@@@@@@@@@@@@@@@@@@@@@@@@@@@@@@@@@@@
%----------------------------------------------------------
\par{\it Proof.} Let $x\in K$. Since $x\in \RN\setminus E$, by the extension formula \rf{DEF-F1} and by properties (b) and (c) of Lemma \reff{P-U}, we have
%----------------------------------------------------------
\be
\|\nabla f_1(x)\|&=&\|\nabla(f_1(x)-\tilde{f}_1(a_K))\|
=\left\|\,\nabla \left(\,\sum_{Q\in W_E}\varphi_Q(x)(\tilde{f}_1(a_Q)-\tilde{f}_1(a_K))
\right)\right\|\nn\\
&=&\left\|\sum_{Q\in W_E}(\tilde{f}_1(a_Q)-\tilde{f}_1(a_K))\nabla
\varphi_Q(x)\right\|\nn\\
&=&
\left\|\sum\left\{(\tilde{f}_1(a_Q)-\tilde{f}_1(a_K))\nabla \varphi_Q(x):
Q^*\cap K\ne \emptyset, Q\in W_E\right\}\right\|\nn\\
&\le&
\sum\left\{|\tilde{f}_1(a_Q)-\tilde{f}_1(a_K)|~\|\nabla \varphi_Q(x)\| :Q^*\cap K\ne \emptyset,Q\in W_E\right\}.\nn
\ee
%----------------------------------------------------------
Hence, by property (d) of Lemma \reff{P-U},
%----------------------------------------------------------
\be
\|\nabla f_1(x)\|&\le& C\sbig\left\{\frac{|\tilde{f}_1(a_Q)-\tilde{f}_1(a_K)|}
{\diam Q}:Q^*\cap K\ne \emptyset,Q\in W_E\right\}. \nn
\ee
%----------------------------------------------------------
By Lemma \reff{Wadd}, $Q^*\cap K\ne \emptyset$ iff $Q\cap K\ne \emptyset$. Also, by this lemma, $\diam Q\sim \diam K$. Hence
%----------------------------------------------------------
$$
\|\nabla f_1(x)\|\le C\sbig\left\{\frac{|\tilde{f}_1(a_Q)-\tilde{f}_1(a_K)|}
{\diam K}:Q\cap K\ne \emptyset,Q\in W_E\right\}.
$$
%----------------------------------------------------------
Integrating this inequality over the cube $K$, we obtain
%----------------------------------------------------------
$$
\intl_K\|\nabla f_1(x)\|^pdx\le C\sbig\left\{\frac{|\tilde{f}_1(a_Q)-\tilde{f}_1(a_K)|^p}
{(\diam K)^{p-n}}:Q\cap K\ne \emptyset,Q\in W_E\right\}
$$
%----------------------------------------------------------
proving the lemma.\bx\bigskip
%----------------------------------------------------------
%@@@@@@@@@@@@@@@@@@@@@@@@@@@@@@@@@@@@@@@@@@@@@@@@@@@@@@@@@@
%@@@@@@@@@@@@@@@@@@@@@@@@@@@@@@@@@@@@@@@@@@@@@@@@@@@@@@@@@@
%@@@@@@@@@@@@@@@@@@@@@@@@@@@@@@@@@@@@@@@@@@@@@@@@@@@@@@@@@@
%@@@@@@@@@@@@@@@@@@@@@@@@@@@@@@@@@@@@@@@@@@@@@@@@@@@@@@@@@@
%----------------------------------------------------------
\par Recall that the set $E$ consists of {\it isolated points} of $\RN$ so that the function  $f_1\in C^{\infty}(\RN)$. This observation and Lemma \reff{GR-C} enable us to estimate its Sobolev seminorm as follows:
%----------------------------------------------------------
\bel{FE-K}
\|\nabla f_1\|^p_{\LPRN}\le C(n)\smed_{K\in\,W_E}
\,\,\smed_{Q\in V_K} \frac{|\tilde{f}_1(a_Q)-\tilde{f}_1(a_K)|^p}{(\diam K)^{p-n}}.
\ee
%----------------------------------------------------------
\par Let us slightly simplify this inequality.
By $\tK$ we denote a cube which maximize  the quantity $|\tilde{f}_1(a_Q)-\tilde{f}_1(a_K)|$ on the family $V_K$; thus
%----------------------------------------------------------
\bel{D-KW}
\max_{Q\in V_K}|\tilde{f}_1(a_Q)-\tilde{f}_1(a_K)|=
|\tilde{f}_1(a_{\tK})-\tilde{f}_1(a_K)|.
\ee
%----------------------------------------------------------
(Of course, $\tK$ depends on $K$ and $f$.)
%----------------------------------------------------------
By part (2) of Lemma \reff{Wadd}, $\# V_K\le N(n)$ so that, by \rf{FE-K},
%----------------------------------------------------------
\bel{F-IL}
\|\nabla f_1\|^p_{\LPRN}\le C(n)\smed_{K\in\,W_E}
\frac{|\tilde{f}_1(a_{\tK})-\tilde{f}_1(a_K)|^p}
{(\diam K)^{p-n}}.
\ee
%----------------------------------------------------------
\par Let us show that we can omit in the right hand side of this inequality those cubes $K\in W_E$ which lie inside of cubes from the family $\Kc_E$.
%----------------------------------------------------------
\par Let
%----------------------------------------------------------
\bel{ETA}
\eta:=\tfrac{1}{21}\, \tau.
\ee
%----------------------------------------------------------
%@@@@@@@@@@@@@@@@@@@@@@@@@@@@@@@@@@@@@@@@@@@@@@@@@@@@@@@@@@
%@@@@@@@@@@@@@@@@@@@@@@@@@@@@@@@@@@@@@@@@@@@@@@@@@@@@@@@@@@
%@@@@@@@@@@@@@@@@@@@@@@@@@@@@@@@@@@@@@@@@@@@@@@@@@@@@@@@@@@
%@@@@@@@@@@@@@@@@@@@@@@@@@@@@@@@@@@@@@@@@@@@@@@@@@@@@@@@@@@
%----------------------------------------------------------
\begin{lemma}\lbl{SMC} Let $x\in E$ and let $K\in W_E$ be a cube such that
%----------------------------------------------------------
\bel{KKV}
K\cap (\eta K^{(x)})\ne\emp.
\ee
%----------------------------------------------------------
Then
%----------------------------------------------------------
$$
\tau Q\subset K^{(x)}~~~\text{for every}~~~Q\in V_K.
$$
%----------------------------------------------------------
Furthermore, $a_{Q}=x$ for every $Q\in V_K$.
%----------------------------------------------------------
%@@@@@@@@@@@@@@@@@@@@@@@@@@@@@@@@@@@@@@@@@@@@@@@@@@@@@@@@@@
\end{lemma}
%@@@@@@@@@@@@@@@@@@@@@@@@@@@@@@@@@@@@@@@@@@@@@@@@@@@@@@@@@@
%----------------------------------------------------------
\par{\it Proof.} Let $K=Q(x_K,r_K)$ and let $Q=Q(x_{Q},r_{Q})$. Since $K\in W_E$, we have
%----------------------------------------------------------
$$
\diam K\le 4\dist(K,E).
$$
%----------------------------------------------------------
Hence $\diam K\le\,4\dist(K,\{x\})$ so that, by \rf{KKV},
%----------------------------------------------------------
\bel{DGV}
\diam K\le 4(\tfrac12\diam(\eta K^{(x)}))=
2\eta\diam K^{(x)}.
\ee
%----------------------------------------------------------
Since $Q\in V_K$, we have $Q\cap K\ne\emp$ and $Q\in W_E$, so that, by Lemma \reff{Wadd},
$$\diam Q\le 4\diam K.$$ Hence
%----------------------------------------------------------
\bel{D-KV}
\diam Q\le 8\eta\diam K^{(x)}.
\ee
%----------------------------------------------------------
Furthermore, by \rf{KKV} and \rf{DGV},
%----------------------------------------------------------
\be
\|x-x_K\|&\le& \tfrac12\diam(\eta K^{(x)}) +\tfrac12\diam K\nn\\&\le& (\eta/2)\diam K^{(x)}+\tfrac12 (2\eta\diam K^{(x)})=
\tfrac32\eta\diam K^{(x)}\nn
\ee
%----------------------------------------------------------
so that
%----------------------------------------------------------
$$
K\subset (3\eta+2\eta) K^{(x)}=5\eta K^{(x)}.
$$
%----------------------------------------------------------
Hence
%----------------------------------------------------------
$$
\tau K\subset 5\eta \tau K^{(x)}\subset K^{(x)}.
$$
%----------------------------------------------------------
(Recall that $\eta=1/(21\tau)$, see \rf{ETA}.) Since the cubes of the family $\Kc_E=\{K^{(x)}:x\in E\}$ are pairwise disjoint,
%----------------------------------------------------------
$$
(\tau K)\cap E \subset K^{(x)}\cap E=\{x\}.
$$
%----------------------------------------------------------
Hence $a_K=x$.
%----------------------------------------------------------
\par In the same fashion we show that $a_{Q}=x$. In fact, by \rf{D-KV},
%----------------------------------------------------------
\be
\|x-x_{Q}\|&\le& \|x-x_{K}\|+\|x_K-x_{Q}\|
\le \tfrac32\eta\diam K^{(x)}+\tfrac12\diam K+\tfrac12\diam Q\nn\\&\le& \tfrac32\eta\diam K^{(x)}+\tfrac12(2\eta\diam K^{(x)})+\tfrac12(8\eta\diam K^{(x)})=
\tfrac{13}{2}\,\eta\diam K^{(x)}.\nn
\ee
%----------------------------------------------------------
This inequality and \rf{D-KV} imply the following:
%----------------------------------------------------------
$$
Q\subset (13\eta+8\eta) K^{(x)}=21\eta K^{(x)}.
$$
%----------------------------------------------------------
Hence
%----------------------------------------------------------
$$
\tau Q\subset 21\tau\eta K^{(x)}\subset K^{(x)}
$$
%----------------------------------------------------------
so that
%----------------------------------------------------------
$$
(\tau Q)\cap E \subset K^{(x)}\cap E=\{x\}.
$$
%----------------------------------------------------------
Thus $a_{Q}=a_K=x$, and the proof is finished.\bx\bigskip
%----------------------------------------------------------
%@@@@@@@@@@@@@@@@@@@@@@@@@@@@@@@@@@@@@@@@@@@@@@@@@@@@@@@@@@
%@@@@@@@@@@@@@@@@@@@@@@@@@@@@@@@@@@@@@@@@@@@@@@@@@@@@@@@@@@
%@@@@@@@@@@@@@@@@@@@@@@@@@@@@@@@@@@@@@@@@@@@@@@@@@@@@@@@@@@
%@@@@@@@@@@@@@@@@@@@@@@@@@@@@@@@@@@@@@@@@@@@@@@@@@@@@@@@@@@
%----------------------------------------------------------
\par The lemma motivates us to introduce two subsets of $\RN$ defined by the following formulas:
%----------------------------------------------------------
\bel{TE-D}
T_E:=\bigcup_{x\in E} ~K^{(x)}=\bigcup\{K:K\in\Kc_E\}
\ee
%----------------------------------------------------------
and
%----------------------------------------------------------
\bel{TET}
T_{E,\tau}:=\bigcup_{x\in E} ~\eta K^{(x)}=\bigcup\{\eta K:K\in\Kc_E\}.
\ee
%----------------------------------------------------------
\par Let us also introduce a collection of cubes
%----------------------------------------------------------
\bel{A-D}
\Ac:=\{K\in W_E:K\cap T_{E,\tau}=\emp\}.
\ee
%----------------------------------------------------------
\par By the lemma
%----------------------------------------------------------
$$
\frac{|\tilde{f}_1(a_{\tK})-\tilde{f}_1(a_K)|^p}
{(\diam K)^{p-n}}=0
$$
%----------------------------------------------------------
provided $K\cap T_{E,\tau}\ne\emp$ or, equivalently, $K\in W_E\setminus\Ac$. Combining this with inequality \rf{F-IL} we obtain the following
%@@@@@@@@@@@@@@@@@@@@@@@@@@@@@@@@@@@@@@@@@@@@@@@@@@@@@@@@@@
%----------------------------------------------------------
%@@@@@@@@@@@@@@@@@@@@@@@@@@@@@@@@@@@@@@@@@@@@@@@@@@@@@@@@@@
\begin{corollary}\lbl{NF1-1} We have
%----------------------------------------------------------
%@@@@@@@@@@@@@@@@@@@@@@@@@@@@@@@@@@@@@@@@@@@@@@@@@@@@@@@@@@
%----------------------------------------------------------
$$
\|\nabla f_1\|^p_{\LPRN}\le C(n)\smed_{K\in\Ac}\,
\frac{|\tilde{f}_1(a_{\tK})-\tilde{f}_1(a_K)|^p}
{(\diam K)^{p-n}}.
$$
%----------------------------------------------------------
\end{corollary}
%@@@@@@@@@@@@@@@@@@@@@@@@@@@@@@@@@@@@@@@@@@@@@@@@@@@@@@@@@@
%----------------------------------------------------------
\par We introduce two subfamilies of the family $\Ac$:
%----------------------------------------------------------
\bel{A1-D}
\Ac_1:=\{K\in W_E:K\cap T_{E}=\emp\}.
\ee
%----------------------------------------------------------
and
%----------------------------------------------------------
$$
\Ac_2:=\{K\in W_E:K\cap T_{E}\ne\emp, K\cap T_{E,\tau}=\emp\}.
$$
%----------------------------------------------------------
Clearly, $\Ac_1,\Ac_2$ is a partition of $\Ac$, i.e., $\Ac_1\cup \Ac_2=A$, $\Ac_1\cap \Ac_2=\emp$. Hence, by Corollary \reff{NF1-1},
%----------------------------------------------------------
\be
\|\nabla f_1\|^p_{\LPRN}&\le& C(n)\smed_{K\in \Ac}\,
\frac{|\tilde{f}_1(a_{\tK})-\tilde{f}_1(a_K)|^p}
{(\diam K)^{p-n}}\nn\\
&=&
C(n)\smed_{K\in \Ac_1}\,
\frac{|\tilde{f}_1(a_{\tK})-\tilde{f}_1(a_K)|^p}
{(\diam K)^{p-n}}+C(n)\smed_{K\in \Ac_2}\,
\frac{|\tilde{f}_1(a_{\tK})-\tilde{f}_1(a_K)|^p}
{(\diam K)^{p-n}}.\nn
\ee
%----------------------------------------------------------
%@@@@@@@@@@@@@@@@@@@@@@@@@@@@@@@@@@@@@@@@@@@@@@@@@@@@@@@@@@
%----------------------------------------------------------
\par Recall that the function $f:\RN\to \R$ satisfies the sufficiency condition of Theorem \reff{S-V2}. Thus
there exists a constant $\lambda>0$ such that for every finite family $\Qc$ of pairwise disjoint cubes and arbitrary mappings $\Qc\ni Q\mapsto Q'\in\Qc$ and
$\Qc\ni Q\mapsto Q''\in\Qc$ such that $Q'\cup Q''\subset\gamma Q$, inequality \rf{CR-V2} holds.
%----------------------------------------------------------
%@@@@@@@@@@@@@@@@@@@@@@@@@@@@@@@@@@@@@@@@@@@@@@@@@@@@@@@@@@
%@@@@@@@@@@@@@@@@@@@@@@@@@@@@@@@@@@@@@@@@@@@@@@@@@@@@@@@@@@
\begin{lemma}\lbl{TQ} Let $Q\in \Ac$. Then for every  cube $Q'\in W_E$ such that $Q'\cap Q\ne\emp$ we have
%----------------------------------------------------------
$$
K^{(a_{Q'})}\subset (22\tau^2)\, Q.
$$
%----------------------------------------------------------
%@@@@@@@@@@@@@@@@@@@@@@@@@@@@@@@@@@@@@@@@@@@@@@@@@@@@@@@@@@
%----------------------------------------------------------
%@@@@@@@@@@@@@@@@@@@@@@@@@@@@@@@@@@@@@@@@@@@@@@@@@@@@@@@@@@
\end{lemma}
%@@@@@@@@@@@@@@@@@@@@@@@@@@@@@@@@@@@@@@@@@@@@@@@@@@@@@@@@@@
\par {\it Proof.} Recall that $\eta=1/(21\tau)$, see \rf{ETA}, and
%----------------------------------------------------------
$$
T_{E,\tau}=\{Q\in W_E:
~Q\cap(\eta K^{(x)})=\emp~~\text{for every}~~x\in E\},
$$
%----------------------------------------------------------
see \rf{TET}. Since $Q\in \Ac$, we have $Q\cap T_{E,\tau}=\emp$ so that %----------------------------------------------------------
$$
\dist(a_Q,Q)>\tfrac12\eta\diam K^{(a_Q)}.
$$
%----------------------------------------------------------
\par Let $Q=Q(x_Q,r_Q), Q'=Q(x_{Q'},r_{Q'})$. Since $a_Q\in \tau Q$, we have $\|a_Q-x_Q\|\le\tau r_Q$ so that
%----------------------------------------------------------
$$
\dist(K^{(a_Q)},Q)\le \|a_Q-x_Q\|\le\tau r_Q
$$
%----------------------------------------------------------
proving that
%----------------------------------------------------------
$$
\diam K^{(a_Q)}\le 2\tau\,r_Q/\eta=(21\tau^2)\diam Q.
$$
%----------------------------------------------------------
Hence
%----------------------------------------------------------
$$
K^{(a_Q)}\subset (\tau+21\tau^2)Q\subset 22\tau^2\,Q.
$$
%----------------------------------------------------------
\par Now if $a_{Q'}=a_Q$, then $K^{(a_{Q'})}=K^{(a_Q)}$ so that in this case $K^{(a_{Q'})}=K^{(a_Q)}\subset 22\tau^2\,Q.$\medskip
%----------------------------------------------------------
\par Suppose that $a_{Q'}\ne a_Q$. Then, by part (ii) of Lemma \reff{PR-KE},
%----------------------------------------------------------
\bel{DKA}
\diam K^{(a_{Q'})}\le \tau(\diam Q+\diam Q')\le \tau(\diam Q+4\diam Q)=10\tau r_Q.
\ee
%----------------------------------------------------------
Since $a_{Q'}\in\tau Q'$, we have $\|a_{Q'}-x_{Q'}\|\le\tau r_{Q'}$. Since $Q\cap Q'\ne\emp$ and $r_{Q'}\le 4r_Q$,
%----------------------------------------------------------
$$
\|x_{Q'}-x_{Q}\|\le r_{Q}+r_{Q'}\le 5r_Q
$$
%----------------------------------------------------------
so that
%----------------------------------------------------------
$$
\|a_{Q'}-x_{Q}\|\le 5r_{Q}+\tau r_{Q'}\le (4\tau+5)r_Q.
$$
%----------------------------------------------------------
Combining this inequality with \rf{DKA} we obtain
%----------------------------------------------------------
$$
K^{(a_{Q'})}\subset (4\tau+5+10\tau)Q\subset 22\tau^2 Q.
$$
%----------------------------------------------------------
\par The lemma is proved.\bx\bigskip
%----------------------------------------------------------
%@@@@@@@@@@@@@@@@@@@@@@@@@@@@@@@@@@@@@@@@@@@@@@@@@@@@@@@@@@
%@@@@@@@@@@@@@@@@@@@@@@@@@@@@@@@@@@@@@@@@@@@@@@@@@@@@@@@@@@
%@@@@@@@@@@@@@@@@@@@@@@@@@@@@@@@@@@@@@@@@@@@@@@@@@@@@@@@@@@
%@@@@@@@@@@@@@@@@@@@@@@@@@@@@@@@@@@@@@@@@@@@@@@@@@@@@@@@@@@
%----------------------------------------------------------
\par We are needed the following combinatorial
%----------------------------------------------------------
%@@@@@@@@@@@@@@@@@@@@@@@@@@@@@@@@@@@@@@@@@@@@@@@@@@@@@@@@@@
\begin{proposition}\lbl{C-PT} Let $N\in\N$ and let $\Bc=\{Q\}$ be a collection of cubes in $\RN$. Suppose that for every cube $Q\in\Bc$ there exist at most $N$ cubes from $\Bc$ which have common points with $Q$.
%----------------------------------------------------------
Then the family $\Bc$ can be partitioned into at most $N+1$ families of pairwise disjoint cubes.
\end{proposition}
%@@@@@@@@@@@@@@@@@@@@@@@@@@@@@@@@@@@@@@@@@@@@@@@@@@@@@@@@@@
\par {\it Proof.} The proposition immediately follows from the next well-known result in the graph theory (see, e.g., \cite{JT}): {\it Every graph can be colored with one more color than the maximum vertex degree.}\bx\medskip
%----------------------------------------------------------
%@@@@@@@@@@@@@@@@@@@@@@@@@@@@@@@@@@@@@@@@@@@@@@@@@@@@@@@@@@
%@@@@@@@@@@@@@@@@@@@@@@@@@@@@@@@@@@@@@@@@@@@@@@@@@@@@@@@@@@
%@@@@@@@@@@@@@@@@@@@@@@@@@@@@@@@@@@@@@@@@@@@@@@@@@@@@@@@@@@
%@@@@@@@@@@@@@@@@@@@@@@@@@@@@@@@@@@@@@@@@@@@@@@@@@@@@@@@@@@
%----------------------------------------------------------
\begin{proposition}\lbl{EA1} Suppose that the hypothesis of Theorem \reff{S-V2} holds with $\gamma=22\tau^2$. Then   %@@@@@@@@@@@@@@@@@@@@@@@@@@@@@@@@@@@@@@@@@@@@@@@@@@@@@@@@@@
%----------------------------------------------------------
\bel{RE-A1}
\smed_{K\in \Ac_1}\,
\frac{|\tilde{f}_1(a_{\tK})-\tilde{f}_1(a_K)|^p}
{(\diam K)^{p-n}}\le C(n)\lambda.
\ee
%----------------------------------------------------------
%@@@@@@@@@@@@@@@@@@@@@@@@@@@@@@@@@@@@@@@@@@@@@@@@@@@@@@@@@@
\end{proposition}
%@@@@@@@@@@@@@@@@@@@@@@@@@@@@@@@@@@@@@@@@@@@@@@@@@@@@@@@@@@
%----------------------------------------------------------
\par {\it Proof.} Let $\tAc\subset \Ac_1$ be a {\it finite family of pairwise disjoint cubes}. Prove that
%----------------------------------------------------------
\bel{S-TA}
\smed_{K\in \tAc}\,
\frac{|\tilde{f}_1(a_{\tK})-\tilde{f}_1(a_K)|^p}
{(\diam K)^{p-n}}\le \lambda
\ee
%----------------------------------------------------------
Let
%----------------------------------------------------------
\bel{D-A1-N}
\Qc:=\tAc\cup \Kc_E.
\ee
%----------------------------------------------------------
Since the cubes of the family $\Ac_1$ and  the cubes of the family $\Kc_E$ have no common points, see \rf{TE-D} and \rf{A1-D}, the cubes of the family $\Qc$ are pairwise disjoint. Furthermore, these cubes satisfy inequality \rf{WKE-S}.
%----------------------------------------------------------
\par Let $Q\in\Qc$. We define two cubes $Q',Q''\in\Qc,$ $Q'\cup Q''\subset \gamma Q,$ as follows. If $Q=K\in \tAc$, we put
%----------------------------------------------------------
\bel{QQP1}
Q':=K^{(a_K)},~Q'':=K^{(a_{\tK})}.
\ee
%----------------------------------------------------------
Then, by definition \rf{DEF-TF1},
%----------------------------------------------------------
$$
\tf_1(a_K):=\av{Q'},~~~~
\tf_1(a_{\tK}):=\av{Q''}.
$$
%----------------------------------------------------------
Furthermore, by Lemma \reff{TQ}, $Q'\cup Q''\subset \gamma Q$
with $\gamma=22\tau^2$.
%----------------------------------------------------------
\par If $Q\in\Kc_E$, i.e., $Q=K^{(x)}$ for some $x\in E$, we put $Q'=Q'':=Q$. Clearly in this case $Q'\cup Q''\subset \gamma Q$ with $\gamma=1$, and $\av{Q'}=\av{Q''}.$
%----------------------------------------------------------
\par By these equalities,
%----------------------------------------------------------
$$
I:=\smed_{K\in\, \tAc}\,
\frac{|\tilde{f}_1(a_{\tK})-\tilde{f}_1(a_K)|^p}
{(\diam K)^{p-n}}=
\smed_{Q\in\tAc}\,
\frac{|\av{Q'}-\av{Q''}|^p}
{(\diam Q)^{p-n}}+\smed_{Q\in\Kc_E}\,
\frac{|\av{Q'}-\av{Q''}|^p}
{(\diam Q)^{p-n}},
$$
%----------------------------------------------------------
so that
%----------------------------------------------------------
\bel{IS-L}
I=\smed_{Q\in\Qc}\,
\frac{|\av{Q'}-\av{Q''}|^p}
{(\diam Q)^{p-n}}.
\ee
%----------------------------------------------------------
Since $Q',Q''\in \Kc_E$ for every $Q\in \Qc$, by  \rf{KXM}, %----------------------------------------------------------
$$
\mu(Q')\ge 2^{p-n}(\diam Q')^{n-p}, ~~~~\mu(Q'')\ge 2^{p-n}(\diam Q'')^{n-p}
$$
%----------------------------------------------------------
proving that
%----------------------------------------------------------
\bel{ER}
\frac{1}{\mu(Q')\mu(Q'')}\le
2^{2(n-p)}(\diam Q'\diam Q'')^{p-n} \le (\diam Q'\diam Q'')^{p-n}.
\ee
%----------------------------------------------------------
Also
%----------------------------------------------------------
\bel{F-AVER}
|\av{Q'}-\av{Q''}|^p\le
\frac{1}{\mu(Q')\mu(Q'')}\iint \limits_{Q'\times Q''}
|f(x)-f(y)|^p\, d\mu(x)d\mu(y).
\ee
%----------------------------------------------------------
Combining this inequality with \rf{IS-L} and \rf{ER}, we obtain
%----------------------------------------------------------
$$
I\le \sbig_{Q\in\Qc}\,\,
\left(\frac{\diam Q' \diam Q''}{\diam Q}\right)^{p-n} \iint \limits_{Q'\times Q''}
|f(x)-f(y)|^p\, d\mu(x)d\mu(y).
$$
%----------------------------------------------------------
By the assumption (see inequality \rf{CR-V2}) and in view of Corollary \reff{P-WKE},  $I\le \lambda$ proving inequality \rf{S-TA}.
%----------------------------------------------------------
\par Since all the terms of the sum in the left hand side of \rf{S-TA} are non-negative, this inequality holds for
an {\it arbitrary} (not necessarily finite) subfamily $\tAc$ of $\Ac_1$ consisting of pairwise disjoint cubes.
%----------------------------------------------------------
\par To prove inequality \rf{S-TA} for the family $\Ac_1$
itself (and consequently to prove the proposition) it remains to make use of Proposition \reff{C-PT}. In fact,  every Whitney cube touches at most $N(n)$ Whitney cubes, see part (2) of Lemma \reff{Wadd}. Since $\Ac_1\subset W_E$, the same is true for cubes of the family $\Ac_1$. Hence, by Proposition \reff{C-PT}, $\Ac_1$ can be partitioned into at most $N(n)+1$ families of pairwise disjoint cubes. Applying to every such a family inequality \rf{S-TA} we obtain the required estimate \rf{RE-A1}.
%----------------------------------------------------------
\par The proposition is proved. \bx\bigskip
%----------------------------------------------------------
%@@@@@@@@@@@@@@@@@@@@@@@@@@@@@@@@@@@@@@@@@@@@@@@@@@@@@@@@@@
%@@@@@@@@@@@@@@@@@@@@@@@@@@@@@@@@@@@@@@@@@@@@@@@@@@@@@@@@@@
%@@@@@@@@@@@@@@@@@@@@@@@@@@@@@@@@@@@@@@@@@@@@@@@@@@@@@@@@@@
%@@@@@@@@@@@@@@@@@@@@@@@@@@@@@@@@@@@@@@@@@@@@@@@@@@@@@@@@@@
%----------------------------------------------------------
\par Let us prove an analog of Proposition \reff{EA1} for the family $\Ac_2$. Recall that
%----------------------------------------------------------
$$
\Ac_2=\{Q\in W_E:Q\cap K\ne\emp~~\text{for some}~K\in \Kc_E,~\text{and}~~Q\cap(\eta H)=\emp~~\text{for every}~H\in \Kc_E\}.
$$
%----------------------------------------------------------
%@@@@@@@@@@@@@@@@@@@@@@@@@@@@@@@@@@@@@@@@@@@@@@@@@@@@@@@@@@
%@@@@@@@@@@@@@@@@@@@@@@@@@@@@@@@@@@@@@@@@@@@@@@@@@@@@@@@@@@
\par Let us fix a cube $K\in\Kc_E$ and consider a family of cubes
%----------------------------------------------------------
$$
J_K:=\{Q\in \Ac_2:Q\cap K\ne \emptyset\}.
$$
%----------------------------------------------------------
Thus
%----------------------------------------------------------
$$
J_K:=\{Q\in W_E:Q\cap K\ne \emptyset,~Q\cap(\eta H)=\emp~~\text{for every}~H\in \Kc_E\}.
$$
%----------------------------------------------------------
%@@@@@@@@@@@@@@@@@@@@@@@@@@@@@@@@@@@@@@@@@@@@@@@@@@@@@@@@@@
%@@@@@@@@@@@@@@@@@@@@@@@@@@@@@@@@@@@@@@@@@@@@@@@@@@@@@@@@@@
%@@@@@@@@@@@@@@@@@@@@@@@@@@@@@@@@@@@@@@@@@@@@@@@@@@@@@@@@@@
%@@@@@@@@@@@@@@@@@@@@@@@@@@@@@@@@@@@@@@@@@@@@@@@@@@@@@@@@@@
%----------------------------------------------------------
\begin{lemma}\lbl{K-IS} (i). If $K,K'\in\Kc_E$ and $K\ne K'$, then $(5K)\cap (5K')=\emp$;
%----------------------------------------------------------
\medskip
%----------------------------------------------------------
%@@@@@@@@@@@@@@@@@@@@@@@@@@@@@@@@@@@@@@@@@@@@@@@@@@@@@@@@@@
\par (ii). Let $Q\in W_E, K\in\Kc_E,$ and let $Q\cap K\ne\emp.$ Then $\diam Q\le 2\diam K$ and $Q\subset 5K$;
%----------------------------------------------------------
\medskip
%----------------------------------------------------------
%@@@@@@@@@@@@@@@@@@@@@@@@@@@@@@@@@@@@@@@@@@@@@@@@@@@@@@@@@@
\par (iii). Let $K,K'\in\Kc_E$ and let $Q\in W_E$. Suppose that $K\cap Q\ne\emp$ and $K'\cap Q\ne\emp$. Then $K=K'$.
%----------------------------------------------------------
%@@@@@@@@@@@@@@@@@@@@@@@@@@@@@@@@@@@@@@@@@@@@@@@@@@@@@@@@@@
\end{lemma}
%----------------------------------------------------------
%@@@@@@@@@@@@@@@@@@@@@@@@@@@@@@@@@@@@@@@@@@@@@@@@@@@@@@@@@@
\par {\it Proof.} (i). This property immediately follows from part (i) of Proposition \reff{RPR}, see \rf{DRS}.
%----------------------------------------------------------
\medskip
%----------------------------------------------------------
%@@@@@@@@@@@@@@@@@@@@@@@@@@@@@@@@@@@@@@@@@@@@@@@@@@@@@@@@@@
%----------------------------------------------------------
\par (ii). Let $K=K^{(x)}$ for some $x\in E$. Since $Q\cap K\ne\emp$, we have $\dist(Q,E)\le r_K$. Since $Q\in W_E$, we have $\diam Q\le 4\dist(Q,E)$ so that
%----------------------------------------------------------
$$
\diam Q\le 4\dist(Q,E)\le 4r_K=2\diam K.
$$
%----------------------------------------------------------
\par  Now let $z\in Q$. Since $Q\cap K\ne\emp$,
%----------------------------------------------------------
$$
\|z-x\|\le r_K+\diam Q\le r_K+4r_K=5r_K
$$
%----------------------------------------------------------
proving the required inclusion $Q\subset 5K$.
%----------------------------------------------------------
\medskip
%----------------------------------------------------------
\par (iii). If $K\ne K'$, then, by part (ii) of the lemma
$Q\subset 5K$ and $Q\subset 5K'$ so that $(5K)\cap (5K')\ne\emp$ which contradicts the property (i).
%@@@@@@@@@@@@@@@@@@@@@@@@@@@@@@@@@@@@@@@@@@@@@@@@@@@@@@@@@@
\par The lemma is proved.\bx\medskip
%----------------------------------------------------------
%@@@@@@@@@@@@@@@@@@@@@@@@@@@@@@@@@@@@@@@@@@@@@@@@@@@@@@@@@@
%@@@@@@@@@@@@@@@@@@@@@@@@@@@@@@@@@@@@@@@@@@@@@@@@@@@@@@@@@@
%@@@@@@@@@@@@@@@@@@@@@@@@@@@@@@@@@@@@@@@@@@@@@@@@@@@@@@@@@@
%@@@@@@@@@@@@@@@@@@@@@@@@@@@@@@@@@@@@@@@@@@@@@@@@@@@@@@@@@@
%----------------------------------------------------------
\begin{corollary}\lbl{JK-S} (i). For every $K\in\Kc_E$ we have
%----------------------------------------------------------
$$
\bigcup_{Q\in J_K}\, Q\subset \,5K;
$$
%----------------------------------------------------------
%@@@@@@@@@@@@@@@@@@@@@@@@@@@@@@@@@@@@@@@@@@@@@@@@@@@@@@@@@@
\par (ii). $J_K\cap J_{K'}=\emp$ provided $K,K'\in\Kc_E$ and $K\ne K'$.
%----------------------------------------------------------
%@@@@@@@@@@@@@@@@@@@@@@@@@@@@@@@@@@@@@@@@@@@@@@@@@@@@@@@@@@
\end{corollary}
%----------------------------------------------------------
%@@@@@@@@@@@@@@@@@@@@@@@@@@@@@@@@@@@@@@@@@@@@@@@@@@@@@@@@@@
%@@@@@@@@@@@@@@@@@@@@@@@@@@@@@@@@@@@@@@@@@@@@@@@@@@@@@@@@@@
%@@@@@@@@@@@@@@@@@@@@@@@@@@@@@@@@@@@@@@@@@@@@@@@@@@@@@@@@@@
%----------------------------------------------------------
\medskip
%----------------------------------------------------------
\begin{lemma} For every cube $K\in \Ac_2$
%----------------------------------------------------------
\bel{NJ}
\# J_K\le (45/\eta)^n.
\ee
%----------------------------------------------------------
\par In addition for every cube $Q\in J_K$ %@@@@@@@@@@@@@@@@@@@@@@@@@@@@@@@@@@@@@@@@@@@@@@@@@@@@@@@@@@
%----------------------------------------------------------
\bel{CQK}
\tfrac12\diam Q\le\diam K \le 9\eta\diam Q.
\ee
%----------------------------------------------------------
%@@@@@@@@@@@@@@@@@@@@@@@@@@@@@@@@@@@@@@@@@@@@@@@@@@@@@@@@@@
\end{lemma}
%@@@@@@@@@@@@@@@@@@@@@@@@@@@@@@@@@@@@@@@@@@@@@@@@@@@@@@@@@@
\par {\it Proof.} Note that, by definition of the family $J_K$, for every $Q\in J_K$ we have
%----------------------------------------------------------
\bel{IE1}
Q\cap K\ne \emptyset~~~\text{and}~~~Q\cap(\eta K)=\emp.
\ee
%----------------------------------------------------------
\par First prove inequality \rf{CQK}.
The first inequality in \rf{CQK} follows from part (ii) of Lemma \reff{K-IS}.
%----------------------------------------------------------
\par Let us prove the second one. By \rf{9QE}, $(9Q)\cap E\ne\emp$. Thus there exists a cube $K'=Q(x_{K'},r_{K'})$ such that $x_{K'}\in 9Q$.
%----------------------------------------------------------
\par Let us consider two cases. First suppose that $x_{K'}\ne x_K.$ Then $K'\ne K$. But
%----------------------------------------------------------
$$
(9Q)\cap K'\ne\emp~~~\text{and}~~~(9Q)\cap K\ne\emp
$$
%----------------------------------------------------------
so that
%----------------------------------------------------------
$$
\dist(K,K')\le\diam(9Q)=9\diam Q.
$$
%----------------------------------------------------------
By part (i) of Lemma \reff{PR-KE},
%----------------------------------------------------------
$$
\diam K +\diam K'\le \dist(K,K')
$$
%----------------------------------------------------------
so that
%----------------------------------------------------------
$$
\diam K \le 9\diam Q.
$$
%----------------------------------------------------------
\par Now suppose that $x_{K'}=x_K$ so that $x_K\in 9Q$.
Hence $\|x_K-x_Q\|\le 9r_Q$. On the other hand $Q\cap(\eta K)=\emp$, see \rf{IE1}, so that
%----------------------------------------------------------
$$
\|x_K-x_Q\|\ge \eta r_K.
$$
%----------------------------------------------------------
Hence $r_K\le 9\eta r_Q.$
%----------------------------------------------------------
\par We have proved that in the both cases
$\diam K \le 9\eta\diam Q$ so that the second inequality in \rf{CQK} is satisfied.
%----------------------------------------------------------
\par Let us prove \rf{NJ}. By part (i) of Corollary \reff{JK-S}, $Q\subset 5K$ for every $Q\in J_K.$ Since the cubes of the family $J_K$ are non-overlapping, this inclusion and the second inequality in \rf{CQK} enable us to estimate the cardinality of $J_K$. We have
%----------------------------------------------------------
$$
\# J_K\le |5K|/\min\{|Q|:Q\in J_K\}\le 5^n|K|/(\eta\diam K/9)^n=(45/\eta)^n
$$
%----------------------------------------------------------
proving \rf{NJ} and the lemma.\bx\medskip
%----------------------------------------------------------
%@@@@@@@@@@@@@@@@@@@@@@@@@@@@@@@@@@@@@@@@@@@@@@@@@@@@@@@@@@
%@@@@@@@@@@@@@@@@@@@@@@@@@@@@@@@@@@@@@@@@@@@@@@@@@@@@@@@@@@
%@@@@@@@@@@@@@@@@@@@@@@@@@@@@@@@@@@@@@@@@@@@@@@@@@@@@@@@@@@
%@@@@@@@@@@@@@@@@@@@@@@@@@@@@@@@@@@@@@@@@@@@@@@@@@@@@@@@@@@
%----------------------------------------------------------
\begin{proposition}\lbl{EA2} Suppose that the hypothesis  of Theorem \reff{S-V2} holds with $\gamma=2^{8}\tau^2$. Then %@@@@@@@@@@@@@@@@@@@@@@@@@@@@@@@@@@@@@@@@@@@@@@@@@@@@@@@@@@
%----------------------------------------------------------
$$
\smed_{K\in \Ac_2}\,
\frac{|\tilde{f}_1(a_{\tK})-\tilde{f}_1(a_K)|^p}
{(\diam K)^{p-n}}\le C(n,\tau)\lambda\,.
$$
%----------------------------------------------------------
%@@@@@@@@@@@@@@@@@@@@@@@@@@@@@@@@@@@@@@@@@@@@@@@@@@@@@@@@@@
\end{proposition}
%@@@@@@@@@@@@@@@@@@@@@@@@@@@@@@@@@@@@@@@@@@@@@@@@@@@@@@@@@@
%----------------------------------------------------------
\par {\it Proof.} Following the same scheme of the proof as in Proposition \reff{EA1}, without loss of generality we may assume that $\Ac_2$ is a finite collection of pairwise disjoint cubes. However, for the family $\Ac_2$ we can not define the family $\Qc$ by the same formula as for the family $\Ac_1$, i.e., to put $\Qc:=\Ac_2\cup\Kc_E$. In fact, in this case the cubes of the family $\Ac_2$ intersect cubes of $\Kc_E$ so that the cubes of $\Qc:=\Ac_2\cup\Kc_E$ are not pairwise disjoint.
%----------------------------------------------------------
\par We modify the definition of $\Qc$ as follows. Let
%----------------------------------------------------------
$$
I_K:=\smed_{Q\in J_K}\,
\frac{|\tilde{f}_1(a_{\tQ})-\tilde{f}_1(a_Q)|^p}
{(\diam Q)^{p-n}}
$$
%----------------------------------------------------------
so that
%----------------------------------------------------------
$$
I:=\smed_{K\in \Ac_2}\,
\frac{|\tilde{f}_1(a_{\tK})-\tilde{f}_1(a_K)|^p}
{(\diam K)^{p-n}}=\sum_{K\in \Kc_E}I_K.
$$
%----------------------------------------------------------
Let $H_K\in J_K$ be a cube such that
%----------------------------------------------------------
$$
\max_{Q\in J_K}
\frac{|\tilde{f}_1(a_{\tQ})-\tilde{f}_1(a_Q)|^p}
{(\diam Q)^{p-n}}=
\frac{|\tilde{f}_1(a_{\tH_K})-\tilde{f}_1(a_{H_K})|^p}
{(\diam H_K)^{p-n}}.
$$
%----------------------------------------------------------
(Of course, $H_K$ depends on $f$ as well.) Then, by \rf{NJ},
%----------------------------------------------------------
$$
I_K\le \# J_K\, \frac{|\tilde{f}_1(a_{\tH_K})-\tilde{f}_1(a_{H_K})|^p}
{(\diam H_K)^{p-n}}\le (27/\eta)^n\,
\frac{|\tilde{f}_1(a_{\tH_K})-\tilde{f}_1(a_{H_K})|^p}
{(\diam H_K)^{p-n}}.
$$
%----------------------------------------------------------
\par Let
%----------------------------------------------------------
$$
\hK:=K+3\,r_K\vec{e},~~~\text{where}~~~\vec{e}:=(1,0,...,0).
$$
%----------------------------------------------------------
Then  $\hK\subset 5K\setminus K$, and $5K\subset 7\,\hK$. Since $H_K\in J_K$, by part (i) of Corollary \reff{JK-S}, $H_K\subset 5 K$ so that $H_K\subset 7\,\hK$. Also, by \rf{CQK},
%@@@@@@@@@@@@@@@@@@@@@@@@@@@@@@@@@@@@@@@@@@@@@@@@@@@@@@@@@@
%----------------------------------------------------------
\bel{E-KE}
\diam \hK/(9\eta)=\diam K/(9\eta) \le \diam H_K.
\ee
%----------------------------------------------------------
%@@@@@@@@@@@@@@@@@@@@@@@@@@@@@@@@@@@@@@@@@@@@@@@@@@@@@@@@@@
\par Now we have
%----------------------------------------------------------
$$
I=\sum_{K\in \Kc_E}I_K
\le C\,\smed_{K\in \Kc_E}\,
\frac{|\tilde{f}_1(a_{\tH_K})-\tilde{f}_1(a_{H_K})|^p}
{(\diam H_K)^{p-n}}
$$
%----------------------------------------------------------
so that, by \rf{E-KE},
%----------------------------------------------------------
\bel{S-I5}
I\le C\,\smed_{K\in \Kc_E}\,
\frac{|\tilde{f}_1(a_{\tH_K})-\tilde{f}_1(a_{H_K})|^p}
{(\diam \hK)^{p-n}}
\ee
%----------------------------------------------------------
with $C=C(n,\tau).$
%----------------------------------------------------------
\par We introduce a family of cubes
%----------------------------------------------------------
$$
\Qc:=\left(\bigcup_{K\in \Kc_E}\hK\right)\bigcup \Kc_E.
$$
%----------------------------------------------------------
\par Since $\hK\subset 5K\setminus K$ for every $K\in\Kc_E$, and the cubes $\{5K:K\in \Kc_E\}$ are pairwise disjoint, see part (i) of Lemma \reff{K-IS},
the family $\Qc$ consists of pairwise disjoint cubes.
%----------------------------------------------------------
\par We are in a position to finish the proof of the proposition. Let $Q\in\Qc\setminus \Kc_E$, i.e., there exists a cube $K\in \Kc_E$ such that $Q=\hK$. We put
%----------------------------------------------------------
$$
Q':=K^{(a_{H_K})},~~~Q'':=K^{(a_{\tH_Q})}.
$$
%----------------------------------------------------------
Then, by definition,
%----------------------------------------------------------
$$
\tf_1(a_{H_K})=\av{Q'},~~~
\tf_1(a_{\tH_Q})=\av{Q''}.
$$
%----------------------------------------------------------
Furthermore, by Lemma  \reff{TQ},
%----------------------------------------------------------
$$
Q'\cup Q''\subset (22\tau^2)\, H_K.
$$
%----------------------------------------------------------
%@@@@@@@@@@@@@@@@@@@@@@@@@@@@@@@@@@@@@@@@@@@@@@@@@@@@@@@@@@
But  $H_K\subset 7\,\hK=7\,Q$ so that
%----------------------------------------------------------
$$
Q'\cup Q''\subset 7\,(22\tau^2)\, Q\subset 2^8\tau^2\, Q.
$$
%----------------------------------------------------------  %@@@@@@@@@@@@@@@@@@@@@@@@@@@@@@@@@@@@@@@@@@@@@@@@@@@@@@@@@@
%----------------------------------------------------------
\par If $Q\in\Kc_E$, we put $Q'=Q''=Q$. Hence, by \rf{S-I5},
%----------------------------------------------------------
\bel{FIN-I}
I\le C\,\smed_{Q\in\Qc}\,\frac{|\av{Q'}-\av{Q''}|^p}
{(\diam Q)^{p-n}}.
\ee
%----------------------------------------------------------
\par Note that the cubes $Q'$ and $Q''$ satisfy inequalities \rf{ER} and \rf{F-AVER} of Proposition \reff{EA1}. Combining these inequalities with \rf{FIN-I}, we obtain
%----------------------------------------------------------
$$
I\le C\,\sbig_{Q\in\Qc}\,\,
\left(\frac{\diam Q' \diam Q''}{\diam Q}\right)^{p-n} \iint \limits_{Q'\times Q''}
|f(x)-f(y)|^p\, d\mu(x)d\mu(y),
$$
%----------------------------------------------------------
so that, by the assumption (see inequality \rf{CR-V2}) and in view of Corollary \reff{P-WKE}, we have $I\le C\lambda$. The proposition is completely proved. \bx\medskip
%----------------------------------------------------------
%@@@@@@@@@@@@@@@@@@@@@@@@@@@@@@@@@@@@@@@@@@@@@@@@@@@@@@@@@@
%@@@@@@@@@@@@@@@@@@@@@@@@@@@@@@@@@@@@@@@@@@@@@@@@@@@@@@@@@@
%@@@@@@@@@@@@@@@@@@@@@@@@@@@@@@@@@@@@@@@@@@@@@@@@@@@@@@@@@@
%@@@@@@@@@@@@@@@@@@@@@@@@@@@@@@@@@@@@@@@@@@@@@@@@@@@@@@@@@@
%@@@@@@@@@@@@@@@@@@@@@@@@@@@@@@@@@@@@@@@@@@@@@@@@@@@@@@@@@@
%@@@@@@@@@@@@@@@@@@@@@@@@@@@@@@@@@@@@@@@@@@@@@@@@@@@@@@@@@@
%----------------------------------------------------------
\par This proposition and Proposition \reff{EA2} imply the following
%----------------------------------------------------------
%@@@@@@@@@@@@@@@@@@@@@@@@@@@@@@@@@@@@@@@@@@@@@@@@@@@@@@@@@@
%@@@@@@@@@@@@@@@@@@@@@@@@@@@@@@@@@@@@@@@@@@@@@@@@@@@@@@@@@@
%@@@@@@@@@@@@@@@@@@@@@@@@@@@@@@@@@@@@@@@@@@@@@@@@@@@@@@@@@@
\begin{corollary}\lbl{E-A12} The following inequality
%----------------------------------------------------------
$$
\smed_{K\in \Ac}\,
\frac{|\tilde{f}_1(a_{\tK})-\tilde{f}_1(a_K)|^p}
{(\diam K)^{p-n}}\le C(n,\tau)\lambda
$$
%----------------------------------------------------------
is satisfied provided the hypothesis of Theorem \reff{S-V2}  holds with $\gamma= 2^8\tau^2$.
%----------------------------------------------------------
\end{corollary}
%@@@@@@@@@@@@@@@@@@@@@@@@@@@@@@@@@@@@@@@@@@@@@@@@@@@@@@@@@@
\par {\it Proof.} Recall that
%----------------------------------------------------------
$$
\smed_{K\in \Ac}\,
\frac{|\tilde{f}_1(a_{\tK})-\tilde{f}_1(a_K)|^p}
{(\diam K)^{p-n}}\le
C\smed_{K\in \Ac_1}\,
\frac{|\tilde{f}_1(a_{\tK})-\tilde{f}_1(a_K)|^p}
{(\diam K)^{p-n}}+C\smed_{K\in \Ac_2}\,
\frac{|\tilde{f}_1(a_{\tK})-\tilde{f}_1(a_K)|^p}
{(\diam K)^{p-n}}.
$$
%----------------------------------------------------------
It remains to apply Proposition \reff{EA1} to the first sum, and Proposition \reff{EA2} to the second sum in the right hand side of this inequality, and the corollary follows.\bx\bigskip
%----------------------------------------------------------
%@@@@@@@@@@@@@@@@@@@@@@@@@@@@@@@@@@@@@@@@@@@@@@@@@@@@@@@@@@
%@@@@@@@@@@@@@@@@@@@@@@@@@@@@@@@@@@@@@@@@@@@@@@@@@@@@@@@@@@
%@@@@@@@@@@@@@@@@@@@@@@@@@@@@@@@@@@@@@@@@@@@@@@@@@@@@@@@@@@
%@@@@@@@@@@@@@@@@@@@@@@@@@@@@@@@@@@@@@@@@@@@@@@@@@@@@@@@@@@
%----------------------------------------------------------
\par Finally, combining this corollary with Corollary \reff{NF1-1} we obtain the required inequality
%----------------------------------------------------------
\bel{F1-N}
\|\nabla f_1\|^p_{\LPRN}\le C\,\lambda
\ee
%----------------------------------------------------------
provided $\gamma=2^8\tau^2$ and $C=C(n,\tau)$.
%----------------------------------------------------------
%@@@@@@@@@@@@@@@@@@@@@@@@@@@@@@@@@@@@@@@@@@@@@@@@@@@@@@@@@@
%@@@@@@@@@@@@@@@@@@@@@@@@@@@@@@@@@@@@@@@@@@@@@@@@@@@@@@@@@@
%@@@@@@@@@@@@@@@@@@@@@@@@@@@@@@@@@@@@@@@@@@@@@@@@@@@@@@@@@@
%@@@@@@@@@@@@@@@@@@@@@@@@@      @@@@@@@@@@@@@@@@@@@@@@@@@@@
%@@@@@@@@@@@@@@@@@@@@@@@          @@@@@@@@@@@@@@@@@@@@@@@@@
%@@@@@@@@@@@@@@@@@@@@@              @@@@@@@@@@@@@@@@@@@@@@@
%@@@@@@@@@@@@@@@@@@@     SECTION 5    @@@@@@@@@@@@@@@@@@@@@
%@@@@@@@@@@@@@@@@@@@@@              @@@@@@@@@@@@@@@@@@@@@@@
%@@@@@@@@@@@@@@@@@@@@@@@          @@@@@@@@@@@@@@@@@@@@@@@@@
%@@@@@@@@@@@@@@@@@@@@@@@@@      @@@@@@@@@@@@@@@@@@@@@@@@@@@
%@@@@@@@@@@@@@@@@@@@@@@@@@@@@@@@@@@@@@@@@@@@@@@@@@@@@@@@@@@
%@@@@@@@@@@@@@@@@@@@@@@@@@@@@@@@@@@@@@@@@@@@@@@@@@@@@@@@@@@
%----------------------------------------------------------
%@@@@@@@@@@@@@@@@@@@@@@@@@@@@@@@@@@@@@@@@@@@@@@@@@@@@@@@@@@
%----------------------------------------------------------
\SECT{5. Sufficiency: the $L_p(\RN;\mu)$-norm of the function $f_2$.}{5}
%----------------------------------------------------------
\addtocontents{toc}{5. Sufficiency: the $L_p(\RN;\mu)$-norm of the function $f_2$. \hfill \thepage\\\par}
%----------------------------------------------------------
\indent
%----------------------------------------------------------
%@@@@@@@@@@@@@@@@@@@@@@@@@@@@@@@@@@@@@@@@@@@@@@@@@@@@@@@@@@
%@@@@@@@@@@@@@@@@@@@@@@@@@@@@@@@@@@@@@@@@@@@@@@@@@@@@@@@@@@
%@@@@@@@@@@@@@@@@@@@@@@@@@@@@@@@@@@@@@@@@@@@@@@@@@@@@@@@@@@
\par In this section we prove that
%----------------------------------------------------------
\bel{F2-N}
\|f_2\|^p_{\LPM}\le C\,\lambda.
\ee
%----------------------------------------------------------
%@@@@@@@@@@@@@@@@@@@@@@@@@@@@@@@@@@@@@@@@@@@@@@@@@@@@@@@@@@
%@@@@@@@@@@@@@@@@@@@@@@@@@@@@@@@@@@@@@@@@@@@@@@@@@@@@@@@@@@
%@@@@@@@@@@@@@@@@@@@@@@@@@@@@@@@@@@@@@@@@@@@@@@@@@@@@@@@@@@
\begin{lemma}\lbl{N-K} For every cube $K\in W_E$ and every $c\in\R$ the following inequality
%----------------------------------------------------------
$$
\intl_K |f_1(x)-c|^pd\mu(x)\le \left(\,\sum_{Q\in V_K}|\tilde{f}_1(a_Q)-c|\right)^p\mu(K)
$$
%----------------------------------------------------------
holds.
\end{lemma}
%@@@@@@@@@@@@@@@@@@@@@@@@@@@@@@@@@@@@@@@@@@@@@@@@@@@@@@@@@@
\par {\it Proof.} By the extension formula \rf{DEF-F1},
%----------------------------------------------------------
\be
\intl_K |f_1(x)-c|^p\,d\mu&=&\intl_K |\sum_{Q\in W_E}\varphi_Q(x)\tilde{f}_1(a_Q)-c|^p\,d\mu\nn\\
&=&
\intl_K |\sum_{Q\in W_E}\varphi_Q(x)(\tilde{f}_1(a_Q)-c)|^p\,d\mu\nn\\
&=&
\intl_K\left|\sum\left\{\varphi_Q(x)
(\tilde{f}_1(a_Q)-c):Q^*\cap K\ne \emptyset, Q\in W_E\right\}\right|^p\,d\mu\nn\\
&\le&
\intl_K\left(\sum\left\{|\varphi_Q(x)|\,
|\tilde{f}_1(a_Q)-c|:Q^*\cap K\ne \emptyset, Q\in W_E\right\}\right)^p\,d\mu.\nn
%----------------------------------------------------------
\ee
%----------------------------------------------------------
\par By part (a) of Lemma \reff{P-U} and part (3) of Lemma \reff{Wadd},
%----------------------------------------------------------
\be
\intl_K |f_1(x)-c|^p\,d\mu
&\le&
\intl_K\left(\sum
\left\{|\tilde{f}_1(a_Q)-c|:Q^*\cap K\ne \emptyset, Q\in W_E\right\}\right)^p\,d\mu\nn\\
&=&
\mu(K)\left(\sum
\left\{|\tilde{f}_1(a_Q)-c|:Q\cap K\ne \emptyset, Q\in W_E\right\}\right)^p\nn\\
&=&
\mu(K)\left(\,\sum_{Q\in V_K}
|\tilde{f}_1(a_Q)-c|\right)^p
\nn
%----------------------------------------------------------
\ee
%@@@@@@@@@@@@@@@@@@@@@@@@@@@@@@@@@@@@@@@@@@@@@@@@@@@@@@@@@@@
proving the lemma.\bx\bigskip
%----------------------------------------------------------
%@@@@@@@@@@@@@@@@@@@@@@@@@@@@@@@@@@@@@@@@@@@@@@@@@@@@@@@@@@
%@@@@@@@@@@@@@@@@@@@@@@@@@@@@@@@@@@@@@@@@@@@@@@@@@@@@@@@@@@
%@@@@@@@@@@@@@@@@@@@@@@@@@@@@@@@@@@@@@@@@@@@@@@@@@@@@@@@@@@
%@@@@@@@@@@@@@@@@@@@@@@@@@@@@@@@@@@@@@@@@@@@@@@@@@@@@@@@@@@
%----------------------------------------------------------
\par Let $K\in \Kc_E$ and let
%----------------------------------------------------------
\bel{D-SK}
S_K:=(\eta K)\bigcup \{Q:Q\in W_E,\,\, Q\cap(\eta K)\ne\emp\}.
\ee
%----------------------------------------------------------
(Recall that $\eta=\tfrac{1}{21\tau}$, see \rf{ETA}.)
%@@@@@@@@@@@@@@@@@@@@@@@@@@@@@@@@@@@@@@@@@@@@@@@@@@@@@@@@@@
%@@@@@@@@@@@@@@@@@@@@@@@@@@@@@@@@@@@@@@@@@@@@@@@@@@@@@@@@@@
%@@@@@@@@@@@@@@@@@@@@@@@@@@@@@@@@@@@@@@@@@@@@@@@@@@@@@@@@@@
%@@@@@@@@@@@@@@@@@@@@@@@@@@@@@@@@@@@@@@@@@@@@@@@@@@@@@@@@@@
\begin{lemma}\lbl{F2-S} For every $K\in \Kc_E$ the following inequality
%----------------------------------------------------------
$$
\intl_{S_K}|f_2|^pd\mu\le \intl_{K}|f-\av{K}|^p\,d\mu
$$
%----------------------------------------------------------
holds.
%----------------------------------------------------------
\end{lemma}
%@@@@@@@@@@@@@@@@@@@@@@@@@@@@@@@@@@@@@@@@@@@@@@@@@@@@@@@@@@
\par {\it Proof.} First prove that
%----------------------------------------------------------
\bel{R-SI}
f_1(y)=\av{K}~~~\text{for every}~~~y\in S_K.
\ee
%----------------------------------------------------------
In fact, since $K\in\Kc_E$, there exists $x\in E$ such that $K=K^{(x)}$. Recall that, by definitions \rf{DEF-TF1} and \rf{DEF-F1},
%----------------------------------------------------------
$$
f_1(x)=\tf_1(x)=\av{K}.
$$
%----------------------------------------------------------
\par Let $y\in S_K\setminus\{x\}$ so there exists a cube $H\in W_E$ such that $y\in H$ and $H\cap(\eta K)\ne\emp$. By formula \rf{DEF-F1} and part (c) of Lemma \reff{P-U},
%----------------------------------------------------------
\be
I&:=&f_1(y)-\tf_1(x)=\left(\,\sum_{Q\in W_E}\varphi_Q(y)\tilde{f}_1(a_Q)\right)-\tf_1(x)\nn\\
&=&
\sum_{Q\in W_E}\varphi_Q(y)(\tilde{f}_1(a_Q)-\tf_1(x))\nn\\
&=&
\sum\{\varphi_Q(y)(\tilde{f}_1(a_Q)-\tf_1(x)):Q\in W_E,~
Q^*\cap H\ne\emp\}.
\nn
%----------------------------------------------------------
\ee
%----------------------------------------------------------
Hence, by part (3) of Lemma \reff{Wadd},
%----------------------------------------------------------
\be
I&=&
\sum\{\varphi_Q(y)(\tilde{f}_1(a_Q)-\tf_1(x)):Q\in W_E,~
Q\cap H\ne\emp\}\nn\\
&=&
\sum_{Q\in V_H}\varphi_Q(y)(\tilde{f}_1(a_Q)-\tf_1(x)).
\nn
%----------------------------------------------------------
\ee
%----------------------------------------------------------
\par Since $H\cap (\eta K)\ne\emp$, by Lemma \reff{SMC}, $a_Q=x$ for every $Q\in V_H$ so that
%----------------------------------------------------------
$$
I:=f_1(y)-\tf_1(x)=f_1(y)-\av{K}=0,
$$
%----------------------------------------------------------
and \rf{R-SI} follows.
%----------------------------------------------------------
\par Furthermore, since $Q\in V_Q$, by Lemma \reff{SMC}, $Q\subset K^{(x)}=K$ for every $Q\in W_E$ such that $Q\cap (\eta K)\ne\emp$. Hence $S_K\subset K$.
%----------------------------------------------------------
\par Finally,
%----------------------------------------------------------
$$
\intl_{S_K}|f_2|^pd\mu=\intl_{S_K}|f(y)-f_1(y)|^pd\mu(y)=
\intl_{S_K}|f(y)-\av{K}|^pd\mu(y)
\le \intl_{K}|f-\av{K}|^p\,d\mu
$$
%----------------------------------------------------------
proving the lemma.\bx\bigskip
%----------------------------------------------------------
%@@@@@@@@@@@@@@@@@@@@@@@@@@@@@@@@@@@@@@@@@@@@@@@@@@@@@@@@@@
%@@@@@@@@@@@@@@@@@@@@@@@@@@@@@@@@@@@@@@@@@@@@@@@@@@@@@@@@@@
%@@@@@@@@@@@@@@@@@@@@@@@@@@@@@@@@@@@@@@@@@@@@@@@@@@@@@@@@@@
%@@@@@@@@@@@@@@@@@@@@@@@@@@@@@@@@@@@@@@@@@@@@@@@@@@@@@@@@@@
%----------------------------------------------------------
\par Let $\Ac$ be the family of cubes defined by \rf{A-D}. Thus $K\in\Ac$ $\Leftrightarrow$ $K\in W_E$ and $K\cap(\eta T)=\emp$ for every cube $T\in\Kc_E$. Our next goal is to prove that
%----------------------------------------------------------
\bel{F-G}
\smed_{K\in\Ac}\, \intl_{K}|f_2|^pd\mu\le C\,\lambda.
\ee
%----------------------------------------------------------
\par We have
%----------------------------------------------------------
$$
\smed_{K\in\Ac}\, \intl_{K}|f_2|^pd\mu
=\smed_{K\in\Ac}\, \intl_{K}|f-f_1|^pd\mu
$$
%----------------------------------------------------------
so that
%----------------------------------------------------------
\bel{F2-R}
\smed_{K\in\Ac}\, \intl_{K}|f_2|^pd\mu
\le 2^p\left\{\smed_{K\in\Ac}\,\intl_K |f_1-\tf_1(a_K)|^p\,d\mu+
\smed_{K\in\Ac}\,\intl_K |f-\tf_1(a_K)|^p\,d\mu \right\}.
\ee
%----------------------------------------------------------
%@@@@@@@@@@@@@@@@@@@@@@@@@@@@@@@@@@@@@@@@@@@@@@@@@@@@@@@@@@
%@@@@@@@@@@@@@@@@@@@@@@@@@@@@@@@@@@@@@@@@@@@@@@@@@@@@@@@@@@
%@@@@@@@@@@@@@@@@@@@@@@@@@@@@@@@@@@@@@@@@@@@@@@@@@@@@@@@@@@
%@@@@@@@@@@@@@@@@@@@@@@@@@@@@@@@@@@@@@@@@@@@@@@@@@@@@@@@@@@
%----------------------------------------------------------
\begin{proposition}\lbl{EA3} Suppose that the hypothesis of Theorem \reff{S-V2} holds with $\gamma=2^8\tau^2$. Then   %@@@@@@@@@@@@@@@@@@@@@@@@@@@@@@@@@@@@@@@@@@@@@@@@@@@@@@@@@@
%----------------------------------------------------------
$$
\smed_{K\in\Ac}\,\intl_K |f_1-\tf_1(a_K)|^p\,d\mu\le C(n,p,\tau)\,\lambda.
$$
%----------------------------------------------------------
%@@@@@@@@@@@@@@@@@@@@@@@@@@@@@@@@@@@@@@@@@@@@@@@@@@@@@@@@@@
\end{proposition}
%@@@@@@@@@@@@@@@@@@@@@@@@@@@@@@@@@@@@@@@@@@@@@@@@@@@@@@@@@@
%----------------------------------------------------------
\par {\it Proof.} By Lemma \reff{N-K},
%----------------------------------------------------------
$$
\intl_K |f_1(x)-\tf_1(a_K)|^p\,d\mu(x)\le \left(\,\sum_{Q\in V_K}|\tilde{f}_1(a_Q)-\tf_1(a_K)|\right)^p\mu(K)
$$
%----------------------------------------------------------
\par Recall that $V_K=\{Q\in W_E:Q\cap Q\ne\emp\}$, and
$\#V_K\le N(n)$, see part (2) of Lemma \reff{Wadd}. We also recall that by $\tK\in V_K$ we denote a cube satisfying \rf{D-KW}. Hence,
%----------------------------------------------------------
\be
\left(\sum_{Q\in V_K} |\tilde{f}_1(a_Q)-\tf_1(a_K)|\right)^p\mu(K)&\le& (\#V_K)^p
\mu(K)|\tilde{f}_1(a_{\tK})-\tf_1(a_K)|^p\nn\\&\le&
C(n)\mu(K)|\tilde{f}_1(a_{\tK})-\tf_1(a_K)|^p.\nn
\ee
%----------------------------------------------------------
\par We obtain
%----------------------------------------------------------
$$
\intl_K |f_1(x)-\tf_1(a_K)|^p\,d\mu(x)\le C(n)\,\mu(K)\,|\tilde{f}_1(a_{\tK})-\tf_1(a_K)|^p.
$$
%----------------------------------------------------------
Since $K\in W_E$, by Lemma \reff{WQ-M}, $\mu(K)\le 2^{15p}(\diam K)^{n-p}$, so that
%----------------------------------------------------------
$$
\intl_K |f_1(x)-\tf_1(a_K)|^p\,d\mu(x)\le C(n,p) \,\frac{|\tilde{f}_1(a_{\tK})-\tf_1(a_K)|^p}{(\diam K)^{p-n}}.
$$
%----------------------------------------------------------
Hence,
%----------------------------------------------------------
\bel{P-A4}
\smed_{K\in\Ac}\,\intl_K |f_1-\tf_1(a_K)|^p\,d\mu\le C(n,p)\smed_{K\in\Ac} \,\frac{|\tilde{f}_1(a_{\tK})-\tf_1(a_K)|^p}{(\diam K)^{p-n}}
\ee
%----------------------------------------------------------
so that, by Corollary \reff{E-A12},
%----------------------------------------------------------
$$
\smed_{K\in\Ac}\,\intl_K |f_1-\tf_1(a_K)|^p\,d\mu\le C(n,p,\tau)\,\lambda
$$
%----------------------------------------------------------
provided the hypothesis of Theorem \reff{S-V2} holds with $\gamma=2^8\tau^2$.\bx\medskip
%----------------------------------------------------------
\par Let
%----------------------------------------------------------
\bel{I2-D}
I_2:=\smed_{K\in\Ac}\,\intl_K |f-\tf_1(a_K)|^p\,d\mu+\smed_{K\in\Kc_E}\, \intl_{K}|f-\av{K}|^p\,d\mu.
\ee
%----------------------------------------------------------
\par We define a collection of cubes
%----------------------------------------------------------
\bel{Q-U}
\Qc:=\Ac\cup \Kc_E.
\ee
%----------------------------------------------------------
As we have noted in the proof of Proposition \reff{EA1}, see a remark after \rf{D-A1-N}, the cubes of $\Qc$ are pairwise disjoint. Furthermore, these cubes satisfy inequality \rf{WKE-S} so that the conditions \rf{G-MD} of Theorem \reff{S-V2} hold.
%----------------------------------------------------------
\par Let $Q\in\Qc$. We define two cubes $Q',Q''\in\Qc$, $Q'\cup Q''\subset \gamma Q$, as follows. If $Q=K\in \Ac$, we put
%----------------------------------------------------------
\bel{RM-2}
Q':=Q=K,~~~Q'':=K^{(a_{K})},
\ee
%----------------------------------------------------------
(c.f. definition \rf{QQP1}). Then, by definition \rf{DEF-TF1},
%----------------------------------------------------------
$$
\tf_1(a_K):=\av{Q''}.
$$
%----------------------------------------------------------
Furthermore, by Lemma \reff{TQ}, $Q'\cup Q''\subset \gamma Q$
with $\gamma=22\tau^2$.
%----------------------------------------------------------
\par If $Q\in\Kc_E$, i.e., $Q=K^{(x)}$ for some $x\in E$, we put
%----------------------------------------------------------
\bel{R1-QP}
Q'=Q'':=Q.
\ee
%----------------------------------------------------------
Thus in this case
%----------------------------------------------------------
$$
Q'\cup Q''\subset \gamma Q~~~\text{ with}~~~\gamma=1.
$$
%----------------------------------------------------------
%@@@@@@@@@@@@@@@@@@@@@@@@@@@@@@@@@@@@@@@@@@@@@@@@@@@@@@@@@@
%@@@@@@@@@@@@@@@@@@@@@@@@@@@@@@@@@@@@@@@@@@@@@@@@@@@@@@@@@@
%@@@@@@@@@@@@@@@@@@@@@@@@@@@@@@@@@@@@@@@@@@@@@@@@@@@@@@@@@@
%@@@@@@@@@@@@@@@@@@@@@@@@@@@@@@@@@@@@@@@@@@@@@@@@@@@@@@@@@@
%----------------------------------------------------------
\begin{proposition}\lbl{W-LI} We have
%----------------------------------------------------------
%@@@@@@@@@@@@@@@@@@@@@@@@@@@@@@@@@@@@@@@@@@@@@@@@@@@@@@@@@@
%----------------------------------------------------------
$$
I_2\le
\smed_{Q\in\Qc}\left(\frac{\diam Q' \diam Q''}{\diam Q}\right)^{p-n} \iint \limits_{Q'\times Q''}
|f(x)-f(y)|^p\, d\mu(x)d\mu(y).
$$
%----------------------------------------------------------
%@@@@@@@@@@@@@@@@@@@@@@@@@@@@@@@@@@@@@@@@@@@@@@@@@@@@@@@@@@
\end{proposition}
%@@@@@@@@@@@@@@@@@@@@@@@@@@@@@@@@@@@@@@@@@@@@@@@@@@@@@@@@@@
%----------------------------------------------------------
\par {\it Proof.} The proof of this proposition is very similar to that of Proposition \reff{EA1}. First, as in that proof, without loss of generality we may assume that $\Ac$ is a finite collection of pairwise disjoint cubes.
%----------------------------------------------------------
%@@@@@@@@@@@@@@@@@@@@@@@@@@@@@@@@@@@@@@@@@@@@@@@@@@@@@@@@@@
%----------------------------------------------------------
\par Note that $Q''\in\Kc_E$ for every $Q\in\Qc$ so that, by Corollary \reff{PR-5K}, see \rf{M-DK},
%----------------------------------------------------------
\bel{DQ2}
(\diam Q'')^{n-p}\le 2^{n-p}\mu(Q'')\le \mu(Q'').
\ee
%----------------------------------------------------------
Also, since  $Q'=Q$  for every $Q\in\Qc$, we have
%----------------------------------------------------------
\bel{QPN}
\diam Q'=\diam Q,~~~Q\in\Qc.
\ee
%----------------------------------------------------------
\par Hence,
%----------------------------------------------------------
$$
I_2:=\smed_{K\in\Ac}\,\intl_K |f-\tf_1(a_K)|^p\,d\mu+\smed_{K\in\Kc_E}\, \intl_{K}|f-\av{K}|^p\,d\mu=
\smed_{Q\in\Qc}\, \intl_{Q'} |f-f_{Q''}|^p\,d\mu.
$$
%----------------------------------------------------------
But
%----------------------------------------------------------
$$
\intl_{Q'} |f-f_{Q''}|^p\,d\mu\le
\frac{1}{\mu(Q'')}\iint \limits_{Q'\times Q''}
|f(x)-f(y)|^p\, d\mu(x)d\mu(y)
$$
%----------------------------------------------------------
so that, by \rf{DQ2},
%----------------------------------------------------------
$$
\intl_{Q'} |f-f_{Q''}|^p\,d\mu\le
(\diam Q'')^{p-n}\iint \limits_{Q'\times Q''}
|f(x)-f(y)|^p\, d\mu(x)d\mu(y).
$$
%----------------------------------------------------------
Hence, by \rf{QPN},
%----------------------------------------------------------
$$
\intl_{Q'} |f-f_{Q''}|^p\,d\mu\le
\left(\frac{\diam Q' \diam Q''}{\diam Q}\right)^{p-n} \iint \limits_{Q'\times Q''}
|f(x)-f(y)|^p\, d\mu(x)d\mu(y)
$$
%----------------------------------------------------------
so that
%----------------------------------------------------------
$$
I_2=\smed_{Q\in\Qc}\, \intl_{Q'} |f-f_{Q''}|^p\,d\mu\le
\smed_{Q\in\Qc}\left(\frac{\diam Q' \diam Q''}{\diam Q}\right)^{p-n} \iint \limits_{Q'\times Q''}
|f(x)-f(y)|^p\, d\mu(x)d\mu(y)
$$
%----------------------------------------------------------
proving the proposition.\bx
%----------------------------------------------------------
%@@@@@@@@@@@@@@@@@@@@@@@@@@@@@@@@@@@@@@@@@@@@@@@@@@@@@@@@@@
%@@@@@@@@@@@@@@@@@@@@@@@@@@@@@@@@@@@@@@@@@@@@@@@@@@@@@@@@@@
%@@@@@@@@@@@@@@@@@@@@@@@@@@@@@@@@@@@@@@@@@@@@@@@@@@@@@@@@@@
%@@@@@@@@@@@@@@@@@@@@@@@@@@@@@@@@@@@@@@@@@@@@@@@@@@@@@@@@@@
%----------------------------------------------------------
\begin{corollary}\lbl{EA4} Suppose that the hypothesis of Theorem \reff{S-V2} holds with $\gamma=22\tau^2$. Then   %@@@@@@@@@@@@@@@@@@@@@@@@@@@@@@@@@@@@@@@@@@@@@@@@@@@@@@@@@@
%----------------------------------------------------------
$$
I_2=\smed_{K\in\Ac}\,\intl_K |f-\tf_1(a_K)|^p\,d\mu+\smed_{K\in\Kc_E}\, \intl_{K}|f-\av{K}|^p\,d\mu\le \,\lambda\,.
$$
%----------------------------------------------------------
%@@@@@@@@@@@@@@@@@@@@@@@@@@@@@@@@@@@@@@@@@@@@@@@@@@@@@@@@@@
\end{corollary}
%@@@@@@@@@@@@@@@@@@@@@@@@@@@@@@@@@@@@@@@@@@@@@@@@@@@@@@@@@@
%----------------------------------------------------------
\par {\it Proof.} Note that the family $\Qc$ and the cubes $Q',Q''$ satisfy all the conditions of Theorem \reff{S-V2}. Then, by the theorem's hypothesis, $I_2\le\lambda$, proving the corollary.\bx\bigskip
%----------------------------------------------------------
%@@@@@@@@@@@@@@@@@@@@@@@@@@@@@@@@@@@@@@@@@@@@@@@@@@@@@@@@@@
%@@@@@@@@@@@@@@@@@@@@@@@@@@@@@@@@@@@@@@@@@@@@@@@@@@@@@@@@@@
%@@@@@@@@@@@@@@@@@@@@@@@@@@@@@@@@@@@@@@@@@@@@@@@@@@@@@@@@@@
%@@@@@@@@@@@@@@@@@@@@@@@@@@@@@@@@@@@@@@@@@@@@@@@@@@@@@@@@@@
%----------------------------------------------------------
\par Combining this corollary with Proposition \reff{EA3} we obtain the required inequality \rf{F-G}. In fact,
%----------------------------------------------------------
\be
\smed_{K\in\Ac}\, \intl_{K}|f_2|^pd\mu
&\le& 2^p\left\{\smed_{K\in\Ac}\,\intl_K |f_1-\tf_1(a_K)|^p\,d\mu+
\smed_{K\in\Ac}\,\intl_K |f-\tf_1(a_K)|^p\,d\mu \right\}\nn\\
&\le& C(n,p,\tau) \lambda\nn
\ee
%----------------------------------------------------------
provided the hypothesis of Theorem \reff{S-V2} hold with $\gamma=2^8\tau^2$.
%----------------------------------------------------------
\par Note that Proposition \reff{EA4} implies the following
%@@@@@@@@@@@@@@@@@@@@@@@@@@@@@@@@@@@@@@@@@@@@@@@@@@@@@@@@@@
\begin{corollary}\lbl{I2-SK} The following inequality
%----------------------------------------------------------
$$
\smed_{K\in\Kc_E}\,\intl_{S_K}|f_2|^pd\mu\le \lambda
$$
%----------------------------------------------------------
holds.
%----------------------------------------------------------
\end{corollary}
%@@@@@@@@@@@@@@@@@@@@@@@@@@@@@@@@@@@@@@@@@@@@@@@@@@@@@@@@@@
%@@@@@@@@@@@@@@@@@@@@@@@@@@@@@@@@@@@@@@@@@@@@@@@@@@@@@@@@@@
%@@@@@@@@@@@@@@@@@@@@@@@@@@@@@@@@@@@@@@@@@@@@@@@@@@@@@@@@@@
%@@@@@@@@@@@@@@@@@@@@@@@@@@@@@@@@@@@@@@@@@@@@@@@@@@@@@@@@@@
%----------------------------------------------------------
\par {\it Proof.}  By Corollary \reff{EA4},
%----------------------------------------------------------
$$
\smed_{K\in\Kc_E}\, \intl_{K}|f-\av{K}|^p\,d\mu\le \,\lambda
$$
%----------------------------------------------------------
so that, by Lemma \reff{F2-S},
%----------------------------------------------------------
$$
\smed_{K\in\Kc_E}\,\intl_{S_K}|f_2|^pd\mu\le \smed_{K\in\Kc_E}\, \intl_{K}|f-\av{K}|^p\,d\mu\le \lambda
$$
%----------------------------------------------------------
proving the corollary.\bx\bigskip
%----------------------------------------------------------
%@@@@@@@@@@@@@@@@@@@@@@@@@@@@@@@@@@@@@@@@@@@@@@@@@@@@@@@@@@
%@@@@@@@@@@@@@@@@@@@@@@@@@@@@@@@@@@@@@@@@@@@@@@@@@@@@@@@@@@
%@@@@@@@@@@@@@@@@@@@@@@@@@@@@@@@@@@@@@@@@@@@@@@@@@@@@@@@@@@
%@@@@@@@@@@@@@@@@@@@@@@@@@@@@@@@@@@@@@@@@@@@@@@@@@@@@@@@@@@
%----------------------------------------------------------
\par To finish the proof of Theorem \reff{S-V2} it remains to note that the collection of subsets $\{S_K:K\in\Kc_E\}$, see \rf{D-SK}, and the family of cubes $\Ac$, see \rf{A-D}, cover $\RN$. Hence
%----------------------------------------------------------
\bel{F-S5}
\intl_{\RN}|f_2|^pd\mu\le \smed_{K\in\Kc_E}\, \intl_{S_K}|f_2|^pd\mu+\smed_{K\in\Ac}\, \intl_{K}|f_2|^pd\mu.
\ee
%----------------------------------------------------------
Combining this inequality with Corollary \reff{I2-SK} and Corollary \reff{F-G} we obtain the required inequality \rf{F2-N} provided $\gamma=2^8\tau^2$ and $C=C(n,\tau)$.
%----------------------------------------------------------
\par In turn, combining \rf{F2-N} with inequality \rf{F1-N}, we finally obtain the required estimate \rf{F12}.
%----------------------------------------------------------
\par Theorem \reff{S-V2} is completely proved.\bx
%----------------------------------------------------------
%@@@@@@@@@@@@@@@@@@@@@@@@@@@@@@@@@@@@@@@@@@@@@@@@@@@@@@@@@@
%----------------------------------------------------------
%@@@@@@@@@@@@@@@@@@@@@@@@@@@@@@@@@@@@@@@@@@@@@@@@@@@@@@@@@@
%@@@@@@@@@@@@@@@@@@@@@@@@@@@@@@@@@@@@@@@@@@@@@@@@@@@@@@@@@@
%@@@@@@@@@@@@@@@@@@@@@@@@@@@@@@@@@@@@@@@@@@@@@@@@@@@@@@@@@@
%@@@@@@@@@@@@@@@@@@@@@@@@@      @@@@@@@@@@@@@@@@@@@@@@@@@@@
%@@@@@@@@@@@@@@@@@@@@@@@          @@@@@@@@@@@@@@@@@@@@@@@@@
%@@@@@@@@@@@@@@@@@@@@@              @@@@@@@@@@@@@@@@@@@@@@@
%@@@@@@@@@@@@@@@@@@@     SECTION 6    @@@@@@@@@@@@@@@@@@@@@
%@@@@@@@@@@@@@@@@@@@@@              @@@@@@@@@@@@@@@@@@@@@@@
%@@@@@@@@@@@@@@@@@@@@@@@          @@@@@@@@@@@@@@@@@@@@@@@@@
%@@@@@@@@@@@@@@@@@@@@@@@@@      @@@@@@@@@@@@@@@@@@@@@@@@@@@
%@@@@@@@@@@@@@@@@@@@@@@@@@@@@@@@@@@@@@@@@@@@@@@@@@@@@@@@@@@
%@@@@@@@@@@@@@@@@@@@@@@@@@@@@@@@@@@@@@@@@@@@@@@@@@@@@@@@@@@
%----------------------------------------------------------
%@@@@@@@@@@@@@@@@@@@@@@@@@@@@@@@@@@@@@@@@@@@@@@@@@@@@@@@@@@
%----------------------------------------------------------
\SECT{6. Refinements of the criterion for the norm in $\SUM$.}{6}
%----------------------------------------------------------
\addtocontents{toc}{6. Refinements of the criterion for the norm in $\SUM$. \hfill \thepage\vspace*{3mm}\par}
%----------------------------------------------------------
\indent
%----------------------------------------------------------
%@@@@@@@@@@@@@@@@@@@@@@@@@@@@@@@@@@@@@@@@@@@@@@@@@@@@@@@@@@
%----------------------------------------------------------
\par {\bf 6.1 A refinement of Theorem \reff{MAIN-CR}.}
%----------------------------------------------------------
\addtocontents{toc}{~~~~6.1. A refinement of Theorem \reff{MAIN-CR}.\hfill \thepage\par}
%----------------------------------------------------------
%@@@@@@@@@@@@@@@@@@@@@@@@@@@@@@@@@@@@@@@@@@@@@@@@@@@@@@@@@@
%@@@@@@@@@@@@@@@@@@@@@@@@@@@@@@@@@@@@@@@@@@@@@@@@@@@@@@@@@@
%@@@@@@@@@@@@@@@@@@@@@@@@@@@@@@@@@@@@@@@@@@@@@@@@@@@@@@@@@@
%@@@@@@@@@@@@@@@@@@@@@@@@@@@@@@@@@@@@@@@@@@@@@@@@@@@@@@@@@@
%@@@@@@@@@@@@@@@@@@@@@@@@@@@@@@@@@@@@@@@@@@@@@@@@@@@@@@@@@@
%@@@@@@@@@@@@@@@@@@@@@@@@@@@@@@@@@@@@@@@@@@@@@@@@@@@@@@@@@@
%----------------------------------------------------------
\begin{theorem}\lbl{REF-MAIN-CR} Let $n<p<\infty$ and let $\mu$ be a non-trivial non-negative Borel measure on $\RN$.  There exist constants $\gamma=\gamma(n)>0$ and  $N=N(n)\in\N$, a family $\Qc$ consisting of pairwise disjoint cubes and a family $\tQc$ of cubes in $\RN$ with covering multiplicity  $M(\tQc)\le N$, mappings
%----------------------------------------------------------
\bel{M-12}
\Qc\ni Q\mapsto Q'\in\tQc~~~\text{and}~~~\Qc\ni Q\mapsto Q''\in\tQc
\ee
%----------------------------------------------------------
satisfying the condition
%----------------------------------------------------------
\bel{G-16}
Q'\cup Q''\subset \gamma Q~~~\text{for all}~~~Q\in\Qc,
\ee
%----------------------------------------------------------
such that for every function $f\in L_{p,loc}(\RN;\mu)$ the following equivalence
%----------------------------------------------------------
\bel{REF-NRM}
\|f\|_{\sum}\sim\left(\shuge_{Q\in\Qc}\,\,
\frac{(\diam Q)^{n-p}\iint \limits_{Q'\times Q''}
|f(x)-f(y)|^p\, d\mu(x)d\mu(y)}
{ \{(\diam Q')^{n-p}+\mu(Q')\} \{(\diam Q'')^{n-p}+\mu(Q'')\}}\right)^{\frac1p}
\ee
%----------------------------------------------------------
holds. The constants of this equivalence depend only on $n$ and $p$.
%----------------------------------------------------------
\end{theorem}
%----------------------------------------------------------
%@@@@@@@@@@@@@@@@@@@@@@@@@@@@@@@@@@@@@@@@@@@@@@@@@@@@@@@@@@
%@@@@@@@@@@@@@@@@@@@@@@@@@@@@@@@@@@@@@@@@@@@@@@@@@@@@@@@@@@
%@@@@@@@@@@@@@@@@@@@@@@@@@@@@@@@@@@@@@@@@@@@@@@@@@@@@@@@@@@
%@@@@@@@@@@@@@@@@@@@@@@@@@@@@@@@@@@@@@@@@@@@@@@@@@@@@@@@@@@
%@@@@@@@@@@@@@@@@@@@@@@@@@@@@@@@@@@@@@@@@@@@@@@@@@@@@@@@@@@
%----------------------------------------------------------
\par {\it Proof.} We let $I(f;\Qc)$ denote the quantity from the right-hand side of the equivalence \rf{REF-NRM}.
Then inequality
%----------------------------------------------------------
$$
I(f;\Qc)\le C(n,p)\|f\|_{\sum}
$$
%----------------------------------------------------------
follows from Proposition \reff{OP-W} where one can put $\Sc=\Qc$. In fact, it can be easily seen that $I(f;\Qc)$
does not exceed the quantity in the left-hand side of inequality \rf{N-11}.
%----------------------------------------------------------
\par Prove that for certain families $\Qc$ and $\tQc$ and mappings from \rf{M-12} with condition \rf{G-16} each depending only on $p,n,$ and the measure $\mu$, we have
%----------------------------------------------------------
$$
\|f\|_{\sum}\le C(n,p)I(f;\Qc)
$$
%----------------------------------------------------------
provided $f\in L_{p,loc}(\RN;\mu)$ is an arbitrary function.
%----------------------------------------------------------
\par We construct these objects using the method of proof of the sufficiency part of Theorem \reff{MAIN-CR}.
%----------------------------------------------------------
\par Let $E$ be the set constructing in Proposition \reff{RPR} and let $\Kc_E$ be the family of cubes defined by \rf{KE-DF}.
%----------------------------------------------------------
\par Let $f_1$ be the function defined by formula \rf{DEF-F1} and let $f_2=f-f_1$. As before given a cube $K\in W_E$ we put
%----------------------------------------------------------
$$
V_K:=\{Q\in W_E:Q\cap K\ne\emp\}.
$$
%----------------------------------------------------------
Then, by \rf{FE-K},
%----------------------------------------------------------
$$
\|f_1\|_{\LOP}^p\le C(n)\,\smed_{K\in\,W_E}
\,\,\smed_{Q\in V_K} \frac{|\tilde{f}_1(a_Q)-\tilde{f}_1(a_K)|^p}
{(\diam K)^{p-n}}.
$$
%----------------------------------------------------------
\par Recall that, by \rf{DEF-TF1},
%----------------------------------------------------------
\bel{DF2-TF}
\tf_1(x):=f_Q=\frac{1}{\mu(Q)}\intl_Q f\,d\mu~~~\text{for every}~~~x\in E,
\ee
%----------------------------------------------------------
where $Q=K^{(x)}=Q(x,R(x))$ is the unique cube from $\Kc_E$ with center at the point $x\in E$.
%----------------------------------------------------------
\par Also recall that, by Lemma \reff{SMC},
%----------------------------------------------------------
\bel{V-A1}
\|f_1\|_{\LOP}^p\le C(n)\,I_1
\ee
%----------------------------------------------------------
where
%----------------------------------------------------------
\bel{IO-1}
I_1:=\smed_{K\in\,\Ac}
\,\,\smed_{Q\in V_K} \frac{|\tilde{f}_1(a_Q)-\tilde{f}_1(a_K)|^p}
{(\diam K)^{p-n}}.
\ee
%----------------------------------------------------------
Here $\Ac\subset W_E$ is the family of cubes defined by \rf{A-D}.
%----------------------------------------------------------
\par Fix a cube $K\in\Ac$. Let $Q\in V_K$. Recall that, by \rf{AQ-T}, $a_Q\in \tau Q$. Since $Q\cap K\ne\emp$ and $\diam Q\le 4 \diam K$, we have $Q\subset 8K$. Hence $a_Q\in 8\tau K$.
%----------------------------------------------------------
\par Since $K\in \Ac$, we have $K\cap K^{(a_Q)}=\emp$. These properties of $a_Q$ and $K^{(a_Q)}$ easily imply that %----------------------------------------------------------
\bel{V-A2}
K^{(a_Q)}\subset \gamma_1 K~~~\text{for every}~~~ Q\in V_K
\ee
%----------------------------------------------------------
with $\gamma_1=156\tau$.
%----------------------------------------------------------
\par Let us divide $K$ into $2^n$ equal cubes of diameter $\tfrac12\diam K$. Let us fix one of these cubes and denote that cube by $\tK$. Let
%----------------------------------------------------------
\bel{TK-D}
T_K:=\tfrac12\tK.
\ee
%----------------------------------------------------------
\par Let $V_K\setminus\{K\}=\{Q_1,Q_2,...,Q_m\}$. We know that $m=m(K)\le C(n)$. Obviously there exists a family of {\it pairwise disjoint} equal cubes
%----------------------------------------------------------
$$
Y_K:=\{\tQ_1,\tQ_2,...,\tQ_m\}
$$
%----------------------------------------------------------
such that $\tQ_i\subset T_K $ and
%----------------------------------------------------------
$$
\diam K\le C(n)\diam \tQ_i,~~~\text{for every}~~~i=1,...,m.
$$
%----------------------------------------------------------
Here $C=C(n)$ is a constant depending only on $n$. Then, clearly, $K\subset \gamma_2 \tQ_i$ for every $1\le i\le m$ with certain $\gamma_2=\gamma_2(n)$ so that, by  \rf{V-A2},
%----------------------------------------------------------
\bel{KQ-3}
K^{(a_{Q_i})},K^{(a_{K})}\subset \gamma_3(n) \tQ_i~~~\text{for every}~~~ i=1,...,m.
\ee
%----------------------------------------------------------
\par Let $Q\in Y_K$; thus $Q=\tQ_i$ for some $1\le i\le m$. We assign to $Q$ cubes $Q',Q''\in\Kc_E$ as follows:
%----------------------------------------------------------
\bel{QK-1}
Q':=K^{(a_{Q_i})}~~~\text{and}~~~Q'':=K^{(a_{K})}.
\ee
%----------------------------------------------------------
Thus $Q',Q''\in\Kc_E$.
%----------------------------------------------------------
\par Then, by \rf{KQ-3}, $Q',Q''\subset  \gamma_3 Q$.
\par Furthermore, since $\diam Q\sim\diam K$ for each $Q\in Y_K$, by \rf{DF2-TF},
%----------------------------------------------------------
$$
\smed_{Q\in V_K} \frac{|\tilde{f}_1(a_Q)-\tilde{f}_1(a_K)|^p}
{(\diam K)^{p-n}}\sim \smed_{Q\in \Qc_K} \frac{|f_{Q'}-f_{Q''}|^p}
{(\diam Q)^{p-n}}.
$$
%----------------------------------------------------------
\par Let $\Qc_1=\cup\{Y_K: K\in\Ac\}$. We obtain
%----------------------------------------------------------
$$
I_1\sim\smed_{Q\in\Qc_1}
\frac{|f_{Q'}-f_{Q''}|^p}
{(\diam Q)^{p-n}}.
$$
%----------------------------------------------------------
\par But
%----------------------------------------------------------
$$
|\av{Q'}-\av{Q''}|^p\le
\frac{1}{\mu(Q')\mu(Q'')}\iint \limits_{Q'\times Q''}
|f(x)-f(y)|^p\, d\mu(x)d\mu(y)
$$
%----------------------------------------------------------
so that
%----------------------------------------------------------
$$
I_1\le C(n)\smed_{Q\in\Qc_1}\,\,
\frac{(\diam Q)^{p-n}}{\mu(Q')\mu(Q'')}\iint \limits_{Q'\times Q''}
|f(x)-f(y)|^p\, d\mu(x)d\mu(y).
$$
%----------------------------------------------------------
Since $Q',Q''\in\Kc_E$, by \rf{M-DK}, $\mu(Q')\sim (\diam Q')^{n-p}$ and $\mu(Q'')\sim (\diam Q'')^{n-p}$ so that
%----------------------------------------------------------
\bel{L-I1}
I_1\le C(n)\,\shuge_{Q\in\Qc_1}\,\,
\frac{(\diam Q)^{n-p}\iint \limits_{Q'\times Q''}
|f(x)-f(y)|^p\, d\mu(x)d\mu(y)}
{ \{(\diam Q')^{n-p}+\mu(Q')\} \{(\diam Q'')^{n-p}+\mu(Q'')\}}.
\ee
%----------------------------------------------------------
Combining this inequality with \rf{V-A1} we obtain
%----------------------------------------------------------
\bel{F1-R}
\|f_1\|_{\LOP}^p\le C(n)\,\shuge_{Q\in\Qc_1}\,\,
\frac{(\diam Q)^{n-p}\iint \limits_{Q'\times Q''}
|f(x)-f(y)|^p\, d\mu(x)d\mu(y)}
{ \{(\diam Q')^{n-p}+\mu(Q')\} \{(\diam Q'')^{n-p}+\mu(Q'')\}}.
\ee
%----------------------------------------------------------
\bigskip
%----------------------------------------------------------
\par Let us estimate $\|f_2\|_{\LPM}$ using the scheme of the proof of the sufficiency part of Theorem \reff{MAIN-CR} given in Section 5. We will also use the settings of this section.
%----------------------------------------------------------
\par By \rf{P-A4},
%----------------------------------------------------------
$$
\smed_{K\in\Ac}\,\intl_K |f_1-\tf_1(a_K)|^p\,d\mu\le C(n,p)\smed_{K\in\Ac}\smed_{Q\in V_K} \,\frac{|\tilde{f}_1(a_{Q})-\tf_1(a_K)|^p}{(\diam K)^{p-n}}
$$
%----------------------------------------------------------
so that, by \rf{IO-1},
%----------------------------------------------------------
\bel{S-I1}
\smed_{K\in\Ac}\,\intl_K |f_1-\tf_1(a_K)|^p\,d\mu\le C(n,p)\,I_1.
\ee
%----------------------------------------------------------
In turn, by \rf{F-S5},
%----------------------------------------------------------
$$
\intl_{\RN}|f_2|^pd\mu\le \smed_{K\in\Kc_E}\, \intl_{S_K}|f_2|^pd\mu+\smed_{K\in\Ac}\, \intl_{K}|f_2|^pd\mu.
$$
%----------------------------------------------------------
Combining this inequality with \rf{F2-R}, Lemma \reff{F2-S} and \rf{S-I1}, we obtain
%----------------------------------------------------------
\bel{F2-I12}
\intl_{\RN}|f_2|^pd\mu\le C\{I_1+I_2\}.
\ee
%----------------------------------------------------------
Recall that $I_2$ is defined by \rf{I2-D}.
%----------------------------------------------------------
\par Let $\Qc$ be the family of cubes defined by \rf{Q-U}.
Let us slightly modify this family as follows.  Recall for each cube $K\in W_E$ we have introduced a cube $\tK$ as one of the cubes from partition of $K$ into the family of $2^n$ equal cubes. See \rf{TK-D}.
%----------------------------------------------------------
\par Let us introduce another cube from this partition and denote this cube by $\hat{K}$. Thus $\diam \hat{K}=\tfrac12\diam K$ and $\tK\ne\hat{K}$ for every cube $K\in W_E$.
%----------------------------------------------------------
\par Let
%----------------------------------------------------------
$$
\hat{\Ac}=\{\tfrac12\hat{Q}: Q\in\Ac\}
$$
%----------------------------------------------------------
and let
%----------------------------------------------------------
$$
\Qc_2:=\hat{\Ac}\cup\Kc_E.
$$
%----------------------------------------------------------
In other words we replace in definition \rf{Q-U} the family $\Ac$ with the family $\hat{\Ac}$. Since $\diam \hat{Q}\sim \diam Q$ and $Q\subset 3\hat{Q}$ the result of Proposition \reff{W-LI} remains true after such a modification, i.e.,
%----------------------------------------------------------
$$
I_2\le
\smed_{Q\in\Qc_2}\left(\frac{\diam Q' \diam Q''}{\diam Q}\right)^{p-n} \iint \limits_{Q'\times Q''}
|f(x)-f(y)|^p\, d\mu(x)d\mu(y).
$$
%----------------------------------------------------------
Note that for each cube $Q\in\Qc_2$ we have $\mu(Q')\le C(\diam Q')^{n-p}$ and the same is true for $Q''$. Hence
%----------------------------------------------------------
\bel{I2-Y}
I_2\le C\,
\shuge_{Q\in\Qc_2}\,\,
\frac{(\diam Q)^{n-p}\iint \limits_{Q'\times Q''}
|f(x)-f(y)|^p\, d\mu(x)d\mu(y)}
{ \{(\diam Q')^{n-p}+\mu(Q')\} \{(\diam Q'')^{n-p}+\mu(Q'')\}}.
\ee
%----------------------------------------------------------
\par Also we remark that by definition $\Qc_1\cap\Qc_2=\emp$, i.e., the family
%----------------------------------------------------------
$$
\Qc:=\Qc_1\cup\Qc_2
$$
%----------------------------------------------------------
consists of {\it pairwise disjoint} cubes. Furthermore, given $Q\in\Qc$ the cubes $Q'$ and $Q''$ belong to the family
%----------------------------------------------------------
$$
\tQc:=\Ac\cup\Kc_E
$$
%----------------------------------------------------------
See \rf{QK-1}, \rf{RM-2} and \rf{R1-QP}. Since the cubes of the family $\Kc_E$ are pairwise disjoint and $\Ac\subset W_E$, covering multiplicity of the family $\tQc$ satisfies the following inequality
%----------------------------------------------------------
$$
M(\tQc)=M(\Ac\cup\Kc_E)\le M(\Ac)+M(\Kc_E)\le M(W_E)+1\le N(n).
$$
%---------------------------------------------------------- %@@@@@@@@@@@@@@@@@@@@@@@@@@@@@@@@@@@@@@@@@@@@@@@@@@@@@@@@@@
%----------------------------------------------------------
\par Now combining inequality \rf{I2-Y} with inequalities \rf{F2-I12} and \rf{L-I1}, we obtain
%----------------------------------------------------------
$$
\|f_2\|_{\LPM}^p\le C\,
\shuge_{Q\in\Qc}\,\,
\frac{(\diam Q)^{n-p}\iint \limits_{Q'\times Q''}
|f(x)-f(y)|^p\, d\mu(x)d\mu(y)}
{ \{(\diam Q')^{n-p}+\mu(Q')\} \{(\diam Q'')^{n-p}+\mu(Q'')\}}.
$$
%----------------------------------------------------------
\par Finally, this inequality and \rf{F1-R} imply that
%----------------------------------------------------------
\be
\|f\|_{\Sigma}^p&\le& 2^p\{\|f_1\|_{\LOP}^p+\|f_2\|_{\LPM}^p\}\nn\\&\le&
C\,\shuge_{Q\in\Qc}\,\,
\frac{(\diam Q)^{n-p}\iint \limits_{Q'\times Q''}
|f(x)-f(y)|^p\, d\mu(x)d\mu(y)}
{ \{(\diam Q')^{n-p}+\mu(Q')\} \{(\diam Q'')^{n-p}+\mu(Q'')\}}\nn
\ee
%----------------------------------------------------------
proving the theorem.\bx
%----------------------------------------------------------
%@@@@@@@@@@@@@@@@@@@@@@@@@@@@@@@@@@@@@@@@@@@@@@@@@@@@@@@@@@
%@@@@@@@@@@@@@@@@@@@@@@@@@@@@@@@@@@@@@@@@@@@@@@@@@@@@@@@@@@
%@@@@@@@@@@@@@@@@@@@@@@@@@@@@@@@@@@@@@@@@@@@@@@@@@@@@@@@@@@
%@@@@@@@@@@@@@@@@@@@@@@@@@@@@@@@@@@@@@@@@@@@@@@@@@@@@@@@@@@
%@@@@@@@@@@@@@@@@@@@@@@@@@@@@@@@@@@@@@@@@@@@@@@@@@@@@@@@@@@
%@@@@@@@@@@@@@@@@@@@@@@@@@@@@@@@@@@@@@@@@@@@@@@@@@@@@@@@@@@
%----------------------------------------------------------
\medskip
%----------------------------------------------------------
\par {\bf 6.2. Lacunae of Whitney's cubes.}
%----------------------------------------------------------
\addtocontents{toc}{~~~~6.2. Lacunae of Whitney's cubes. \hfill \thepage\par} In the next subsection we present another refinement of Theorem \reff{MAIN-CR}. We obtain this refinement with the help of a modification of the classical Whitney extension method which we described and used at the beginning of Section 4. See formula \rf{DEF-F1}.
%----------------------------------------------------------
\par As we have noted in Section 1 the main idea of this approach is to use certain {\it families} of Whitney's cubes rather than to treat each Whitney cube separately. We call these families of Whitney cubes {\it lacunae}.
%----------------------------------------------------------
\par In this subsection we present main definitions and main properties of lacunae. For the proof of these properties we refer the reader to paper \cite{S6}, Sections 4-5.
%----------------------------------------------------------
\par Let $E$ be a closed subset of $\RN$ and let $W_E$ be a Whitney decomposition of its complement $\RN\setminus E$, see Theorem \reff{Wcov} and Lemma \reff{Wadd}. As we have already noted in Section 4, see \rf{9QE},
%----------------------------------------------------------
\bel{INT-W}
(9Q)\cap E\ne\emp~~~ \text{for every}~~~ Q\in W_E.
\ee
%----------------------------------------------------------
%----------------------------------------------------------
\par By $LW_E$ we denote a subfamily of Whitney cubes satisfying the following condition:
%----------------------------------------------------------
\bel{L-PR}
(10Q)\cap E=(\q Q)\cap E.
\ee
%----------------------------------------------------------
\par Then we introduce a binary relation $\sim$ on $LW_E$: for every $Q_1,Q_2\in LW_E$
%----------------------------------------------------------
$$
Q_1\sim Q_2 ~~\Longleftrightarrow~ (10Q_1)\cap E= (10Q_2)\cap E.
$$
%----------------------------------------------------------
\par It can be easily seen that $\sim$ satisfies the axioms of equivalence relations, i.e., it is reflexive, symmetric and transitive. Given a cube $Q\in LW_E$ by
%----------------------------------------------------------
$$
[Q]:=\{K\in LW_E: K\sim Q\}
$$
%----------------------------------------------------------
we denote the equivalence class of $Q$. We refer to this equivalence class as {\it a true lacuna} with respect to the set $E$.
%----------------------------------------------------------
\par Let
%----------------------------------------------------------
$$
\tL_E=LW_E\backslash\sim\,=\{[Q]: Q\in LW_E\}
$$
%----------------------------------------------------------
be the corresponding quotient set of $LW_E$ by $\sim$\,, i.e., the set of all possible equivalence classes (lacunae) of $LW_E$ by $\sim$\,.
%----------------------------------------------------------
\par Thus for every pair of Whitney cubes $Q_1,Q_2\in W_E$ which belong to a true lacuna $L\in\tL_E$ we have
%----------------------------------------------------------
\bel{I-L}
(10Q_1)\cap E=(\q Q_1)\cap E=(10Q_2)\cap E=(\q Q_2)\cap E.
\ee
%----------------------------------------------------------
By $V_L$ we denote the associated set of the lacuna $L$
%----------------------------------------------------------
\bel{D-VL}
V_L:=(\q Q)\cap E.
\ee
%----------------------------------------------------------
Here $Q$ is an arbitrary cube from $L$. By \rf{I-L}, any choice of a cube $Q\in L$ provides the same set $V_L$ so that $V_L$ is well-defined. Also note that for each cube $Q$ which belong to a true lacuna $L$ we have $V_L=(10Q)\cap E.$\medskip
%----------------------------------------------------------
\par We extend the family $\tL_E$ of true lacunae to a family $\LE$ of {\it all lacunae} in the following way. Suppose that $Q\in W_E\setminus LW_E$, see \rf{L-PR}, i.e.,
%----------------------------------------------------------
\bel{A-L}
(10Q)\cap E\ne(\q Q)\cap E.
\ee
%----------------------------------------------------------
In this case to the cube $Q$ we assign a lacuna $L:=\{Q\}$ consisting of a unique cube - the cube $Q$.  We also put $V_L:=(\q Q)\cap E$ as in \rf{D-VL}.
%----------------------------------------------------------
\par We refer to such a lacuna $L:=\{Q\}$ as an {\it elementary lacuna} with respect to the set $E$. By $\hL_E$ we denote the family of all elementary lacunae with respect to $E$:
%----------------------------------------------------------
$$
\hL_E:=\{L=\{Q\}:Q\in W_E\setminus LW_E\}
$$
%----------------------------------------------------------
\par We note that property \rf{A-L} implies the existence of a point
%----------------------------------------------------------
$$
a\in (E\setminus (10Q))\cap (\q Q).
$$
%----------------------------------------------------------
On the other hand, by \rf{INT-W}, there exists a point $$b\in(9Q)\cap E.$$ Hence
%----------------------------------------------------------
$$
\|a-b\|\ge r_Q=(1/2)\diam Q
$$
%----------------------------------------------------------
so that
%----------------------------------------------------------
$$
\diam V_L=\diam ((\q Q)\cap E)\ge \tfrac12\diam Q
$$
%----------------------------------------------------------
provided
%----------------------------------------------------------
$$
L=\{Q\}\in \hL_E
$$
%----------------------------------------------------------
is an elementary lacuna.
%----------------------------------------------------------
\par Finally, by $\LE$ we denote the family of all lacunae with respect to $E$:
%----------------------------------------------------------
$$
\LE=\tL_E\cup \hL_E.
$$
%----------------------------------------------------------
\par We turn to description of main properties of lacunae. Recall that the detailed proofs of these properties one can find in \cite{S6}, Sections 4-5.
%----------------------------------------------------------
%@@@@@@@@@@@@@@@@@@@@@@@@@@@@@@@@@@@@@@@@@@@@@@@@@@@@@@@@@@
%@@@@@@@@@@@@@@@@@@@@@@@@@@@@@@@@@@@@@@@@@@@@@@@@@@@@@@@@@@
%@@@@@@@@@@@@@@@@@@@@@@@@@@@@@@@@@@@@@@@@@@@@@@@@@@@@@@@@@@
%----------------------------------------------------------
\begin{proposition} Let $L\in\LE$ be a lacuna.
If $\diam V_L>0$, then there exists a cube $Q_L\in L$ such that
%----------------------------------------------------------
$$
\diam Q_L=\min\{\diam Q: Q\in L\}.
$$
%----------------------------------------------------------
Furthermore,
%----------------------------------------------------------
$$
\tfrac{1}{\q}\diam V_L\le \diam Q_L\le \gamma \diam V_L
$$
%----------------------------------------------------------
where $\gamma$ is an absolute constant.
%@@@@@@@@@@@@@@@@@@@@@@@@@@@@@@@@@@@@@@@@@@@@@@@@@@@@@@@@@@
%----------------------------------------------------------
\end{proposition}
%----------------------------------------------------------
%----------------------------------------------------------
%@@@@@@@@@@@@@@@@@@@@@@@@@@@@@@@@@@@@@@@@@@@@@@@@@@@@@@@@@@
%@@@@@@@@@@@@@@@@@@@@@@@@@@@@@@@@@@@@@@@@@@@@@@@@@@@@@@@@@@
%@@@@@@@@@@@@@@@@@@@@@@@@@@@@@@@@@@@@@@@@@@@@@@@@@@@@@@@@@@
%@@@@@@@@@@@@@@@@@@@@@@@@@@@@@@@@@@@@@@@@@@@@@@@@@@@@@@@@@@
%----------------------------------------------------------
\par Given a lacuna  $L\in\LE$ we let $U_L$ denote the union of all cubes which belong to the lacuna:
%----------------------------------------------------------
$$
U_L:=\cup\{Q:Q\in L\}.
$$
%----------------------------------------------------------
By $\diam L$ we denote the diameter of the set $U_L$:
%----------------------------------------------------------
$$
\diam L:=\diam U_L=\sup\{\|a-b\|:a,b\in U_L\}.
$$
%----------------------------------------------------------
\par We say that $L$ is bounded if $\diam L<\infty$. If $\diam L=\infty$ we say that $L$ is an unbounded lacuna.
%----------------------------------------------------------
\begin{proposition} (i). For every lacuna $L\in\LE$
%----------------------------------------------------------
$$
\diam L\sim \sup\{\diam Q: Q\in L\}\sim \dist(V_L,E\setminus V_L)
$$
%----------------------------------------------------------
with absolute constants in the equivalences;
%----------------------------------------------------------
%@@@@@@@@@@@@@@@@@@@@@@@@@@@@@@@@@@@@@@@@@@@@@@@@@@@@@@@@@@
\par (ii). If $E$ is an unbounded set then  every lacuna $L\in\LE$ is bounded;
%----------------------------------------------------------
\par (iii). If $E$ is bounded, there exists the unique unbounded lacuna $L^{\max}\in\LE$. The lacuna $L^{\max}$ is a true lacuna for which $V_{L^{\max}}=E$.
%----------------------------------------------------------
%@@@@@@@@@@@@@@@@@@@@@@@@@@@@@@@@@@@@@@@@@@@@@@@@@@@@@@@@@@
%----------------------------------------------------------
\end{proposition}
%@@@@@@@@@@@@@@@@@@@@@@@@@@@@@@@@@@@@@@@@@@@@@@@@@@@@@@@@@@
%@@@@@@@@@@@@@@@@@@@@@@@@@@@@@@@@@@@@@@@@@@@@@@@@@@@@@@@@@@
%@@@@@@@@@@@@@@@@@@@@@@@@@@@@@@@@@@@@@@@@@@@@@@@@@@@@@@@@@@
%@@@@@@@@@@@@@@@@@@@@@@@@@@@@@@@@@@@@@@@@@@@@@@@@@@@@@@@@@@
%----------------------------------------------------------
%@@@@@@@@@@@@@@@@@@@@@@@@@@@@@@@@@@@@@@@@@@@@@@@@@@@@@@@@@@
%@@@@@@@@@@@@@@@@@@@@@@@@@@@@@@@@@@@@@@@@@@@@@@@@@@@@@@@@@@
%@@@@@@@@@@@@@@@@@@@@@@@@@@@@@@@@@@@@@@@@@@@@@@@@@@@@@@@@@@
%@@@@@@@@@@@@@@@@@@@@@@@@@@@@@@@@@@@@@@@@@@@@@@@@@@@@@@@@@@
%----------------------------------------------------------
%@@@@@@@@@@@@@@@@@@@@@@@@@@@@@@@@@@@@@@@@@@@@@@@@@@@@@@@@@@
%----------------------------------------------------------
\begin{proposition} Let $L\in \LE$ be a bounded lacuna. Then there exists a cube  $\QL\in L$ such that
%----------------------------------------------------------
$$
\diam \QL=\max\{\diam K: K\in L\}.
$$
%----------------------------------------------------------
Furthermore,
%----------------------------------------------------------
\bel{CL-D}
\diam \QL\sim \diam L\sim\dist(V_L,E\setminus V_L),
\ee
%----------------------------------------------------------
and
%----------------------------------------------------------
$$
V_L\cup U_L\subset \gamma \QL.
$$
%----------------------------------------------------------
Here the constant $\gamma$ and constants in the equivalences of \rf{CL-D} are absolute .
%@@@@@@@@@@@@@@@@@@@@@@@@@@@@@@@@@@@@@@@@@@@@@@@@@@@@@@@@@@
%----------------------------------------------------------
\end{proposition}
%----------------------------------------------------------
%@@@@@@@@@@@@@@@@@@@@@@@@@@@@@@@@@@@@@@@@@@@@@@@@@@@@@@@@@@
%@@@@@@@@@@@@@@@@@@@@@@@@@@@@@@@@@@@@@@@@@@@@@@@@@@@@@@@@@@
%@@@@@@@@@@@@@@@@@@@@@@@@@@@@@@@@@@@@@@@@@@@@@@@@@@@@@@@@@@
%@@@@@@@@@@@@@@@@@@@@@@@@@@@@@@@@@@@@@@@@@@@@@@@@@@@@@@@@@@
%@@@@@@@@@@@@@@@@@@@@@@@@@@@@@@@@@@@@@@@@@@@@@@@@@@@@@@@@@@
%@@@@@@@@@@@@@@@@@@@@@@@@@@@@@@@@@@@@@@@@@@@@@@@@@@@@@@@@@@
%@@@@@@@@@@@@@@@@@@@@@@@@@@@@@@@@@@@@@@@@@@@@@@@@@@@@@@@@@@
%----------------------------------------------------------
\begin{proposition} Let $L\in\LE$ be a lacuna and let $Q\in L$. Suppose that there exist a lacuna $L'\in \LE$, $L\ne L'$, and a cube $Q'\in L'$ such that $Q\cap Q'\ne\emp$. Then:
%----------------------------------------------------------
\par (i). If $L$ is a true lacuna, then $L'$ is an elementary lacuna, i.e., $L'\in \hL_E=\LE\setminus\tL_E$;
%----------------------------------------------------------
\par (ii). Either
%----------------------------------------------------------
$$
\diam Q\sim \diam V_L\sim \diam Q_L
$$
%----------------------------------------------------------
or
%----------------------------------------------------------
$$
\diam Q\sim \dist(V_L,E\setminus V_L)\sim \diam\QL
$$
%----------------------------------------------------------
with absolute constants in the equivalences.
%----------------------------------------------------------
\end{proposition}
%----------------------------------------------------------
%@@@@@@@@@@@@@@@@@@@@@@@@@@@@@@@@@@@@@@@@@@@@@@@@@@@@@@@@@@
%@@@@@@@@@@@@@@@@@@@@@@@@@@@@@@@@@@@@@@@@@@@@@@@@@@@@@@@@@@
%@@@@@@@@@@@@@@@@@@@@@@@@@@@@@@@@@@@@@@@@@@@@@@@@@@@@@@@@@@
%@@@@@@@@@@@@@@@@@@@@@@@@@@@@@@@@@@@@@@@@@@@@@@@@@@@@@@@@@@
%----------------------------------------------------------
\begin{proposition}\lbl{INT-L1} Let $L\in\LE$ be a lacuna and let
%----------------------------------------------------------
$$
\Ic_L:=\{K\in W_E\setminus L: \exists~Q\in L~~\text{such that}~~K\cap Q\ne\emp\}.
$$
%----------------------------------------------------------
Then $\card\Ic_L\le\gamma(n).$
%----------------------------------------------------------
\end{proposition}
%----------------------------------------------------------
%@@@@@@@@@@@@@@@@@@@@@@@@@@@@@@@@@@@@@@@@@@@@@@@@@@@@@@@@@@
%@@@@@@@@@@@@@@@@@@@@@@@@@@@@@@@@@@@@@@@@@@@@@@@@@@@@@@@@@@
%@@@@@@@@@@@@@@@@@@@@@@@@@@@@@@@@@@@@@@@@@@@@@@@@@@@@@@@@@@
%@@@@@@@@@@@@@@@@@@@@@@@@@@@@@@@@@@@@@@@@@@@@@@@@@@@@@@@@@@
%----------------------------------------------------------
\par One of the main ingredient of the lacunary approach is a mapping $\LE\ni L\mapsto \PRL(L)\in E$ whose properties are described by the following
%----------------------------------------------------------
%@@@@@@@@@@@@@@@@@@@@@@@@@@@@@@@@@@@@@@@@@@@@@@@@@@@@@@@@@@
%@@@@@@@@@@@@@@@@@@@@@@@@@@@@@@@@@@@@@@@@@@@@@@@@@@@@@@@@@@
%@@@@@@@@@@@@@@@@@@@@@@@@@@@@@@@@@@@@@@@@@@@@@@@@@@@@@@@@@@
%@@@@@@@@@@@@@@@@@@@@@@@@@@@@@@@@@@@@@@@@@@@@@@@@@@@@@@@@@@
%@@@@@@@@@@@@@@@@@@@@@@@@@@@@@@@@@@@@@@@@@@@@@@@@@@@@@@@@@@
%@@@@@@@@@@@@@@@@@@@@@@@@@@@@@@@@@@@@@@@@@@@@@@@@@@@@@@@@@@
%----------------------------------------------------------
\begin{proposition}\lbl{L-PE} There exist an absolute constant $\gamma>0$ and a mapping
%----------------------------------------------------------
$$
\LE\ni L~~\longrightarrow~~\PRL(L)\in E
$$
%----------------------------------------------------------
such that:
%----------------------------------------------------------
\par (i). For every lacuna $L\in\LE$ we have
%----------------------------------------------------------
\bel{PR-DSE}
\PRL(L)\in (\gamma\,Q_L)\cap E~;
\ee
%----------------------------------------------------------
\par (ii). For every $a\in E$
%----------------------------------------------------------
$$
\card\{L\in\LE:\PRL(L)=a\}\le C(n).
$$
%----------------------------------------------------------
\end{proposition}
%@@@@@@@@@@@@@@@@@@@@@@@@@@@@@@@@@@@@@@@@@@@@@@@@@@@@@@@@@@
%@@@@@@@@@@@@@@@@@@@@@@@@@@@@@@@@@@@@@@@@@@@@@@@@@@@@@@@@@@
%----------------------------------------------------------
\par We refer to the mapping $\PRL$ as a ``projection'' of $\Lc_E$ into the set $E$.\smallskip
%----------------------------------------------------------
\par Let $L\in\LE$ be a lacuna. Recall that
%----------------------------------------------------------
$$
U_L=\cup\{Q: Q\in L\}.
$$
%----------------------------------------------------------
%@@@@@@@@@@@@@@@@@@@@@@@@@@@@@@@@@@@@@@@@@@@@@@@@@@@@@@@@@@
%@@@@@@@@@@@@@@@@@@@@@@@@@@@@@@@@@@@@@@@@@@@@@@@@@@@@@@@@@@
%@@@@@@@@@@@@@@@@@@@@@@@@@@@@@@@@@@@@@@@@@@@@@@@@@@@@@@@@@@
%@@@@@@@@@@@@@@@@@@@@@@@@@@@@@@@@@@@@@@@@@@@@@@@@@@@@@@@@@@
%----------------------------------------------------------
\begin{definition} {\em Let $L,L'\in\LE$ be lacunae.  We say that $L$ and $L'$ are {\it contacting
lacunae} if $U_L\cap U_{L'}\ne\emp.$ In this case we write $L\leftrightarrow L'$.}
%----------------------------------------------------------
\end{definition}
%----------------------------------------------------------
%@@@@@@@@@@@@@@@@@@@@@@@@@@@@@@@@@@@@@@@@@@@@@@@@@@@@@@@@@@
%@@@@@@@@@@@@@@@@@@@@@@@@@@@@@@@@@@@@@@@@@@@@@@@@@@@@@@@@@@
%----------------------------------------------------------
\par Thus $L\leftrightarrow L'$ whenever there exist cubes $Q\in L$ and $Q'\in L'$ such that $Q\cap Q'\ne\emp$. We refer to the pair of such cubes as {\it contacting cubes}. Let us present several properties of contacting lacunae and contacting cubes.
%----------------------------------------------------------
%@@@@@@@@@@@@@@@@@@@@@@@@@@@@@@@@@@@@@@@@@@@@@@@@@@@@@@@@@@
%@@@@@@@@@@@@@@@@@@@@@@@@@@@@@@@@@@@@@@@@@@@@@@@@@@@@@@@@@@
%@@@@@@@@@@@@@@@@@@@@@@@@@@@@@@@@@@@@@@@@@@@@@@@@@@@@@@@@@@
%@@@@@@@@@@@@@@@@@@@@@@@@@@@@@@@@@@@@@@@@@@@@@@@@@@@@@@@@@@
%----------------------------------------------------------
\begin{proposition} (i). Every lacuna $L\in\LE$ contacts with at most $C(n)$ lacunae, i.e.,
%----------------------------------------------------------
$$
\card\{L'\in\LE: L'\lr L\}\le C(n);
$$
%----------------------------------------------------------
\par (ii). Every true lacuna contacts only with elementary lacunae.
%----------------------------------------------------------
\end{proposition}
%----------------------------------------------------------
%@@@@@@@@@@@@@@@@@@@@@@@@@@@@@@@@@@@@@@@@@@@@@@@@@@@@@@@@@@
%@@@@@@@@@@@@@@@@@@@@@@@@@@@@@@@@@@@@@@@@@@@@@@@@@@@@@@@@@@
%----------------------------------------------------------
%@@@@@@@@@@@@@@@@@@@@@@@@@@@@@@@@@@@@@@@@@@@@@@@@@@@@@@@@@@
%@@@@@@@@@@@@@@@@@@@@@@@@@@@@@@@@@@@@@@@@@@@@@@@@@@@@@@@@@@
%@@@@@@@@@@@@@@@@@@@@@@@@@@@@@@@@@@@@@@@@@@@@@@@@@@@@@@@@@@
%@@@@@@@@@@@@@@@@@@@@@@@@@@@@@@@@@@@@@@@@@@@@@@@@@@@@@@@@@@
%----------------------------------------------------------
\begin{proposition} Let $L\in\LE$ be a lacuna and let $Q\in L$ be a contacting cube. (I.e., there exist a lacuna $L'\in\LE$ and a cube $Q'\in L'$ such that $Q\cap Q'\ne\emp$.) Then either
%----------------------------------------------------------
$$
\diam Q\sim \diam V_L\sim
\min\{\diam K: K\in L\}=\diam Q_L
$$
%----------------------------------------------------------
or
%----------------------------------------------------------
$$
\diam Q\sim \dist(V_L,E\setminus V_L)\sim
\max\{\diam K: K\in L\}=\diam \QL
$$
%----------------------------------------------------------
with absolute constants in the equivalences.
%@@@@@@@@@@@@@@@@@@@@@@@@@@@@@@@@@@@@@@@@@@@@@@@@@@@@@@@@@@
%----------------------------------------------------------
\end{proposition}
%----------------------------------------------------------
%@@@@@@@@@@@@@@@@@@@@@@@@@@@@@@@@@@@@@@@@@@@@@@@@@@@@@@@@@@
%@@@@@@@@@@@@@@@@@@@@@@@@@@@@@@@@@@@@@@@@@@@@@@@@@@@@@@@@@@
%%@@@@@@@@@@@@@@@@@@@@@@@@@@@@@@@@@@@@@@@@@@@@@@@@@@@@@@@@@
%@@@@@@@@@@@@@@@@@@@@@@@@@@@@@@@@@@@@@@@@@@@@@@@@@@@@@@@@@@
%@@@@@@@@@@@@@@@@@@@@@@@@@@@@@@@@@@@@@@@@@@@@@@@@@@@@@@@@@@
%%@@@@@@@@@@@@@@@@@@@@@@@@@@@@@@@@@@@@@@@@@@@@@@@@@@@@@@@@@
%@@@@@@@@@@@@@@@@@@@@@@@@@@@@@@@@@@@@@@@@@@@@@@@@@@@@@@@@@@
%@@@@@@@@@@@@@@@@@@@@@@@@@@@@@@@@@@@@@@@@@@@@@@@@@@@@@@@@@@
%@@@@@@@@@@@@@@@@@@@@@@@@@@@@@@@@@@@@@@@@@@@@@@@@@@@@@@@@@@
%----------------------------------------------------------
\medskip
%----------------------------------------------------------
\par {\bf 6.3. A lacunary refinement of the main result.}
%----------------------------------------------------------
\addtocontents{toc}{~~~~6.3. A lacunary refinement of the main result. \hfill \thepage\\\par}
%----------------------------------------------------------
\par The next refinement of Theorem \reff{MAIN-CR} is motivated by its important applications to characterization of the restrictions of Sobolev functions to closed subsets in $\RT$. See \cite{S6}.
%----------------------------------------------------------
%@@@@@@@@@@@@@@@@@@@@@@@@@@@@@@@@@@@@@@@@@@@@@@@@@@@@@@@@@@
%@@@@@@@@@@@@@@@@@@@@@@@@@@@@@@@@@@@@@@@@@@@@@@@@@@@@@@@@@@
%@@@@@@@@@@@@@@@@@@@@@@@@@@@@@@@@@@@@@@@@@@@@@@@@@@@@@@@@@@
%@@@@@@@@@@@@@@@@@@@@@@@@@@@@@@@@@@@@@@@@@@@@@@@@@@@@@@@@@@
%@@@@@@@@@@@@@@@@@@@@@@@@@@@@@@@@@@@@@@@@@@@@@@@@@@@@@@@@@@
%@@@@@@@@@@@@@@@@@@@@@@@@@@@@@@@@@@@@@@@@@@@@@@@@@@@@@@@@@@
%----------------------------------------------------------
\begin{theorem}\lbl{REF-2} Let $\mu$ be a non-trivial non-negative Borel measure on $\RN$, $n<p<\infty$, and let
%----------------------------------------------------------
$$
\ssmall=\SUM.
$$
%----------------------------------------------------------
\par There exist families of closed sets $\{G_1,G_2,...\}$ and $\{H_1,H_2,...\}$ with covering multiplicity $M(\{G_i\}),M(\{H_i\})\le C(n)$, and a family $\{\lambda_1,\lambda_2,...\}$ of positive numbers such that for every function $f\in L_{1,loc}(\RN;\mu)$ the following equivalence
%----------------------------------------------------------
$$
\|f\|_\Sigma^p\sim \sum_{i=1}^\infty \,\,
\lambda_i \iint\limits_{G_i\times H_i}
|f(x)-f(y)|^p\, d\mu(x)\,d\mu(y)
$$
%----------------------------------------------------------
holds. The constants of this equivalence depend only on $n$ and $p$.
%----------------------------------------------------------
\end{theorem}
%----------------------------------------------------------
%@@@@@@@@@@@@@@@@@@@@@@@@@@@@@@@@@@@@@@@@@@@@@@@@@@@@@@@@@@
%@@@@@@@@@@@@@@@@@@@@@@@@@@@@@@@@@@@@@@@@@@@@@@@@@@@@@@@@@@
%@@@@@@@@@@@@@@@@@@@@@@@@@@@@@@@@@@@@@@@@@@@@@@@@@@@@@@@@@@
%@@@@@@@@@@@@@@@@@@@@@@@@@@@@@@@@@@@@@@@@@@@@@@@@@@@@@@@@@@
%@@@@@@@@@@@@@@@@@@@@@@@@@@@@@@@@@@@@@@@@@@@@@@@@@@@@@@@@@@
%----------------------------------------------------------
\par {\it Proof.} We follow the scheme of the proof of Theorem \reff{REF-MAIN-CR}. Let $L\in\Lc_E$ be a lacuna. For the sake of brevity we put
%----------------------------------------------------------
$$
A_L:=\PRL(L)
$$
%----------------------------------------------------------
where $\PRL$ denotes the ``projection'' of $L$ into $E$ , see Proposition \reff{L-PE}.
%----------------------------------------------------------
\par We modify the Whitney extension formula as follows: We put
%----------------------------------------------------------
$$
a_Q=A_L~~~~\text{for all}~~~~Q\in L.
$$
%----------------------------------------------------------
Cf. \rf{AQ-T}.
%----------------------------------------------------------
\par Note that, by  \rf{PR-DSE},
%----------------------------------------------------------
$$
A_L\in(\gamma Q_L)\cap E
$$
%----------------------------------------------------------
where $Q_L$ is a cube in $L$ of minimal diameter. Since for every $Q\in L$
%----------------------------------------------------------
$$
(90 Q)\cap E=(90 Q_L)\cap E,
$$
%----------------------------------------------------------
we have $(90 Q)\cap(90 Q_L)\ne\emp$. But $\diam Q_L\le\diam Q$ so that
%----------------------------------------------------------
$$
Q_L\subset \gamma_1 Q~~~\text{for every}~~~Q\in L
$$
%----------------------------------------------------------
with some absolute $\gamma_1>0$. Hence
%----------------------------------------------------------
\bel{AL-2}
a_Q=A_L\subset  (\gamma_2 Q)\cap E~~~~\text{for every}~~~~Q\in L.
\ee
%----------------------------------------------------------
\par This shows that we can construct the component $f_1$ using the extension formula \rf{DEF-F1}. Then the functions $f_1$ and $f_2=f-f_1$ will provide an almost optimal decomposition of $f$, i.e.,
%----------------------------------------------------------
$$
\|f_1\|_{\LOP}+\|f_2\|_{\LPM}\sim \|f\|_{\Sigma}.
$$
%----------------------------------------------------------
\par Let us construct the required families $\{G_i\}$ and $\{H_i\}$ using the approach suggested in the proof of Theorem \reff{REF-MAIN-CR}.
%----------------------------------------------------------
\par We begin with the estimate of the quantity $\|f_1\|_{\LOP}$. First we modify the estimate \rf{V-A1} and definition \rf{IO-1}.
%----------------------------------------------------------
\par Let $L\in\Lc_E$. We note that for every $K\in L$ and every $Q\in V_K\cap L$ we have $a_Q=a_K=A_L$ so that
%----------------------------------------------------------
$$
|\tf_1(a_Q)-\tf_1(a_K)|=0.
$$
%----------------------------------------------------------
\par Let
%----------------------------------------------------------
$$
T(L):=\{K\in L:\exists\, Q\in W_E\setminus L, Q\cap K\ne\emp\}.
$$
%----------------------------------------------------------
Thus $T(L)$ is the family of contacting cubes of $L$. By Proposition \reff{INT-L1} and part (2) of Lemma \reff{Wadd},
%----------------------------------------------------------
\bel{TL-C}
\#T(L)\le C(n),~~~L\in\Lc_E.
\ee
%----------------------------------------------------------
Thee observations enable us to modify inequalities \rf{V-A1} and \rf{IO-1} as follows:
%----------------------------------------------------------
$$
\|f_1\|_{\LOP}^p\le C(n)\,\tI_1
$$
%----------------------------------------------------------
where
%----------------------------------------------------------
\bel{I-N}
\tI_1:=\smed_{L\in\Lc_E}
\,\,\smed_{K\in\, T(L)\cap \Ac}\,\,\,\smed_{Q\in V_K\setminus L} \frac{|\tilde{f}_1(a_Q)-\tilde{f}_1(a_K)|^p}
{(\diam K)^{p-n}}.
\ee
%----------------------------------------------------------
\par We construct the family $\Qc$ and the mappings
%----------------------------------------------------------
$$
\Qc\ni Q\mapsto Q'\in\Kc_E~~~\text{and}~~~
\Qc\ni Q\mapsto Q''\in\Kc_E
$$
%----------------------------------------------------------
precisely as in Theorem \reff{REF-MAIN-CR}, see formulas \rf{V-A2}-\rf{QK-1}, but only for cubes $K,Q$ from inequality \rf{I-N}, i.e., for contacting cubes.
%----------------------------------------------------------
\par As a result, we again obtain inequality \rf{F1-R}. Note that, by Proposition  \reff{OP-W}, we have the following:
%----------------------------------------------------------
\bel{OP-F}
\shuge_{Q\in\Qc_1}\,\,
\frac{(\diam Q)^{n-p}\iint \limits_{Q'\times Q''}
|f(x)-f(y)|^p\, d\mu(x)d\mu(y)}
{ \{(\diam Q')^{n-p}+\mu(Q')\} \{(\diam Q'')^{n-p}+\mu(Q'')\}}\le C(n)\,\|f\|_{\Sigma}^p.
\ee
%----------------------------------------------------------
\par The crucial point of this construction is as follows:
{\it the mappings
%----------------------------------------------------------
$$
\Qc\ni Q\mapsto Q'\in\Kc_E~~~\text{and}~~~
\Qc\ni Q\mapsto Q''\in\Kc_E
$$
%----------------------------------------------------------
are ``almost'' one-to-one.} Thus for every $K\in\Kc_E$ there exist at most $C(n)$ cubes $Q\in\Qc_1$ such that $Q'=K$. The same is true for the mapping $Q\mapsto Q''$.
%----------------------------------------------------------
\par This statement easily follows from a similar property of the ``projection'' operator $\PRL:\Lc_E\to E$ and inequality \rf{TL-C}. In fact, let $c_K$ be the center of $K$. Then, by  part (ii) Proposition \reff{L-PE}, there are at most $C_1(n)$ lacunae  $L\in\Lc_E$ such that $A_L=c_K$. Each lacuna $L$ from this family contains at most $C_2(n)$ contacting cubes. For each such a cube, say $H$, there are at most $C_3(n)$ cubes from other lacunae which contact with $H$. Finally, we obtain at most
%----------------------------------------------------------
$$
\#\{Q\in\Qc_1: Q'=K\}\le C_1(n)\,C_2(n)\,C_3(n).
$$
%----------------------------------------------------------
The same estimate is true for $Q''$.\bigskip
%----------------------------------------------------------
\par Let us estimate $\|f_2\|_{\LPM}$. We again follow the prove of Theorem \reff{REF-MAIN-CR}, see \rf{F2-I12}. We modify this inequality in the same fashion as we did this for the quantity $I_1$. We obtain
%----------------------------------------------------------
\bel{F2-RN}
\intl_{\RN}|f_2|^p\,d\mu\le C\{\tI_1+I_2\}
\ee
%----------------------------------------------------------
where the quantities $\tI_1$ and $I_2$ are defined by \rf{I-N} and \rf{I2-D} respectively.
%----------------------------------------------------------
\par We know that $\tI_1$ can be estimated via the family $\Qc_1$ which we have constructed below. Let us estimate the quantity
%----------------------------------------------------------
\bel{I2-ST}
I_2=\smed_{K\in\Ac}\,\intl_K |f-\tf_1(a_K)|^p\,d\mu+\smed_{K\in\Kc_E}\, \intl_{K}|f-\av{K}|^p\,d\mu=I_3+I_4.
\ee
%----------------------------------------------------------
\par First we estimate the quantity
%----------------------------------------------------------
$$
I_4:=\smed_{K\in \Kc_E}\,\intl_{K}|f-\av{K}|^p\,d\mu.
$$
%----------------------------------------------------------
Let $\Qc_2:=\Kc_E.$ To each $Q\in\Qc_2$ we assign cubes $Q'$ and $Q''$ by letting $Q'=Q''=Q$. Then
%----------------------------------------------------------
$$
I_4:=\smed_{Q\in\Qc_2}\,\intl_{Q'} |f-f_{Q''}|^p\,d\mu
\le \smed_{Q\in\Qc_2}\,\frac{1}{\mu(Q')}
\iint\limits_{Q'\times Q''}|f(x)-f(y)|^p\,d\mu(x)d\mu(y).
$$
%----------------------------------------------------------
Since $Q',Q''\in\Kc_E$, we have
%----------------------------------------------------------
$$
\mu(Q')\sim (\diam Q')^{n-p},~~~\mu(Q'')\sim (\diam Q'')^{n-p},
$$
%----------------------------------------------------------
see Corollary \reff{PR-5K}. Hence
%----------------------------------------------------------
$$
I_4:=\smed_{Q\in\Qc_2}\,\intl_{Q'} |f-f_{Q''}|^p\,d\mu
\le \smed_{Q\in\Qc_2}\,\frac{1}{\mu(Q')}
\iint\limits_{Q'\times Q''}|f(x)-f(y)|^p\,d\mu(x)d\mu(y)
$$
%----------------------------------------------------------
so that
%----------------------------------------------------------
\bel{I4-TY}
I_4\le C\shuge_{Q\in\Qc_2}\,\,
\frac{(\diam Q)^{n-p}\iint \limits_{Q'\times Q''}
|f(x)-f(y)|^p\, d\mu(x)d\mu(y)}
{ \{(\diam Q')^{n-p}+\mu(Q')\} \{(\diam Q'')^{n-p}+\mu(Q'')\}}.
\ee
%----------------------------------------------------------
Clearly, $Q',Q''\subset\gamma Q$ for every $Q\in\Qc_2$ with $\gamma=1$. Then, by Proposition \reff{OP-W}, inequality \rf{OP-F} remains true after replacement the family $\Qc_1$ by $\Qc_2$.
%----------------------------------------------------------
\par It is also clear that $\Qc_1\cap\Qc_2=\emp$, and the (identical) mappings $\Qc_2\ni Q\mapsto Q'\in\Kc_E$ and $\Qc_2\ni Q\mapsto Q''\in\Kc_E$ are one-to-one mappings.
%----------------------------------------------------------
\par We turn to the quantity
%----------------------------------------------------------
$$
I_3=\smed_{K\in\Ac}\,\intl_K |f-\tf_1(a_K)|^p\,d\mu.
$$
%----------------------------------------------------------
Since covering multiplicity $M(W_E)\le C(n)$, we have
%----------------------------------------------------------
$$
I_3\le C(n)\smed_{L\in\Lc_E}\,\,\intl_{\tU_L} |f-\tf_1(A_L)|^p\,d\mu
$$
%----------------------------------------------------------
where
%----------------------------------------------------------
$$
\tU_L:=\cup\{Q: Q\in L\cap \Ac\}.
$$
%----------------------------------------------------------
We also recall that
%----------------------------------------------------------
$$
\tf_1(A_L)=\tfrac{1}{\mu(K_L)}
\intl_{K_L}f\,d\mu
$$
%----------------------------------------------------------
where $K_L\in\Kc_E$ is the unique cube such that $A_L=c_{K_L}$. Hence
%----------------------------------------------------------
\bel{I3-J3}
I_3\le C\,J_3
\ee
%----------------------------------------------------------
where
%----------------------------------------------------------
\bel{J3-DF}
J_3:=\smed_{L\in\Lc_E}\,\frac{1}{\mu(K_L)}
\iint\limits_{\tU_L\times K_L}|f(x)-f(y)|^p\,d\mu(x)d\mu(y).
\ee
%----------------------------------------------------------
\par By definition of the family $\Ac$, see \rf{A-D}, for every cube $Q\in\tU_L$ we have $(\eta K_L)\cap Q=\emp$. Hence
%----------------------------------------------------------
\bel{TR-4}
\dist(Q,A_L)\ge\eta\diam K_L.
\ee
%----------------------------------------------------------
On the other hand, by \rf{AL-2}, $A_L\subset \gamma Q$ so that
%----------------------------------------------------------
\bel{P-7}
\dist(Q,A_L)\le\gamma\diam Q
\ee
%----------------------------------------------------------
proving that
%----------------------------------------------------------
$$
\diam K_L\le C\diam Q.
$$
%----------------------------------------------------------
These inequalities also imply that
%----------------------------------------------------------
\bel{KG3}
K_L\subset \gamma_3 Q,~~~Q\in \tU_L.
\ee
%----------------------------------------------------------
\par Let us note that the family of sets
%----------------------------------------------------------
$$
J:=\{\tU_L\cup K_L:L\in\Lc_E\}
$$
%----------------------------------------------------------
has covering multiplicity $M(J)\le C(n)$.\medskip
%----------------------------------------------------------
\par To complete the proof of the theorem we need the following
%----------------------------------------------------------
%@@@@@@@@@@@@@@@@@@@@@@@@@@@@@@@@@@@@@@@@@@@@@@@@@@@@@@@@@@
%@@@@@@@@@@@@@@@@@@@@@@@@@@@@@@@@@@@@@@@@@@@@@@@@@@@@@@@@@@
%@@@@@@@@@@@@@@@@@@@@@@@@@@@@@@@@@@@@@@@@@@@@@@@@@@@@@@@@@@
%@@@@@@@@@@@@@@@@@@@@@@@@@@@@@@@@@@@@@@@@@@@@@@@@@@@@@@@@@@
%@@@@@@@@@@@@@@@@@@@@@@@@@@@@@@@@@@@@@@@@@@@@@@@@@@@@@@@@@@
%@@@@@@@@@@@@@@@@@@@@@@@@@@@@@@@@@@@@@@@@@@@@@@@@@@@@@@@@@@
%----------------------------------------------------------
\begin{lemma}\lbl{N-L7} Let $L\in\Lc_E$ and let $f\in\Sigma$. Then
$J_3\le C(n,p)\|f\|_{\Sigma}^p.$
%----------------------------------------------------------
\end{lemma}
%----------------------------------------------------------
%@@@@@@@@@@@@@@@@@@@@@@@@@@@@@@@@@@@@@@@@@@@@@@@@@@@@@@@@@@
%@@@@@@@@@@@@@@@@@@@@@@@@@@@@@@@@@@@@@@@@@@@@@@@@@@@@@@@@@@
%@@@@@@@@@@@@@@@@@@@@@@@@@@@@@@@@@@@@@@@@@@@@@@@@@@@@@@@@@@
%@@@@@@@@@@@@@@@@@@@@@@@@@@@@@@@@@@@@@@@@@@@@@@@@@@@@@@@@@@
%@@@@@@@@@@@@@@@@@@@@@@@@@@@@@@@@@@@@@@@@@@@@@@@@@@@@@@@@@@
%@@@@@@@@@@@@@@@@@@@@@@@@@@@@@@@@@@@@@@@@@@@@@@@@@@@@@@@@@@
%----------------------------------------------------------
\par {\it Proof.} Let $f=f_1+f_2$ where $f_1\LOP$ and $f_2\in\LPM$. Then
%----------------------------------------------------------
\be
J(f;L)&:=&\frac{1}{\mu(K_L)}
\iint\limits_{\tU_L\times K_L}
|f(x)-f(y)|^p\,d\mu(x)d\mu(y)\nn\\
&\le& 2^p\left\{\frac{1}{\mu(K_L)}
\iint\limits_{\tU_L\times K_L}
|f_1(x)-f_1(y)|^p\,d\mu(x)d\mu(y)\right.\nn\\
&+&\left.\frac{1}{\mu(K_L)}
\iint\limits_{\tU_L\times K_L}|f_2(x)-f_2(y)|^p\,d\mu(x)d\mu(y)\right\}\nn\\
&=&2^p\{I_1(L)+I_2(L)\}.\nn
\ee
%----------------------------------------------------------
\par Let $q:=(p+n)/2$. By \rf{KG3}, $K_L\subset\gamma Q$ for every $Q\in \tU_L$ so that $Q\cup K_L\subset \gamma_1 Q$. Then, by the Sobolev-Poincar\'e inequality, see Proposition \reff{EDIFF}, for every $x,y\in Q\cup K_L$ we have
%----------------------------------------------------------
$$
|f(x)-f(y)|^p\le C(n,q)\, (\diam Q)^p \left(\frac{1}{|\gamma_1 Q|}\intl_{\gamma_1 Q}
\|\nabla f(z)\|^q\,dz\right) ^{\frac{p}{q}}.
$$
%----------------------------------------------------------
Hence
%----------------------------------------------------------
\be
S_1(Q)&:=&\frac{1}{\mu(K_L)}
\iint\limits_{Q\times K_L}
|f_1(x)-f_1(y)|^p\,d\mu(x)d\mu(y)\nn\\&\le& C\,\mu(Q)
(\diam Q)^p \left(\frac{1}{|\gamma_1 Q|}\intl_{\gamma_1 Q}
\|\nabla f(z)\|^q\,dz\right) ^{\frac{p}{q}}.\nn
\ee
%----------------------------------------------------------
Recall that $K_L\in\Kc_E$ so that, by Corollary \reff{PR-5K},
%----------------------------------------------------------
$$
\mu(K_L)\sim (\diam K_L)^{n-p}.
$$
%----------------------------------------------------------
Since $Q\in W_E$, by Corollary \reff{P-WKE},
%----------------------------------------------------------
\bel{M-QP}
\mu(Q)\le C(\diam Q)^{n-p}.
\ee
%----------------------------------------------------------
Hence
%----------------------------------------------------------
\be
S_1(Q)&\le& C\,(\diam Q)^{n-p}
(\diam Q)^p \left(\frac{1}{|\gamma_1 Q|}\intl_{\gamma_1 Q}
\|\nabla f(z)\|^q\,dz\right) ^{\frac{p}{q}}\nn\\
&\le& C\,|Q|\,\left(\frac{1}{|\gamma_1 Q|}
\intl_{\gamma_1 Q}
\|\nabla f(z)\|^q\,dz\right) ^{\frac{p}{q}}.\nn
\ee
%----------------------------------------------------------
By \rf{MA-S1},
%----------------------------------------------------------
$$
S_1(Q)\le C\,|Q|\,\left(\frac{1}{|\gamma_1 Q|}
\intl_{\gamma_1 Q}
\|\nabla f(z)\|^q\,dz\right) ^{\frac{p}{q}}\le C\,\intl_{\tfrac12 Q}\Mc[\,\|\nabla f_1(z)\|^q]^{\frac{p}{q}}(z)\,dz.
$$
%----------------------------------------------------------
Hence
%----------------------------------------------------------
\be
I_1(L)&:=&\frac{1}{\mu(K_L)}
\iint\limits_{\tU_L\times K_L}
|f_1(x)-f_1(y)|^p\,d\mu(x)d\mu(y)\nn\\&\le&
\smed_{Q\in\tU_L}\,S_1(Q)\le C\intl_{U_L}\Mc[\,\|\nabla f_1(z)\|^q]^{\frac{p}{q}}(z)\,dz.\nn
\ee
%----------------------------------------------------------
Recall that $U_L=\cup\{Q:Q\in L\}$.
%----------------------------------------------------------
\par Since the sets $\{U_L:L\in\Lc_E\}$ are pairwise disjoint, we obtain
%----------------------------------------------------------
$$
A_1:=\smed_{L\in\Lc_E}I_1(L)\le C\, \smed_{L\in\Lc_E}\,
\intl_{U_L}\Mc[\,\|\nabla f_1(z)\|^q]^{\frac{p}{q}}(z)\,dz\le C\,\intl_{\RN}\Mc[\,\|\nabla f_1(z)\|^q]^{\frac{p}{q}}(z)\,dz.
$$
%----------------------------------------------------------
Since $p>q$, by the Hardy-Littlewood maximal theorem,
%----------------------------------------------------------
$$
A_1\le C\,\intl_{\RN}
(\|\nabla f_1(z)\|^q)^{\frac{p}{q}}(z)\,dz=
C\intl_{\RN}
\|\nabla f_1(z)\|^p(z)\,dz=\|f_1\|_{\LOP}^p.
$$
%----------------------------------------------------------
\medskip
%----------------------------------------------------------
\par Let us estimate the quantity
%----------------------------------------------------------
$$
I_2(L):=\frac{1}{\mu(K_L)}
\iint\limits_{\tU_L\times K_L}
|f_2(x)-f_2(y)|^p\,d\mu(x)d\mu(y).
$$
%----------------------------------------------------------
We have
%----------------------------------------------------------
$$
I_2(L)\le \frac{2^p}{\mu(K_L)}\left\{
\iint\limits_{\tU_L\times K_L}
|f_2(x)|^p\,d\mu(x)d\mu(y)+\iint\limits_{\tU_L\times K_L}
|f_2(y)|^p\,d\mu(x)d\mu(y)\right\}.
$$
%----------------------------------------------------------
Hence
%----------------------------------------------------------
$$
I_2(L)\le 2^p\left\{
\intl_{\tU_L}
|f_2(x)|^p\,d\mu(x)+\frac{\mu(\tU_L)}{\mu(K_L)}\intl_{K_L}
|f_2(y)|^p\,d\mu(x)\right\}.
$$
%----------------------------------------------------------
\par Prove that
%----------------------------------------------------------
\bel{M-9}
\mu(\tU_L)\le C\mu(K_L).
\ee
%----------------------------------------------------------
By \rf{M-QP},
%----------------------------------------------------------
$$
\mu(\tU_L)\le \smed_{Q\in L\cap\Ac}\mu(Q)\le C\smed_{Q\in L\cap\Ac}(\diam Q)^{n-p}.
$$
%----------------------------------------------------------
\par Note that for every $Q\in L\cap\Ac$ and every $x\in Q$ we have
%----------------------------------------------------------
\bel{X-AQ}
\|x-A_L\|\sim \diam Q.
\ee
%----------------------------------------------------------
In fact, since $Q\in W_E$, its diameter $\diam Q\sim \dist(Q,E)$ so that $\diam Q\le C\dist(Q,E)$. Hence
%----------------------------------------------------------
$$
\diam Q\le C\|x-A_L\|.
$$
%----------------------------------------------------------
On the other hand, by \rf{P-7},
%----------------------------------------------------------
$$
\|x-A_L\|\le \diam Q+\dist(Q,A_L)\le C\,\diam Q
$$
%----------------------------------------------------------
proving \rf{X-AQ}.
%----------------------------------------------------------
\par Hence
%----------------------------------------------------------
$$
\mu(\tU_L)\le C\smed_{Q\in L\cap\Ac}(\diam Q)^{n-p}
\le C\smed_{Q\in L\cap\Ac}\,\,\intl_Q\frac{dx}{\|x-A_L\|^p}.
$$
%----------------------------------------------------------
Note that, by \rf{TR-4}, we have
%----------------------------------------------------------
$$
\|x-A_L\|\ge \eta\diam K_L,~~~x\in\tU_L.
$$
%----------------------------------------------------------
We also note that covering multiplicity $M(W_E)\le N(n)$. Hence
%----------------------------------------------------------
$$
\mu(\tU_L)\le C N(n)\intl_{\|x-A_L\|\ge \eta\diam K_L}\frac{dx}{\|x-A_L\|^p}.
$$
%----------------------------------------------------------
We obtain
%----------------------------------------------------------
$$
\mu(\tU_L)\le C (\diam K_L)^{n-p}\sim \mu(K_L)
$$
%----------------------------------------------------------
proving \rf{M-9}.
%----------------------------------------------------------
\par Hence
%----------------------------------------------------------
$$
I_2(L)\le C\left\{
\intl_{\tU_L}
|f_2(x)|^p\,d\mu(x)+\intl_{K_L}
|f_2(y)|^p\,d\mu(x)\right\}.
$$
%----------------------------------------------------------
Recall that the mapping $L\mapsto K_L$ is an ``almost'' one-to-one, so that covering multiplicity of the family $\{K_L:L\in\Lc_E\}$ is bounded by a constant $N(n)$. Since the sets $\{\tU_L:L\in\Lc_E\}$ are pairwise disjoint, we obtain
%----------------------------------------------------------
\be
A_2&:=&\smed_{L\in\Lc_E}I_2(L)\le
C\smed_{L\in\Lc_E}\intl_{\tU_L}
|f_2(x)|^p\,d\mu(x)+C\smed_{L\in\Lc_E}\intl_{K_L}
|f_2(y)|^p\,d\mu(x)\nn\\
&\le& C\,\intl_{\RN}
|f_2(x)|^p\,d\mu(x)+C\,N(n)\intl_{\RN}
|f_2(x)|^p\,d\mu(x)\nn\\
&\le& C\intl_{\RN}
|f_2(x)|^p\,d\mu(x)\nn
\ee
%----------------------------------------------------------
proving that
%----------------------------------------------------------
$$
A_2\le\|f_2\|_{\LPM}^p.
$$
%----------------------------------------------------------
Finally,
%----------------------------------------------------------
\be
J_3&=&\smed_{L\in\Lc_E}\,\frac{1}{\mu(K_L)}
\iint\limits_{\tU_L\times K_L}
|f(x)-f(y)|^p\,d\mu(x)d\mu(y)\nn\\&\le& C(A_1+A_2)\le C(\|f_1\|_{\LOP}+\|f_2\|_{\LPM})^p.\nn
\ee
%----------------------------------------------------------
taking the infimum in this inequality over all functions $f_1\in\LOP$ and $f_2\in\LPM$ such that $f=f_1+f_2$ we obtain the statement of the lemma.\bx
%----------------------------------------------------------
%@@@@@@@@@@@@@@@@@@@@@@@@@@@@@@@@@@@@@@@@@@@@@@@@@@@@@@@@@@
%@@@@@@@@@@@@@@@@@@@@@@@@@@@@@@@@@@@@@@@@@@@@@@@@@@@@@@@@@@
%@@@@@@@@@@@@@@@@@@@@@@@@@@@@@@@@@@@@@@@@@@@@@@@@@@@@@@@@@@
%----------------------------------------------------------
\medskip
%----------------------------------------------------------
\par Let us finish the proof of the theorem.
%----------------------------------------------------------
%@@@@@@@@@@@@@@@@@@@@@@@@@@@@@@@@@@@@@@@@@@@@@@@@@@@@@@@@@@
%@@@@@@@@@@@@@@@@@@@@@@@@@@@@@@@@@@@@@@@@@@@@@@@@@@@@@@@@@@
%@@@@@@@@@@@@@@@@@@@@@@@@@@@@@@@@@@@@@@@@@@@@@@@@@@@@@@@@@@
%----------------------------------------------------------
\par Let us enumerate the cubes of the family $\Qc_1$:
%----------------------------------------------------------
$$
\Qc_1=\{Q_i:i\in\N\}.
$$
%----------------------------------------------------------
Let
%----------------------------------------------------------
$$
D_1(f):=\sum_{i=1}^\infty \,\,
\alpha_i \iint\limits_{Q'_i\times Q''_i}
|f(x)-f(y)|^p\, d\mu(x)\,d\mu(y)
$$
%----------------------------------------------------------
where
%----------------------------------------------------------
$$
\alpha_i:= \frac{(\diam Q_i)^{n-p}}
{\{(\diam Q'_i)^{n-p}+\mu(Q'_i)\} \{(\diam Q''_i)^{n-p}+\mu(Q''_i)\}}
$$
%----------------------------------------------------------
Then, by \rf{F1-R},
%----------------------------------------------------------
$$
\|f_1\|_{\LOP}^p\le C\,D_1(f).
$$
%----------------------------------------------------------
Recall that we have constructed the families of cubes
%----------------------------------------------------------
$$
\Gc_1:=\{Q'_i:i\in\N\}~~~\text{and}~~~
\Hc_1:=\{Q''_i:i\in\N\}
$$
%----------------------------------------------------------
in such a way that $M(\Gc_1),M(\Hc_1)\le C(n)$. We have also proved that
%----------------------------------------------------------
$$
D_1(f)\le C\|f\|_{\Sigma}^p.
$$
%----------------------------------------------------------
\par We have also defined a family $\Qc_2$ of cubes
%----------------------------------------------------------
$$
\Qc_2=\{K_i:i\in\N\}.
$$
%----------------------------------------------------------
Let
%----------------------------------------------------------
$$
D_2(f):=\sum_{i=1}^\infty \,\,
\beta_i \iint\limits_{K'_i\times K''_i}
|f(x)-f(y)|^p\, d\mu(x)\,d\mu(y)
$$
%----------------------------------------------------------
where
%----------------------------------------------------------
$$
\beta_i:= \frac{(\diam K_i)^{n-p}}
{\{(\diam K'_i)^{n-p}+\mu(K'_i)\} \{(\diam K''_i)^{n-p}+\mu(K''_i)\}}.
$$
%----------------------------------------------------------
By \rf{I4-TY},
%----------------------------------------------------------
\bel{I4-D}
I_4\le C\,D_2(f).
\ee
%----------------------------------------------------------
\par We know that
%----------------------------------------------------------
$$
D_2(f)\le C\|f\|_{\Sigma}^p.
$$
%----------------------------------------------------------
Furthermore, we know that for the families of cubes
%----------------------------------------------------------
$$
\Gc_2:=\{K'_i:i\in\N\}~~~\text{and}~~~
\Hc_2:=\{K''_i:i\in\N\}
$$
%----------------------------------------------------------
we have $M(\Gc_2),M(\Hc_2)\le C(n)$.
%----------------------------------------------------------
\par Finally, we have constructed two families of sets
%----------------------------------------------------------
$$
\Gc_3:=\{\tU_L:L\in\Lc_E\}~~~\text{and}~~~
\Hc_3:=\{K_L:L\in\Lc_E\}
$$
%----------------------------------------------------------
with certain properties. Let us enumerate the sets of these families with preservation of the correspondence $\tU_L\lr K_L$:
%----------------------------------------------------------
$$
\Gc_3:=\{S_i:i\in\N\}~~~\text{and}~~~
\Hc_3:=\{T_i:i\in\N\}.
$$
%----------------------------------------------------------
We know that $M(\Gc_3),M(\Hc_3)\le C(n)$. Note that the quantity $J_3$ defined by \rf{J3-DF} can be written in the form $J_3=D_3(f)$ where
%----------------------------------------------------------
$$
D_3(f):=\sum_{i=1}^\infty \,\,
\theta_i \iint\limits_{S_i\times T_i}
|f(x)-f(y)|^p\, d\mu(x)\,d\mu(y)
$$
%----------------------------------------------------------
where
%----------------------------------------------------------
$$
\theta_i:=1/\mu(T_i).
$$
%----------------------------------------------------------
Furthermore, by Lemma \reff{N-L7},
%----------------------------------------------------------
$$
J_3=D_3(f)\le C\|f\|_{\Sigma}.
$$
%----------------------------------------------------------
\par Combining inequalities \rf{F2-RN}, \rf{I2-ST}, \rf{I3-J3} and \rf{I4-D} with definition of $J_3$, see  \rf{J3-DF}, we obtain
%----------------------------------------------------------
$$
\|f_2\|_{\LPM}^p\le C\{D_1(f)+D_2(f)+D_3(f)\}.
$$
%----------------------------------------------------------
Hence,
%----------------------------------------------------------
$$
\|f\|_{\Sigma}^p\le 2^p\{\|f_1\|_{\LOP}^p+\|f_2\|_{\LPM}\}^p\le C\{D_1(f)+D_2(f)+D_3(f)\}
$$
%----------------------------------------------------------
so that
%----------------------------------------------------------
$$
\|f\|_{\Sigma}^p\sim D_1(f)+D_2(f)+D_3(f).
$$
%----------------------------------------------------------
\par Theorem \reff{REF-2} is completely proved.\bx
%----------------------------------------------------------
%@@@@@@@@@@@@@@@@@@@@@@@@@@@@@@@@@@@@@@@@@@@@@@@@@@@@@@@@@@
%@@@@@@@@@@@@@@@@@@@@@@@@@@@@@@@@@@@@@@@@@@@@@@@@@@@@@@@@@@
%@@@@@@@@@@@@@@@@@@@@@@@@@@@@@@@@@@@@@@@@@@@@@@@@@@@@@@@@@@
%@@@@@@@@@@@@@@@@@@@@@@@@@      @@@@@@@@@@@@@@@@@@@@@@@@@@@
%@@@@@@@@@@@@@@@@@@@@@@@          @@@@@@@@@@@@@@@@@@@@@@@@@
%@@@@@@@@@@@@@@@@@@@@@              @@@@@@@@@@@@@@@@@@@@@@@
%@@@@@@@@@@@@@@@@@@@     SECTION 7
%@@@@@@@@@@@@@@@@@@@@@
%@@@@@@@@@@@@@@@@@@@@@              @@@@@@@@@@@@@@@@@@@@@@@
%@@@@@@@@@@@@@@@@@@@@@@@          @@@@@@@@@@@@@@@@@@@@@@@@@
%@@@@@@@@@@@@@@@@@@@@@@@@@      @@@@@@@@@@@@@@@@@@@@@@@@@@@
%@@@@@@@@@@@@@@@@@@@@@@@@@@@@@@@@@@@@@@@@@@@@@@@@@@@@@@@@@@
%@@@@@@@@@@@@@@@@@@@@@@@@@@@@@@@@@@@@@@@@@@@@@@@@@@@@@@@@@@
%----------------------------------------------------------
%@@@@@@@@@@@@@@@@@@@@@@@@@@@@@@@@@@@@@@@@@@@@@@@@@@@@@@@@@@
%----------------------------------------------------------
\SECT{7. Further results and comments.}{7}
%----------------------------------------------------------
\addtocontents{toc}{7. Further results and comments. \hfill \thepage\vspace*{3mm}\par}
%----------------------------------------------------------
\indent
%----------------------------------------------------------
\par {\bf 7.1. Modifications of Theorem \reff{MAIN-CR}.}
%----------------------------------------------------------
\addtocontents{toc}{~~~~7.1. Modifications of Theorem \reff{MAIN-CR}. \hfill \thepage\par}
%----------------------------------------------------------
In this subsection we present several versions of the criterion for calculation of the norm of a function in the space $\SUM$. First of them is Theorem \reff{CR-V1} which we have formulated in Section 1. Its proof is very short.\medskip
%----------------------------------------------------------
\par {\it Proof of Theorem \reff{CR-V1}.} The necessity part of the theorem directly follows from the necessity part of Theorem \reff{MAIN-CR}. In fact, the left hand side of inequality  \rf{IN-V1} is majorized (up to an absolute constant) by the left hand side of \rf{CR} provided inequality \rf{M-DM-V1} holds. In turn, the sufficiency part of the theorem immediately follows from Theorem \reff{S-V2}.\bx\medskip
%----------------------------------------------------------
%@@@@@@@@@@@@@@@@@@@@@@@@@@@@@@@@@@@@@@@@@@@@@@@@@@@@@@@@@@
%@@@@@@@@@@@@@@@@@@@@@@@@@@@@@@@@@@@@@@@@@@@@@@@@@@@@@@@@@@
%@@@@@@@@@@@@@@@@@@@@@@@@@@@@@@@@@@@@@@@@@@@@@@@@@@@@@@@@@@
%@@@@@@@@@@@@@@@@@@@@@@@@@@@@@@@@@@@@@@@@@@@@@@@@@@@@@@@@@@
%@@@@@@@@@@@@@@@@@@@@@@@@@@@@@@@@@@@@@@@@@@@@@@@@@@@@@@@@@@
%@@@@@@@@@@@@@@@@@@@@@@@@@@@@@@@@@@@@@@@@@@@@@@@@@@@@@@@@@@
%@@@@@@@@@@@@@@@@@@@@@@@@@@@@@@@@@@@@@@@@@@@@@@@@@@@@@@@@@@
%----------------------------------------------------------
\begin{remark} {\em Inequality \rf{M-DM-V1} of  Theorem \reff{CR-V1} can be replaced by weaker conditions
%----------------------------------------------------------
\bel{UR}
\mu(Q')\le (\diam Q')^{n-p},~~~
\mu(Q'')\le (\diam Q'')^{n-p}\le C\,\mu(Q''),
\ee
%----------------------------------------------------------
and
%---------------------------------------------------------- %@@@@@@@@@@@@@@@@@@@@@@@@@@@@@@@@@@@@@@@@@@@@@@@@@@@@@@@@@@
%----------------------------------------------------------
\bel{C-QQ}
\mu(Q')\le C\,\mu(Q'')
\ee
%----------------------------------------------------------
where $C=C(n,p)$ is a constant depending only on $n$ and $p$.
%----------------------------------------------------------
\par Note that \rf{C-QQ} is equivalent to the inequality
%----------------------------------------------------------
$$
\diam Q''\le C\,\diam Q'
$$
%----------------------------------------------------------
provided the inequalities in \rf{UR} hold.
%----------------------------------------------------------
\par In fact, in all our definitions of the cube $Q''$ (see Propositions \reff{EA1}, \reff{EA2} and \reff{EA4}) we have $Q''\in\Kc_E$ so that $\mu(Q'')\sim (\diam Q'')^{n-p}$, see \rf{M-DK}.
%----------------------------------------------------------
\par On the other hand, in Propositions \reff{EA1} and  \reff{EA2} both $Q'$ and $Q''$  belong to $\Kc_E$ so that in this case without loss of generality we may assume that \rf{C-QQ} holds with $C=1$. In the same way we can define the cubes $Q'$ and $Q''$ in Propositions \reff{EA4} whenever $Q=Q'=Q''\in\Kc_E$, see \rf{R1-QP}. %----------------------------------------------------------
\par In the remaining case, see \rf{RM-2},
%----------------------------------------------------------
$$
Q'=K\in \Ac~~~\text{and}~~~Q''=K^{(a_K)}\in\Kc_E.
$$
%----------------------------------------------------------
Then, by Lemma \reff{TQ}, $Q''\subset(22\tau^2) Q'$ so that $\diam Q''\le C\diam Q'$. But, by \rf{M-DK} and \rf{WKE-S},
%----------------------------------------------------------
$$\mu(Q'')\sim (\diam Q'')^{n-p}~~~\text{and}~~~\mu(Q')\le C(\diam Q')^{n-p}
$$
%----------------------------------------------------------
proving \rf{C-QQ}.\rbx}
%----------------------------------------------------------
\end{remark}
%----------------------------------------------------------
%@@@@@@@@@@@@@@@@@@@@@@@@@@@@@@@@@@@@@@@@@@@@@@@@@@@@@@@@@@
%@@@@@@@@@@@@@@@@@@@@@@@@@@@@@@@@@@@@@@@@@@@@@@@@@@@@@@@@@@
%@@@@@@@@@@@@@@@@@@@@@@@@@@@@@@@@@@@@@@@@@@@@@@@@@@@@@@@@@@
%@@@@@@@@@@@@@@@@@@@@@@@@@@@@@@@@@@@@@@@@@@@@@@@@@@@@@@@@@@
%@@@@@@@@@@@@@@@@@@@@@@@@@@@@@@@@@@@@@@@@@@@@@@@@@@@@@@@@@@
%@@@@@@@@@@@@@@@@@@@@@@@@@@@@@@@@@@@@@@@@@@@@@@@@@@@@@@@@@@
%@@@@@@@@@@@@@@@@@@@@@@@@@@@@@@@@@@@@@@@@@@@@@@@@@@@@@@@@@@
%@@@@@@@@@@@@@@@@@@@@@@@@@@@@@@@@@@@@@@@@@@@@@@@@@@@@@@@@@@
%----------------------------------------------------------
\begin{remark}\lbl{RM-V2} {\em Let us replace inequality \rf{CR} in Theorem \reff{MAIN-CR} with the following one:
%----------------------------------------------------------
\bel{V-TH3}
\shuge_{Q\in\Qc}\,\,
\frac{\frac{1}{\mu(Q')\mu(Q'')}\iint \limits_{Q'\times Q''}
|f(x)-f(y)|^p\, d\mu(x)d\mu(y)}
{(\diam Q)^{p-n}\{1+(\diam Q')^{n-p}/\mu(Q')
+(\diam Q'')^{n-p}/\mu(Q'')\}}
\le \lambda
\ee
%----------------------------------------------------------
%@@@@@@@@@@@@@@@@@@@@@@@@@@@@@@@@@@@@@@@@@@@@@@@@@@@@@@@@@@
%@@@@@@@@@@@@@@@@@@@@@@@@@@@@@@@@@@@@@@@@@@@@@@@@@@@@@@@@@@
%----------------------------------------------------------
Then the result of Theorem \reff{MAIN-CR} remains true after such a modification. Thus we obtain another criterion for calculation of the norm in the space $\SUM$.
%----------------------------------------------------------
\par In fact, the necessity part of this new criterion follows from Proposition \reff{OP-W} (with $\Sc=\Qc$, $S_{Q'}=Q'$ and $S_{Q''}=Q''$). In turn, the sufficiency directly follows from the sufficiency part of Theorem \reff{MAIN-CR} because the left hand side of inequality \rf{CR} does not exceed the left hand side of \rf{V-TH3}. See definition \rf{AQQP} and inequality \rf{INAQ}.\rbx}
%----------------------------------------------------------
\end{remark}\medskip
%----------------------------------------------------------
%@@@@@@@@@@@@@@@@@@@@@@@@@@@@@@@@@@@@@@@@@@@@@@@@@@@@@@@@@@
%@@@@@@@@@@@@@@@@@@@@@@@@@@@@@@@@@@@@@@@@@@@@@@@@@@@@@@@@@@
%----------------------------------------------------------
\par The criterion \rf{V-TH3} and previous results lead us to the following result formulated in the spirit of Theorem \reff{CR-V1}.
%----------------------------------------------------------
%@@@@@@@@@@@@@@@@@@@@@@@@@@@@@@@@@@@@@@@@@@@@@@@@@@@@@@@@@@
%@@@@@@@@@@@@@@@@@@@@@@@@@@@@@@@@@@@@@@@@@@@@@@@@@@@@@@@@@@
%@@@@@@@@@@@@@@@@@@@@@@@@@@@@@@@@@@@@@@@@@@@@@@@@@@@@@@@@@@
%@@@@@@@@@@@@@@@@@@@@@@@@@@@@@@@@@@@@@@@@@@@@@@@@@@@@@@@@@@
%@@@@@@@@@@@@@@@@@@@@@@@@@@@@@@@@@@@@@@@@@@@@@@@@@@@@@@@@@@
%@@@@@@@@@@@@@@@@@@@@@@@@@@@@@@@@@@@@@@@@@@@@@@@@@@@@@@@@@@
%----------------------------------------------------------
\begin{theorem}\lbl{CR-V4} Let us replace inequality \rf {IN-V1} in the formulation of Theorem \reff{CR-V1} by the inequality
%----------------------------------------------------------
\bel{V4-K}
\shuge_{Q\in\Qc}\,\,
\left(\frac{\diam Q' \diam Q''}{\diam Q}\right)^{p-n} \frac{\iint \limits_{Q'\times Q''}
|f(x)-f(y)|^p\, d\mu(x)d\mu(y)}
{(\diam Q')^{p-n}\mu(Q')+(\diam Q'')^{p-n}\mu(Q'')}
\le \lambda.
\ee
%----------------------------------------------------------
Then, after such a modification the result of Theorem \reff{CR-V1} remains true.
%----------------------------------------------------------
\end{theorem}
%----------------------------------------------------------
%@@@@@@@@@@@@@@@@@@@@@@@@@@@@@@@@@@@@@@@@@@@@@@@@@@@@@@@@@@
%@@@@@@@@@@@@@@@@@@@@@@@@@@@@@@@@@@@@@@@@@@@@@@@@@@@@@@@@@@
%@@@@@@@@@@@@@@@@@@@@@@@@@@@@@@@@@@@@@@@@@@@@@@@@@@@@@@@@@@
%@@@@@@@@@@@@@@@@@@@@@@@@@@@@@@@@@@@@@@@@@@@@@@@@@@@@@@@@@@
%@@@@@@@@@@@@@@@@@@@@@@@@@@@@@@@@@@@@@@@@@@@@@@@@@@@@@@@@@@
%@@@@@@@@@@@@@@@@@@@@@@@@@@@@@@@@@@@@@@@@@@@@@@@@@@@@@@@@@@
%----------------------------------------------------------
\par {\it Proof.} Clearly, by \rf{M-DM-V1}, the left hand side of \rf{IN-V1} is smaller than the left hand side of \rf{V4-K} so that the sufficiency follows from the sufficiency part of Theorem \reff{CR-V1}. On the other hand the left hand side of \rf{V4-K} is smaller (up to an absolute constant) than  to the
the left hand side of \rf{V-TH3} provided inequality \rf{M-DM-V1} is satisfied. But as we have seen in Remark \reff{RM-V2}, the necessity of \rf{V-TH3} follows from Proposition \reff{OP-W}. This proves the necessity part of Theorem \reff{CR-V4}.\bx
%----------------------------------------------------------
%@@@@@@@@@@@@@@@@@@@@@@@@@@@@@@@@@@@@@@@@@@@@@@@@@@@@@@@@@@
%@@@@@@@@@@@@@@@@@@@@@@@@@@@@@@@@@@@@@@@@@@@@@@@@@@@@@@@@@@
%----------------------------------------------------------
\medskip
%----------------------------------------------------------
%@@@@@@@@@@@@@@@@@@@@@@@@@@@@@@@@@@@@@@@@@@@@@@@@@@@@@@@@@@
%@@@@@@@@@@@@@@@@@@@@@@@@@@@@@@@@@@@@@@@@@@@@@@@@@@@@@@@@@@
%@@@@@@@@@@@@@@@@@@@@@@@@@@@@@@@@@@@@@@@@@@@@@@@@@@@@@@@@@@
%@@@@@@@@@@@@@@@@@@@@@@@@@@@@@@@@@@@@@@@@@@@@@@@@@@@@@@@@@@
%----------------------------------------------------------
\begin{remark} {\em In all modifications of the main result we may assume that the cubes $Q',Q''$ belong to a certain family $\tQc$ of pairwise disjoint cubes which my be different from the family $\Qc$. For instance, Theorem \reff{CR-V1} can be modified in the following way:
%----------------------------------------------------------
%@@@@@@@@@@@@@@@@@@@@@@@@@@@@@@@@@@@@@@@@@@@@@@@@@@@@@@@@@@
%@@@@@@@@@@@@@@@@@@@@@@@@@@@@@@@@@@@@@@@@@@@@@@@@@@@@@@@@@@
%@@@@@@@@@@@@@@@@@@@@@@@@@@@@@@@@@@@@@@@@@@@@@@@@@@@@@@@@@@
%@@@@@@@@@@@@@@@@@@@@@@@@@@@@@@@@@@@@@@@@@@@@@@@@@@@@@@@@@@
%@@@@@@@@@@@@@@@@@@@@@@@@@@@@@@@@@@@@@@@@@@@@@@@@@@@@@@@@@@
%@@@@@@@@@@@@@@@@@@@@@@@@@@@@@@@@@@@@@@@@@@@@@@@@@@@@@@@@@@
%----------------------------------------------------------
\begin{theorem} A function $f\in L_{p,loc}(\RN;\mu)$ belongs to the space $\LOP+\LPM$, $n<p<\infty$, if and only if there exists a constant $\lambda>0$ which satisfies the following conditions for a certain absolute constant$\gamma$: Let $\Qc$ and $\tQc$ be arbitrary finite families of pairwise disjoint cubes in $\RN$. Suppose that to each cube $Q\in\Qc$ we have arbitrarily assigned two cubes $Q',Q''\in\tQc$ such that
$Q'\cup Q''\subset \gamma Q$ and inequality \rf{M-DM-V1} is satisfied.
%----------------------------------------------------------
\par Then inequality \rf{IN-V1} holds. Furthermore,
$\|f\|_{\sum}\sim \inf \lambda^{\frac{1}{p}}$
with constants of equivalence depending only on $n$ and $p$.
%----------------------------------------------------------
\end{theorem}
%----------------------------------------------------------
%@@@@@@@@@@@@@@@@@@@@@@@@@@@@@@@@@@@@@@@@@@@@@@@@@@@@@@@@@@
%@@@@@@@@@@@@@@@@@@@@@@@@@@@@@@@@@@@@@@@@@@@@@@@@@@@@@@@@@@
%@@@@@@@@@@@@@@@@@@@@@@@@@@@@@@@@@@@@@@@@@@@@@@@@@@@@@@@@@@
%----------------------------------------------------------
\par ({\it Necessity}). We apply  Proposition \reff{OP-W} to $\Sc=\tQc$, $S_{Q'}=Q'$ and  $S_{Q''}=Q''$ and prove that inequality \rf{V-TH3} holds. As we have noted in Remark \reff{RM-V2}, the left hand side of inequality \rf{CR} does not exceed the left hand side of \rf{V-TH3} which proves the necessity.
%----------------------------------------------------------
\par ({\it Sufficiency}). The sufficiency follows from the sufficiency part of Theorem \reff{MAIN-CR} which is proven for the case $\tQc=\Qc$.\rbx}
%----------------------------------------------------------
\end{remark}
%----------------------------------------------------------
%@@@@@@@@@@@@@@@@@@@@@@@@@@@@@@@@@@@@@@@@@@@@@@@@@@@@@@@@@@
%@@@@@@@@@@@@@@@@@@@@@@@@@@@@@@@@@@@@@@@@@@@@@@@@@@@@@@@@@@
%@@@@@@@@@@@@@@@@@@@@@@@@@@@@@@@@@@@@@@@@@@@@@@@@@@@@@@@@@@
%@@@@@@@@@@@@@@@@@@@@@@@@@@@@@@@@@@@@@@@@@@@@@@@@@@@@@@@@@@
%----------------------------------------------------------
%@@@@@@@@@@@@@@@@@@@@@@@@@@@@@@@@@@@@@@@@@@@@@@@@@@@@@@@@@@
%@@@@@@@@@@@@@@@@@@@@@@@@@@@@@@@@@@@@@@@@@@@@@@@@@@@@@@@@@@
%@@@@@@@@@@@@@@@@@@@@@@@@@@@@@@@@@@@@@@@@@@@@@@@@@@@@@@@@@@
%@@@@@@@@@@@@@@@@@@@@@@@@@@@@@@@@@@@@@@@@@@@@@@@@@@@@@@@@@@
%----------------------------------------------------------
\begin{remark} {\em As we have mentioned in Section 1, Theorem \reff{MAIN-CR} and its variants  have important and interesting applications to the Whitney-type problems of characterizations of restrictions of Sobolev functions to subsets of $\RN$. In particular, in \cite{S6}
we need a variant of Theorem \reff{CR-V1} formulated in terms of families of Euclidean balls rather than cubes.
%----------------------------------------------------------
%@@@@@@@@@@@@@@@@@@@@@@@@@@@@@@@@@@@@@@@@@@@@@@@@@@@@@@@@@@
%@@@@@@@@@@@@@@@@@@@@@@@@@@@@@@@@@@@@@@@@@@@@@@@@@@@@@@@@@@
%@@@@@@@@@@@@@@@@@@@@@@@@@@@@@@@@@@@@@@@@@@@@@@@@@@@@@@@@@@
%@@@@@@@@@@@@@@@@@@@@@@@@@@@@@@@@@@@@@@@@@@@@@@@@@@@@@@@@@@
%@@@@@@@@@@@@@@@@@@@@@@@@@@@@@@@@@@@@@@@@@@@@@@@@@@@@@@@@@@
%@@@@@@@@@@@@@@@@@@@@@@@@@@@@@@@@@@@@@@@@@@@@@@@@@@@@@@@@@@
%----------------------------------------------------------
\begin{theorem} Let $n<p<\infty$ and let $\mu$ be a non-trivial non-negative Borel measure on $\RN$.  A function $f\in L_{p,loc}(\RN;\mu)$ belongs to the space $\LOP+\LPM$ if and only if there exists a constant $\lambda>0$ which satisfy all of the following conditions for a certain absolute positive constant $\gamma$: Let $\Bc$ be an arbitrary finite family of pairwise disjoint balls in $\RN$. Suppose that to each ball $B\in\Bc$ we have arbitrarily assigned two balls $B',B''\in\Bc$ such that $B'\cup B''\subset \gamma B$ and
%----------------------------------------------------------
$$
(\diam B')^{p-n}\mu(B')+(\diam B'')^{p-n}\mu(B'')\le 1.
$$
%----------------------------------------------------------
Then the following inequality
%----------------------------------------------------------
$$
\sbig_{B\in\Qc}\,\,
\left(\frac{\diam B' \diam B''}{\diam B}\right)^{p-n} \iint \limits_{B'\times B''}
|f(x)-f(y)|^p\, d\mu(x)d\mu(y)
\le \lambda
$$
%----------------------------------------------------------
holds.
%----------------------------------------------------------
\par Furthermore,
$\|f\|_{\sum}\sim \inf \lambda^{\frac{1}{p}}$ with constants of equivalence depending only on $n$ and $p$.
%----------------------------------------------------------
\end{theorem}
%----------------------------------------------------------
%@@@@@@@@@@@@@@@@@@@@@@@@@@@@@@@@@@@@@@@@@@@@@@@@@@@@@@@@@@
%@@@@@@@@@@@@@@@@@@@@@@@@@@@@@@@@@@@@@@@@@@@@@@@@@@@@@@@@@@
%@@@@@@@@@@@@@@@@@@@@@@@@@@@@@@@@@@@@@@@@@@@@@@@@@@@@@@@@@@
%----------------------------------------------------------
\par {\it A sketch of the proof.} The proof follows precisely the scheme of the proof of Theorem \reff{CR-V1}. \par There is only one place in this scheme where we have to slightly change formulations of corresponding results. We mean an analogue of the Whitney covering Theorem \reff{Wcov} for Euclidean balls. Of course, in this case we can not cover the open set $\RN\setminus E$ by {\it non-overlapping} balls $B$ such that
$\diam B\sim\dist(B,E)$. Nevertheless for our purpose it suffice to cover $\RN\setminus E$ by a family $\widetilde{W}_E$ of balls whose {\it covering multiplicity} is bounded by a constant $N=N(n)$ depending only on $n$. In other words, every point $x\in\RN$ is covered at most $N$ balls from the family $\widetilde{W}_E$.
%----------------------------------------------------------
\par The existence of a Whitney-type covering of such a kind follows from a general result proven by M. Guzman \cite{G}. (Note that this result relies on the Besicovitch  covering theorem \cite{Be}.)\rbx}
%----------------------------------------------------------
\end{remark}
%----------------------------------------------------------
%@@@@@@@@@@@@@@@@@@@@@@@@@@@@@@@@@@@@@@@@@@@@@@@@@@@@@@@@@@
%@@@@@@@@@@@@@@@@@@@@@@@@@@@@@@@@@@@@@@@@@@@@@@@@@@@@@@@@@@
%@@@@@@@@@@@@@@@@@@@@@@@@@@@@@@@@@@@@@@@@@@@@@@@@@@@@@@@@@@
%@@@@@@@@@@@@@@@@@@@@@@@@@@@@@@@@@@@@@@@@@@@@@@@@@@@@@@@@@@
%----------------------------------------------------------
\medskip
%----------------------------------------------------------
\par {\bf 7.2. The $K$-functional for the couple   $\vec{A}=(\LPM,\LOP)$.}
%----------------------------------------------------------
\addtocontents{toc}{~~~~7.2. The $K$-functional for the couple $\vec{A}=(\LPM,\LOP)$. \hfill \thepage\par}
%----------------------------------------------------------
\par Theorem \reff{MAIN-CR} and its modifications presented in the previous subsection enable us to give various explicit formulas for the $K$-functional of the Banach couple
%----------------------------------------------------------
$$
\vec{A}=(\LPM,\LOP).
$$
%----------------------------------------------------------
We recall that, for each $t>0$
%----------------------------------------------------------
$$
K(t; f :\vec{A}):=\inf\{\|f_1\|_{\LPM}+t\|f_2\|_{\LOP}:
f_1+f_2=f, f_1\in\LPM,f_2 \in \LOP\}
$$
%----------------------------------------------------------
so that $\|f\|_{\sum}=K(1;f:\vec{A})$ and
%----------------------------------------------------------
\bel{F-KFN}
K(t;f:\vec{A})= t\|f\|_{\ssmall_t}
\ee
%----------------------------------------------------------
where
%----------------------------------------------------------
\bel{S-KFN}
\mathlarger{\Sigma}_t:=\LOP+L_p(\RN;\tfrac{1}{t^p}\mu).
\ee
%----------------------------------------------------------
See Remark \reff{KF}.
%----------------------------------------------------------
\par In Section 1 we have presented such a formula for $K(\cdot;f:\vec{A})$. This result directly follows from Theorem \reff{CR-V4} and equalities \rf{F-KFN} and \rf{S-KFN}.
%----------------------------------------------------------
\par Let us prove that the $K$-functional of the couple $\vec{A}$ can be {\it quasi-linearized}, i.e., for each $t>0$ there exist continuous linear operators
%----------------------------------------------------------
$$
T_1[t]:\ssmall\to \LOP~~~~\text{and}~~~~T_2[t]:\ssmall\to \LPM
$$
%----------------------------------------------------------
such that
%----------------------------------------------------------
$$
T_1[t]+T_2[t]=Id_{\Sigma}
$$
%----------------------------------------------------------
and for every $f\in\sum$ the following inequality
%----------------------------------------------------------
$$
\|T_2[t](f)\|_{\LPM}+t\|T_1[t](f)\|_{\LOP}\le C(n,p)K(t;f:\vec{A}).
$$
%----------------------------------------------------------
holds. Here as before $\sum=\SUM$.
%----------------------------------------------------------
%@@@@@@@@@@@@@@@@@@@@@@@@@@@@@@@@@@@@@@@@@@@@@@@@@@@@@@@@@@
%@@@@@@@@@@@@@@@@@@@@@@@@@@@@@@@@@@@@@@@@@@@@@@@@@@@@@@@@@@
%@@@@@@@@@@@@@@@@@@@@@@@@@@@@@@@@@@@@@@@@@@@@@@@@@@@@@@@@@@
%@@@@@@@@@@@@@@@@@@@@@@@@@@@@@@@@@@@@@@@@@@@@@@@@@@@@@@@@@@
%----------------------------------------------------------
\medskip
%----------------------------------------------------------
\par This property easily follows from Theorem \reff{MainLinear} and equalities \rf{F-KFN} and \rf{S-KFN}. In fact, let us apply Theorem \reff{MainLinear} to the measure $\mu_t:=\mu/t^p$. By this theorem there exist continuous linear operators
%----------------------------------------------------------
$$
T_1[t]:\mathlarger{\Sigma}_t\to \LOP~~~~\text{and}~~~~T_2[t]:\mathlarger{\Sigma}_t\to L_p(\RN;\mu_t)
$$
%----------------------------------------------------------
such that
%----------------------------------------------------------
$$
T_1[t]+T_2[t]=Id_{\Sigma_t}
$$
%----------------------------------------------------------
and
%----------------------------------------------------------
\bel{K-ST}
\|T_1[t]\|_{\Sigma_t\to\LOP}+\|T_2[t]\|_{\Sigma_t\to\LPM}
\le C(n,p).
\ee
%----------------------------------------------------------
Since
%----------------------------------------------------------
$$
\|f\|_{L_p(\RN;\mu_t)}=\tfrac{1}{t}\,\|f\|_{\LPM}\,,
$$
%----------------------------------------------------------
the Banach space $L_p(\RN;\mu_t)$ coincides with the Banach space $\LPM$ proving that the Banach spaces $\Sigma_t$ and $\Sigma$ coincide as well. Hence
%----------------------------------------------------------
$$
T_1[t]:\ssmall\to \LOP~~~~\text{and}~~~~T_2[t]:\ssmall\to \LPM
$$
%----------------------------------------------------------
and
%----------------------------------------------------------
$$
T_1[t]+T_2[t]=Id_{\Sigma}.
$$
%----------------------------------------------------------
\par Furthermore, by \rf{F-KFN} and \rf{K-ST}, for every $f\in\sum$ we have
%----------------------------------------------------------
$$
\|T_1[t](f)\|_{\LOP}\le C\|f\|_{\ssmall_t}=\,CK(t;f:\vec{A})/t,
$$
%----------------------------------------------------------
and
%----------------------------------------------------------
$$
\|T_2[t](f)\|_{L_p(\RN;\mu_t)}=\|T_2[t](f)\|_{\LPM}/t
\le C\|f\|_{\Sigma_t}=\,CK(t;f:\vec{A})/t.
$$
%----------------------------------------------------------
We obtain
%----------------------------------------------------------
$$
t\|T_1[t](f)\|_{\LOP}\le C\,K(t;f:\vec{A}),
$$
%----------------------------------------------------------
and
%----------------------------------------------------------
$$
\|T_2[t](f)\|_{\LPM}\le C\,K(t;f:\vec{A}).
$$
%----------------------------------------------------------
\par Hence
%----------------------------------------------------------
$$
\|T_2[t](f)\|_{\LPM}+t\|T_1[t](f)\|_{\LOP}\le C(n,p)K(t;f:\vec{A})
$$
%----------------------------------------------------------
proving that the $K$-functional of the couple $\vec{A}=(\LPM,\LOP)$ is quasi-linearizable.\medskip
%----------------------------------------------------------
%@@@@@@@@@@@@@@@@@@@@@@@@@@@@@@@@@@@@@@@@@@@@@@@@@@@@@@@@@@
%@@@@@@@@@@@@@@@@@@@@@@@@@@@@@@@@@@@@@@@@@@@@@@@@@@@@@@@@@@
%@@@@@@@@@@@@@@@@@@@@@@@@@@@@@@@@@@@@@@@@@@@@@@@@@@@@@@@@@@
%@@@@@@@@@@@@@@@@@@@@@@@@@@@@@@@@@@@@@@@@@@@@@@@@@@@@@@@@@@
%----------------------------------------------------------
\par Finally, we remark that Z. Ditzian and V. Totik \cite{DT} have studied a number of variants of the $K$-functional for the Banach couple $\vec{B}=(L_p(\R),L^1_{p}(\R;\mu))$ where   $L^1_{p}(\R;\mu)$ is a homogeneous Sobolev space on $\R$ with respect to the measure $\mu$. This space is defined by the finiteness of the seminorm
%----------------------------------------------------------
$$
\|f\|_{L^1_p(\R;\mu)}:=
\left(\,\intl_{\R}\,|f'(x)|^p\,d\mu(x)\right)^{\frac1p}.
$$
%----------------------------------------------------------
\par At first sight, the couples $\vec{A}=(L_p(\R;\mu),L^1_p(\R))$ and $\vec{B}=(L_p(\R),L^1_{p}(\R;\mu))$ look very similar to each other. Nevertheless, in general, the $K$-functionals of these couples are very different from each other, and their calculations require different methods and ideas.\bigskip
%----------------------------------------------------------
%@@@@@@@@@@@@@@@@@@@@@@@@@@@@@@@@@@@@@@@@@@@@@@@@@@@@@@@@@@
%@@@@@@@@@@@@@@@@@@@@@@@@@@@@@@@@@@@@@@@@@@@@@@@@@@@@@@@@@@
%@@@@@@@@@@@@@@@@@@@@@@@@@@@@@@@@@@@@@@@@@@@@@@@@@@@@@@@@@@
%@@@@@@@@@@@@@@@@@@@@@@@@@@@@@@@@@@@@@@@@@@@@@@@@@@@@@@@@@@
%----------------------------------------------------------
\par {\bf 7.3. Theorem \reff{VRN} and subfamilies of ``minimal'' pairwise disjoint cubes.}
%----------------------------------------------------------
\addtocontents{toc}{~~~~7.3. Theorem \reff{VRN} and subfamilies of ``minimal'' pairwise disjoint cubes. \hfill \thepage\\\\\par}
%----------------------------------------------------------
The result of Theorem \reff{VRN} can be reformulated in a purely geometrical way. In fact, given a function $w:\RN\to(0,\infty)$ consider a family of cubes
%----------------------------------------------------------
$$
\Ac_w=\{Q=Q(x,w(x)):x\in\RN\}
$$
%----------------------------------------------------------
and a subfamily of $\Ac_w$
%----------------------------------------------------------
$$
\Bc=\{K=Q(x,w(x)):x\in S\}
$$
%----------------------------------------------------------
where $S$ is the set determined in Theorem \reff{VRN}.
Then the conditions $(i)$ and $(ii)$ from this proposition are equivalent to the following statements:\medskip
%----------------------------------------------------------
\par $(i')$. For every cube $Q\in\Ac_w$ there exists a cube $K\in\Bc$ such that $K\subset 83 Q$;\medskip
%@@@@@@@@@@@@@@@@@@@@@@@@@@@@@@@@@@@@@@@@@@@@@@@@@@@@@@@@@@
%----------------------------------------------------------
\par $(ii')$. The cubes of the family $\Bc$ are pairwise disjoint.\medskip
%----------------------------------------------------------
\par Thus Theorem \reff{VRN} states that for every function $w$ satisfying condition \rf{AZ} the family $\Ac_w$ contains a subfamily $\Bc$ satisfying conditions $(i')$ and $(ii')$.
%----------------------------------------------------------
\par This geometrical reformulation of the proposition motivates the following
%----------------------------------------------------------
%@@@@@@@@@@@@@@@@@@@@@@@@@@@@@@@@@@@@@@@@@@@@@@@@@@@@@@@@@@
%@@@@@@@@@@@@@@@@@@@@@@@@@@@@@@@@@@@@@@@@@@@@@@@@@@@@@@@@@@
%----------------------------------------------------------
\begin{question}\lbl{Q-I} Let $\Ac$ be a family of cubes in $\RN$. Under what conditions on $\Ac$ there exists a subfamily $\Bc$ of $\Ac$ such that:
%----------------------------------------------------------
\par (a). For every cube $Q\in\Ac$ there exists a cube $K\in\Bc$ such that $K\cap Q\ne\emp$ and $\diam K\le\diam Q$;
%@@@@@@@@@@@@@@@@@@@@@@@@@@@@@@@@@@@@@@@@@@@@@@@@@@@@@@@@@@
%----------------------------------------------------------
\par (b). The cubes of the family $\Bc$ are pairwise disjoint\,?
%----------------------------------------------------------
\end{question}
%----------------------------------------------------------
%@@@@@@@@@@@@@@@@@@@@@@@@@@@@@@@@@@@@@@@@@@@@@@@@@@@@@@@@@@
%@@@@@@@@@@@@@@@@@@@@@@@@@@@@@@@@@@@@@@@@@@@@@@@@@@@@@@@@@@
%@@@@@@@@@@@@@@@@@@@@@@@@@@@@@@@@@@@@@@@@@@@@@@@@@@@@@@@@@@
%@@@@@@@@@@@@@@@@@@@@@@@@@@@@@@@@@@@@@@@@@@@@@@@@@@@@@@@@@@
%@@@@@@@@@@@@@@@@@@@@@@@@@@@@@@@@@@@@@@@@@@@@@@@@@@@@@@@@@@
%@@@@@@@@@@@@@@@@@@@@@@@@@@@@@@@@@@@@@@@@@@@@@@@@@@@@@@@@@@
%----------------------------------------------------------
\begin{remark} {\em Note that the condition (a) implies the inclusion $2Q\supset K$.
%----------------------------------------------------------
\par  Also, let $\gamma\ge 1$ and let $\gamma\Ac=\{\gamma Q:Q\in\Ac\}$. Then the existence of a collection $\Bc$ satisfying the conditions (a) and (b) for the family $\gamma\Ac$  implies the existence of a subfamily $\widetilde{\Bc}$ of the family $\Ac$ such that: $(a').$ For each $Q\in\Ac$ there exists $K\in\widetilde{\Bc}$ such that $K\subset (2\gamma)Q$; $(b').$ the cubes
$\{\gamma K:K\in\widetilde{\Bc}\}$ are pairwise disjoint.
%----------------------------------------------------------
\par Clearly, one can put $\widetilde{\Bc}=\tfrac1\gamma \Bc$.\rbx}
\end{remark}
%----------------------------------------------------------
%@@@@@@@@@@@@@@@@@@@@@@@@@@@@@@@@@@@@@@@@@@@@@@@@@@@@@@@@@@
%@@@@@@@@@@@@@@@@@@@@@@@@@@@@@@@@@@@@@@@@@@@@@@@@@@@@@@@@@@
%@@@@@@@@@@@@@@@@@@@@@@@@@@@@@@@@@@@@@@@@@@@@@@@@@@@@@@@@@@
%@@@@@@@@@@@@@@@@@@@@@@@@@@@@@@@@@@@@@@@@@@@@@@@@@@@@@@@@@@
%@@@@@@@@@@@@@@@@@@@@@@@@@@@@@@@@@@@@@@@@@@@@@@@@@@@@@@@@@@
%@@@@@@@@@@@@@@@@@@@@@@@@@@@@@@@@@@@@@@@@@@@@@@@@@@@@@@@@@@
%----------------------------------------------------------
\par V. Dolnikov kindly drew the author's attention to the fact that a family $\Bc$ satisfying conditions (a) and (b) exists whenever $\Ac$ is an arbitrary {\it finite} collection of cubes. Here is a short Dolnikov's proof of this statement. %----------------------------------------------------------
\par Let $K_1$ be a cube of the minimal diameter among all the cubes of the family $\Ac_1:=\Ac$. By $G_1$ we denote all cubes of $\Ac_1$ which intersect $K_1$.
%----------------------------------------------------------
\par We put $\Ac_2:=\Ac_1\setminus G_1$. If $\Ac_2=\emp$ we stop and put $\Bc=\{K_1\}$.  If  $\Ac_2\ne\emp$, by $K_2$ we denote a cube of the minimal diameter among all the cubes of the family $\Ac_2$. We continue this procedure.  Since $\Ac$ is finite, this process will stop on a certain (finite) step $m$.
%----------------------------------------------------------
\par As a result we obtain a finite collection of pairwise disjoint cubes $\Bc=\{K_1,...,K_m\}$ and a partition $\{G_1,...,G_m\}$ of $\Ac$ such that for each $1\le i\le m$ the following conditions are satisfied: the cube $K_i\in G_i$,  $K_i\cap Q\ne\emp$, and $\diam K_i\le\diam Q$ for every $Q\in G_i$. Clearly, the collection $\Bc$ satisfies the conditions (a) and (b) of Question \reff{Q-I}.
%----------------------------------------------------------
\par Let us also note that for an {\it infinite} collection $\Ac$ of cubes in $\RN$ a family of cubes $\Bc$ satisfying conditions (a) and (b) of Question \reff{Q-I} in general does not exist. For instance, $\Bc$ does not exists whenever $\Ac=\{Q\left(0,\tfrac1n\right):n=1,2,...\}$.
%----------------------------------------------------------
\par These examples show that a certain ``continuity'' condition (apparently in the spirit of the condition \rf{AZ}) should be posed on the collection $\Ac$ to provide the existence of a subfamily $\Bc$ satisfying  conditions $(a)$ and $(b)$ of Question \reff{Q-I}.
%----------------------------------------------------------
%@@@@@@@@@@@@@@@@@@@@@@@@@@@@@@@@@@@@@@@@@@@@@@@@@@@@@@@@@@
%@@@@@@@@@@@@@@@@@@@@@@@@@@@@@@@@@@@@@@@@@@@@@@@@@@@@@@@@@@
%@@@@@@@@@@@@@@@@@@@@@@@@@@@@@@@@@@@@@@@@@@@@@@@@@@@@@@@@@@
%@@@@@@@@@@@@@@@@@@@@@@@@@@@@@@@@@@@@@@@@@@@@@@@@@@@@@@@@@@
%@@@@@@@@@@@@@@@@@@@@@@@@@@@@@@@@@@@@@@@@@@@@@@@@@@@@@@@@@@
%----------------------------------------------------------
%&&&&&&&&&&&&&&&&&&&&&&&&&&&&&&&&&&&&&&&&&&&&&&&&&&&&&&&&&&
%                                                         &
%                      REFERENCES                         &
%_________________________________________________________&
%&&&&&&&&&&&&&&&&&&&&&&&&&&&&&&&&&&&&&&&&&&&&&&&&&&&&&&&&&&

%@@@@@@@@@@@@@@@@@@@@@@@@@@@@@@@@@@@@@@@@@@@@@@@@@@@@@@@@@@
%@@@@@@@@@@@@@@@@@@@@@@@@@@@@@@@@@@@@@@@@@@@@@@@@@@@@@@@@@@
%@@@@@@@@@@@@@@@@@@@@@@@@@@@@@@@@@@@@@@@@@@@@@@@@@@@@@@@@@@

\begin{thebibliography}{AB}
%@@@@@@@@@@@@@@@@@@@@@@@@@@@@@@@@@@@@@@@@@@@@@@@@@@@@@@@@@@
\bibitem {BSh} C. Bennett, R. Sharpley, Interpolation of operators. Pure and Applied Mathematics, 129. Academic Press, Inc., Boston, MA, 1988. xiv+469 pp.
%----------------------------------------------------------
\bibitem {BL} J. Bergh, J. L\"{o}fstr\"{o}m, Interpolation Spaces, An Introduction, Springer-Verlag, 1976.
%----------------------------------------------------------
\bibitem {Be} A. S. Besicovitch, A general form of the covering principle and relative differentiation of
additive functions (I), (II), Proc. Cambridge Philos. Soc. {\bf 41} (1945), 103--110; {\bf 42} (1946) 1--10.
%----------------------------------------------------------
\bibitem {DT} Z. Ditzian, V. Totik, Moduli of smoothness. Springer Verlag, New York, 1987.
%----------------------------------------------------------
\bibitem {FIL} C. Fefferman, A. Israel, G. K. Luli,
Sobolev extension by linear operators, arXiv:1205.2525v2
%----------------------------------------------------------
\bibitem {G} M. de Guzm\'{a}n, Differentiation of
    integrals in $\RN$, Lect. Notes in Math. 481,
    Springer-Verlag, 1975.
%----------------------------------------------------------
\bibitem {Is} A. Israel, A Bounded Linear Extension Operator for $L^{2,p}(\RT)$, arXiv: 1011.0689v1.
%----------------------------------------------------------
\bibitem {JT}  T. R. Jensen and B. Toft, Graph Coloring Problems, Wiley-Interscience, New York, 1995.
%----------------------------------------------------------
\bibitem {M}  V.G. Maz'ja,  Sobolev spaces,
     Springer-Verlag, Berlin, 1985, xix+486 pp.
%----------------------------------------------------------
\bibitem {MP} V. Maz'ya,  S. Poborchi, Differentiable
    Functions on Bad Domains, Word Scientific, River Edge,
    NJ, 1997.
%----------------------------------------------------------
\bibitem {P} J. Peetre, A theory of interpolation of normed spaces. Notas de Matematica, No. 39, Instituto de Matematica Pura e Aplicada, Conselho Nacional de Pesquisas, Rio de Janeiro, 1968, iii+86 pp.
%-----------------------------------------------------------
\bibitem {S6} P. Shvartsman, Sobolev $L^2_p$-functions on closed subsets of $\R^2$ (preprint).
%----------------------------------------------------------
\bibitem {St} E. M. Stein, Singular integrals and
    differentiability properties of functions, Princeton
    Univ. Press, Princeton, New Jersey, 1970.
%----------------------------------------------------------
%@@@@@@@@@@@@@@@@@@@@@@@@@@@@@@@@@@@@@@@@@@@@@@@@@@@@@@@@@@
\end{thebibliography}
\end{document}